\documentclass{amsart}
\title{The R-matrix of the affine Yangian}

\author[A. Appel]{Andrea Appel} 
\address{Dipartimento di Scienze Matematiche, 
	Fisiche e Informatiche, Universit\`a di Parma, 
	INdAM-GNSAGA and INFN Gruppo Collegato di Parma, 
	43124 Parma, Italy}
\email{\href{mailto:andrea.appel@unipr.it}{andrea.appel@unipr.it}}

\author[S. Gautam]{Sachin Gautam}
\address{Department of Mathematics, The Ohio State University,
Columbus, OH 43210, USA}
\email{\href{mailto:gautam.42@osu.edu}{gautam.42@osu.edu}}

\author[C. Wendlandt]{Curtis Wendlandt}
\address{Department of Mathematics and Statistics,
University of Saskatchewan,
Saskatoon, SK S7N 5E6,
Canada}
\email{\href{mailto:wendlandt@math.usask.ca}{wendlandt@math.usask.ca}}


\usepackage{amsmath,amssymb,mathtools}
\usepackage[usenames,dvipsnames]{color}
\usepackage{enumitem}
\usepackage{comment}
\usepackage{upgreek}
\usepackage{hyperref} 
\hypersetup{
    colorlinks=true,       			
    linkcolor=blue, 
    citecolor=blue,        			
    filecolor=blue,      				
    urlcolor=blue           			
}

\usepackage{float}
\usepackage[all]{xy}
\usepackage{tikz-cd} 
\xyoption{arc}
\usepackage{dynkin-diagrams}
\newtheorem*{thm}{Theorem}
\newtheorem*{prop}{Proposition}
\newtheorem*{lem}{Lemma}
\newtheorem*{cor}{Corollary}

\newenvironment{pf}[1][Proof]{
	\begin{proof}[\scshape{#1}]}
	{\end{proof}
	}
\theoremstyle{definition}

\newtheorem*{rem}{Remark}

\numberwithin{equation}{section}

\numberwithin{figure}{section}

\newcommand{\Pseries}[2]{#1[\negthinspace[#2]\negthinspace]}

\newcommand{\Triv}{\mathsf{1}}

\newcommand{\Top}{\operatorname{T}}
\newcommand{\adT}{\mathcal{T}}

\newcommand{\id}{\mathbf{1}}

\newcommand{\veps}{\varepsilon}
\newcommand{\lp}{\left(}
\newcommand{\rp}{\right)}

\newcommand{\g}{\mathfrak{g}}
\newcommand{\h}{\mathfrak{h}}

\newcommand{\Lsl}{\mathfrak{sl}}

\newcommand{\Sym}{\mathfrak{S}}

\newcommand{\vup}[1]{v_+}
\newcommand{\vdown}[1]{v_-}

\newcommand{\bfA}{\mathbf{A}}
\newcommand{\bfB}{\mathbf{B}}

\newcommand{\calO}{\mathcal{O}}
\newcommand{\PP}{\mathcal{P}}

\newcommand{\RR}{\mathcal{R}}

\newcommand{\C}{\mathbb{C}}
\newcommand{\nC}{\mathbb{C}^{\times}}
\newcommand{\N}{\mathbb{Z}_{\geqslant 0}}

\newcommand{\Q}{\mathbb{Q}}
\newcommand{\R}{\mathbb{R}}
\newcommand{\Z}{\mathbb{Z}}

\renewcommand{\Re}{\operatorname{Re}}

\newcommand{\hit}{\operatorname{ht}}

\newcommand {\wh}[1]{\widehat{#1}}
\newcommand {\ol}[1]{\overline{#1}}
\newcommand {\ul}[1]{\underline{#1}}

\newcommand{\reg}{\scriptscriptstyle{\operatorname{reg}}}

\newcommand {\Aut}{\operatorname{Aut}}

\newcommand{\End}{\operatorname{End}}
\newcommand{\Hom}{\operatorname{Hom}}
\newcommand{\ad}{\operatorname{ad}}
\newcommand{\Ad}{\operatorname{Ad}}

\renewcommand{\dim}{\operatorname{dim}}

\newcommand{\op}{\operatorname{op}}

\newcommand {\Rep}{\operatorname{Rep}}

\newcommand{\ds}{\displaystyle}

\newcommand {\Comment}[1]{}
\newcommand {\Omit}[1]{}

\newcommand {\Yhg}{Y_\hbar(\g)}

\newcommand{\Yhsl}[1]{Y_{\hbar}(\sl_{#1})}

\newcommand {\longisom}{\stackrel{\sim}{\longrightarrow}}

\newcommand{\cbin}[2]{\left(\begin{array}{c} #1\\ #2\end{array}\right)}

\renewcommand {\sl}{\mathfrak{sl}}

\newcommand {\sfA}{{\mathsf A}}

\newcommand {\sfD}{{\mathsf D}}

\newcommand {\Id}{\operatorname{Id}}
\newcommand {\ie}{{\it i.e.}, }
\newcommand{\bfI}{{\mathbf I}}

\newcommand{\Borel}[1]{\mathcal{B}\lp #1\rp}
\newcommand{\Laplace}[1]{\mathcal{L}\lp #1\rp}

\newcommand{\hbfI}{\bfI}
\newcommand{\hbfA}{\bfA}
\newcommand{\gkm}{\g}
\newcommand{\gaff}{\g^\prime}
\newcommand{\hkm}{\h}
\newcommand{\haff}{\h^\prime}
\newcommand{\naff}{\mathfrak{n}}

\newcommand{\caff}{\EuScript{C}} 
\newcommand{\sczero}{\ell^{(0)}} 
\newcommand{\qdzero}{\mathsf{d}^{(0)}} 
\newcommand{\ctzero}{\mathsf{g}_0} 
\newcommand{\ctsing}{\mathsf{g}_{\scriptscriptstyle{\mathrm{sing}}}}

\newcommand{\anode}{\circ} 

\newcommand{\Qaff}{\mathsf{Q}}
\newcommand{\raff}{\Phi}

\newcommand{\Ykm}{Y_{\hbar}(\gkm)}
\newcommand{\Ysupkm}[1]{Y^{#1}_{\hbar}(\gkm)}

\newcommand{\OY}{\mathcal{O}(\Ykm)}

\newcommand{\Prim}{\mathrm{Prim}^{\scriptscriptstyle{\mathrm{D}}}} 

\newcommand{\ddelta}[1]{\underset{\scriptscriptstyle{{\operatorname{D}},#1}}{\Delta}}

\newcommand{\ddeltaz}[1]{\underset{\scriptscriptstyle{{\operatorname{D}},#1\,\,\,}}{\Delta^z}}
\newcommand{\kmdeltaz}[1]{\underset{\scriptscriptstyle{{\operatorname{KM}},#1\,\,\,}}{\Delta^z}}
\newcommand{\dtensor}[1]{\underset{\scriptscriptstyle{\operatorname{D},#1}}{\otimes}}
\newcommand{\kmtensor}[1]{\underset{\scriptscriptstyle{\operatorname{ KM},#1}}{\otimes}}

\newcommand{\sop}{\mathsf{p}} 
\newcommand{\Rd}{\mathcal{D}} 
\newcommand{\Rtau}{\mathcal{T}} 
\newcommand{\Rg}{\mathcal{G}} 
\newcommand{\Rl}{\Lambda} 

\DeclareMathAlphabet\EuScript{U}{eus}{m}{n}

\SetMathAlphabet\EuScript{bold}{U}{eus}{b}{n}

\newcommand {\curtiscomment}[1]{} 

\newcommand{\Thz}[3]{\Theta_{#2}(#1)_{#3}}

\newcommand{\dri}{\scriptscriptstyle{\mathrm{D}}}
\newcommand{\gfin}{\mathfrak{a}}
\newcommand{\Yhgfin}{Y_\hbar(\gfin)}
\newcommand{\gKM}{\g_{\scriptstyle{\mathsf{KM}}}}
\newcommand{\YhgKM}{Y_\hbar(\gKM)}
\newcommand{\fml}[1]{[\![#1]\!]}
\newcommand{\ten}{\otimes}
\newcommand{\Rmsb}[2]{\RR^-(#2)_{#1}} 
\newcommand{\hgt}[1]{\mathsf{ht}(#1)}
\newcommand{\Topp}{\Top'}
\newcommand{\sfR}{\mathsf{R}}

\allowdisplaybreaks

\newcommand{\ceq}{=} 

\begin{document}

\maketitle

\begin{abstract}
Let $\g$ be an affine Lie algebra with associated Yangian $\Yhg$.
We prove the existence of two meromorphic $R$--matrices associated to any pair of representations of $\Yhg$ in the category $\mathcal{O}$. 
They are related by a unitary constraint and constructed as products of the form $\RR^{\uparrow/\downarrow}(s)=\RR^+(s)\cdot\RR^{0,\uparrow/\downarrow}(s)\cdot\RR^-(s)$, where $\RR^+(s) = \RR^-_{21}(-s)^{-1}$. 
The factors $\RR^{0,\uparrow/\downarrow}(s)$ are meromorphic, abelian $R$--matrices,
and $\RR^-(s)$ is a rational twist. 
Our proof relies on two novel ingredients.
The first is an irregular, abelian, additive difference equation
whose difference operator is given in terms of 
the $q$--Cartan matrix of $\g$.
The regularization of this difference equation gives rise to 
$\RR^{0,\uparrow/\downarrow}(s)$ as the
exponentials of the two canonical fundamental solutions.
The second key ingredient is
a higher order analogue of the adjoint action of 
the affine Cartan subalgebra $\h\subset\g$ on $\Yhg$. This action has no classical counterpart, and produces
a system of linear equations from which $\RR^-(s)$
is recovered as the unique solution. 
Moreover, we show that both $\RR^{\uparrow/\downarrow}(s)$
give rise to the same rational $R$--matrix 
on the tensor product of any two highest--weight representations.
\end{abstract}

{\setlength{\parskip}{0pt}
\setcounter{tocdepth}{1} 
\tableofcontents
}

\section{Introduction}\label{sec:Intro}


\subsection{}\label{opening}
Let $\g$ be an affine Kac-Moody algebra. In contrast to its
counterpart for a finite--dimensional, semisimple Lie algebra,
the affine Yangian $\Yhg$ is not known to have a universal $R$--matrix.
In particular, it is not known whether an arbitrary representation
$V$ of $\Yhg$ gives rise to a solution of the quantum Yang-Baxter
equation \eqref{qybe} with additive spectral parameter.

A notable exception arises when $\g$ is simply--laced, and $V$ is the
equivariant cohomology of a quiver
variety for the underlying affine Dynkin diagram.
In this setting, a rational solution of the \eqref{qybe} corresponding to $V$ has been constructed by
Maulik and Okounkov in \cite{maulik-okounkov-qgqc} using
stable envelopes.

The main goal of this paper is to construct solutions $\RR(s)$ of the
\eqref{qybe} for an arbitrary affine Yangian $\Yhg$, with the exception\footnote{
	The case $\g=\wh{\Lsl}_2$ needs to be addressed separately.
	The first (incomplete) definition of the Yangian of type $\mathsf{A}_1^{(1)}$ appeared in \cite{boyarchenko-levendorskii}, and was later modified in \cite{Tsy17b,Tsy17}, see also \cite{BerTsy19} and  \cite[Rem.~5.2]{Kod19}.}
of type $\mathsf{A}_1^{(1)}$,
 and a representation $V$ whose restriction
to $\g$ lies in its category $\mathcal{O}$. Our solutions are {\it meromorphic}
with respect to the spectral parameter $s$, natural with respect to
the underlying representations, and compatible with the tensor product. 

We also show
that, on highest weight representations, 
they can
be normalized so as to be rational in $s$, and conjecture that they then
coincide with those constructed in \cite{maulik-okounkov-qgqc}
when $\g$ is simply-laced and $V$ arises from geometry.

\subsection{}\label{opening-2}
Recall that the quantum Yang--Baxter equation is the following
equation for an $\End(V^{\otimes 2})$--valued function (or formal
power series) $\EuScript{R}(s)$:
\begin{equation}\label{qybe}
\tag{QYBE}
\EuScript{R}_{12}(s_1)\EuScript{R}_{13}(s_1+s_2)\EuScript{R}_{23}(s_2)
=
\EuScript{R}_{23}(s_2)\EuScript{R}_{13}(s_1+s_2)\EuScript{R}_{12}(s_1).
\end{equation}
This equation is
$\End(V^{\otimes 3})$--valued, and the subscripts
indicate on which tensor factors $\EuScript{R}$ acts.

Drinfeld's theory of Yangians 
provides a uniform method for
constructing rational solutions of this equation. 
The Yangian $\Yhgfin$ of a finite-dimensional simple Lie algebra 
$\gfin$ was introduced in \cite{drinfeld-qybe}. It is the canonical 
Hopf algebra deformation of the current algebra $U(\gfin[z])$, 
and gives rise to rational solutions of \eqref{qybe} on its irreducible, 
finite--dimensional representations.

To state Drinfeld's results more
precisely, let $\Delta:\Yhgfin\to\Yhgfin^{\otimes 2}$ be the
coproduct and $\tau_s\in\Aut(\Yhgfin)$ the
one--parameter group of Hopf algebra automorphisms,
quantizing the shift automorphisms $z\mapsto z+s$ of $U(\gfin[z])$.
We abbreviate $\Delta_s = \tau_s\otimes\Id\circ\Delta$
and $\Delta^{\op}_s = \tau_s\otimes\Id\circ\Delta^{\op}$.
Then, by \cite[Thm.~3]{drinfeld-qybe}, there exists a unique formal series
\[
\RR^{\dri}(s) \in 1^{\otimes 2} + s^{-1}\Pseries{\Yhgfin^{\otimes 2}}{s^{-1}}
\]
satisfying the intertwining equation
\[
\Delta^{\op}_s(x) = \RR^{\dri}(s)\Delta_s(x)\RR^{\dri}(s)^{-1} \quad \forall\quad x\in \Yhgfin,
\]
 in addition to the cabling identities
\[
\begin{aligned}
\Delta\otimes\Id (\RR^{\dri}(s)) &= \RR^{\dri}_{13}(s)\RR^{\dri}_{23}(s),\\
\Id\otimes\Delta (\RR^{\dri}(s)) &= \RR^{\dri}_{13}(s)\RR^{\dri}_{12}(s).
\end{aligned}
\]
%
These readily imply that $\RR^{\dri}(s)$ is a solution
to \eqref{qybe} on $\Yhgfin^{\otimes 3}$.

Let $V_1,V_2$ be two finite--dimensional irreducible representations
of $\Yhgfin$, and $\RR^{\dri}_{V_1,V_2}(s)$ the evaluation
of $\RR^{\dri}(s)$ on $V_1\otimes V_2$.
Upon normalizing
by its eigenvalue on the tensor product
of highest--weight spaces, the dependence of $\RR^{\dri}_{V_1,V_2}(s)$  on $s$ becomes rational
\cite[Thm.~4]{drinfeld-qybe}.
Hence, in particular, one obtains rational solutions to
\eqref{qybe} valued in any finite--dimensional, irreducible representation
$V$ of $\Yhgfin$.

The functional nature in $s\in\C$ of $\RR^{\dri}_{V_1,V_2}(s)$, for arbitrary
finite--dimensional representations $V_1,V_2$ of $\Yhgfin$, was studied
by the last two authors and Toledano Laredo in \cite{sachin-valerio-III,GTLW}.
It is shown to be a divergent
series in $s$ which admits two canonical resummations. In more succinct
terms, the category
$\Rep_{\scriptscriptstyle{\text{fd}}}(\Yhgfin)$ admits two {\em meromorphic} braidings  (see \cite[Thm.~7.1]{GTLW}).
It was also shown that there is no {\em rational}
braiding, \ie the rational normalization described above cannot be carried out 
uniformly for all representations (see Thm.~7.3 of {\em loc. cit.}).
These results do not rely on the
existence of $\RR^{\dri}(s)\in
\Pseries{\Yhgfin^{\otimes 2}}{s^{-1}}$ and, in fact, can be used
to prove it in a new and
constructive manner{\footnote{Drinfeld's proof was based on cohomological methods and has not yet appeared in the literature.}}
(see Thm. 7.4 of {\em loc. cit.}).

\subsection{}\label{YKM}
In \cite{drinfeld-yangian-qaffine}, a different
presentation of $\Yhgfin$ is given, called Drinfeld's new (or loop)
presentation, which takes as input only the Cartan matrix of $\gfin$. It can be used, more generally, to define the Yangian $\YhgKM$ of any
symmetrizable Kac--Moody algebra $\gKM$, but only as an algebra. In light of Drinfeld's foundational results, it is natural
to ask whether \eqref{qybe} can be solved for a suitable
category of representations of $\YhgKM$ in a functorial manner.
The following two obstacles present themselves almost immediately.

\begin{enumerate}[label=(O\arabic*)]\setlength{\itemsep}{3pt}

\item\label{it:O2} The universal $R$--matrix $\RR^{\dri}(s)$
is only known to exist in finite types. Moreover, our
results suggest, {\em a posteriori}, that perhaps there is no such object outside of finite types
(see Section \ref{what-is-left} \ref{intro-f1}--\ref{intro-f2} below).

\item\label{it:O1} The coproduct $\Delta$ has only been
defined in the case where $\gKM$ is of finite or affine type \cite{guay-nakajima-wendlandt}.\footnote{
The results from \cite{guay-nakajima-wendlandt} require to exclude types 
$\mathsf{A}_1^{(1)}$ and $\mathsf{A}_2^{(2)}$. The latter case, however, has been considered in \cite{Ueda-A2}, and we expect our results to readily extend to type $\mathsf{A}_2^{(2)}$.}

\end{enumerate}

\subsection{}\label{we-do}
We answer in the affirmative the question raised in the previous paragraph
for the Yangian $\Yhg$ associated to an affine Lie algebra $\g$.
We consider the category $\OY$, consisting of the representations of $\Yhg$
whose restriction to $\g$ lies in category $\mathcal{O}$.
We prove that there are two meromorphic braidings on $\OY$, related to each other
by a unitarity relation. On any tensor product of highest--weight representations
they both give rise to the same rational solution of \eqref{qybe}  (see Theorems \ref{thm:main} and \ref{thm:rat} below).

\subsection{}
The techniques employed in the aforementioned works of \cite{sachin-valerio-III, GTLW}
appear to be constrained by the finite--dimensionality of
the underlying Lie algebra
and the representations of the associated Yangian. 
In spite of this, they have
led us to formulate a very general principle, which we 
refer to as the {\em abelianization method}
(see Section \ref{ssec:abel} for a detailed description).
Roughly, this amounts to building $\RR(s)$ in the form
\begin{equation}\label{gauss}
	\RR(s) = \RR^+(s)\cdot \RR^0(s)\cdot \RR^-(s)\ ,
\end{equation}
where $\RR^0(s)$ is a meromorphic braiding with respect to 
the \emph{Drinfeld tensor product} (see Section~\ref{ssec:dr-tensor}), 
$\RR^-(s)$ is a rational twist relating the standard and 
the Drinfeld tensor products, and $\RR^+(s) = \RR^-_{21}(-s)^{-1}$.
Such factorizations, in the context of {\em Yangian doubles},
first appeared as a conjecture
in the work of Khoroshkin and Tolstoy \cite{khoroshkin-tolstoy}, 
later proved by the third author in \cite{WDYhg, WRQD}. 

Thus, we implement the abelianization method in the affine
setting which we believe to be amenable to further generalizations
(to other Kac--Moody types, Maulik--Okounkov Yangians,
their trigonometric and elliptic counterparts etc.)
and is the main contribution of this paper.

\subsection{}\label{what-is-left}
We end this introduction by highlighting three significant
features of our $R$--matrices that arise beyond finite--type concerning their $\hbar$--dependence, universality, and uniqueness.

\begin{enumerate}[leftmargin=2.5em, label=(F\arabic*)]\setlength{\itemsep}{0.25cm}
\item\label{intro-f1} 
Our $R$--matrices do not have a limit as $\hbar\to 0$. Indeed, the
reader can notice the appearance of $\hbar$ in the denominators
of the formulae given in Theorems
\ref{T:map-T} and \ref{thm:R0-main}, which therefore cannot be specialized at $\hbar=0$.
This is consistent with the results of Enriquez
\cite[Rmk. 8]{Enr03} and Yang--Zhao \cite[Sec.~2.1.1]{YaGuPBW} at the classical level. The semiclassical analysis of our work will be carried out in a forthcoming joint work by the third author and Alex Weekes.


In contrast, the $R$--matrices obtained by Maulik and Okounkov \cite{maulik-okounkov-qgqc}
are rational in both deformation and spectral parameters, and tend to
 the identity as $\hbar\to 0$ or $s\to \infty$. This is not in contradiction with
 the conjecture at the end of Section~\ref{opening}, since
the singularities of our $R$--matrices at $\hbar=0$ disappear {\em on highest--weight representations}.
 The precise formulation of this fact is out of the scope of this paper and
 it will be presented in a separate work.

\item\label{intro-f2}
Our $R$--matrices do not correspond to the action of a {\em universal}
element in $\Pseries{\Yhg^{\otimes 2}}{s^{-1}}$. At the core, the problem
is that the Casimir tensor, outside of finite types, 
is not an element of $\g\otimes \g$. It belongs instead to a suitable completion. 
This is also the reason for the appearance of a $z$--parameter in the formulae of \cite{guay-nakajima-wendlandt}. Even with this modification,
the unipotent parts $\RR^{\pm}(s,z)$ of our $R$--matrices can only be viewed as
elements of $\Pseries{\Pseries{\Yhg^{\otimes 2}}{s^{-1}}}{z^{\mp 1}}$
which cannot be multiplied together as such. In this paper, we
only work with category $\mathcal{O}$ representations, which circumvents this
problem.

\item\label{intro-f3} 
Our $R$--matrices are not unique. Their construction is based on the intertwining
and cabling equations, hence they are only unique up to central, primitive elements. 
In finite--type, one can show that there are no such elements in $\Yhgfin$
and the uniqueness follows (see \cite[App.~B]{GTLW}). In our case,
we show that there are at least three such elements, given in Section \ref{ssec:cr}.
We conjecture that uniqueness holds up to multiplication
by $\exp(Xs^{-1})$ where $X$ lies in the tensor--square of their linear span.
\end{enumerate}

\section{Main results}\label{sec:main-results}

This section contains the main theorems of this paper and a detailed explanation of the abelianization method.

\subsection{}
Let $\g$ be a Kac--Moody algebra of affine type, and let $\Yhg$ be the
Yangian associated to $\g$, where $\hbar\in \C^\times$ is a fixed nonzero complex number (see Section \ref{ssec: yangian}). Let $\OY$ be the category consisting of those $\Yhg$-modules  whose restriction to $\g\subset\Yhg$ lies in its category $\calO$ (see Section \ref{ssec:cat-O}).
In \cite{guay-nakajima-wendlandt}, Guay, Nakajima, and the third
author define an algebra homomorphism
\begin{equation*}
\Delta^z:\Yhg\to\Yhg^{\ten 2}[z^{-1};z]\!]
\end{equation*}
called the {\em twisted coproduct} (see Section \ref{ssec:pre-cop}), where $\Yhg^{\ten 2}[z^{-1};z]\!]$ is the space of Laurent series in $z$ with coefficients in $\Yhg^{\ten 2}$.
%
Evaluating $\Delta^z$
at $z=1$ leads to
terms with infinitely many tensors, which is to be interpreted
as a topological coproduct.
Nevertheless, this yields a
well--defined action of $\Yhg$ 
on any tensor product of modules in the category $\OY$. We denote the tensor product
of $V_1,V_2\in\OY$ by $V_1\kmtensor{0} V_2$ and, for any $s\in \C$, we set 
\[
V_1\kmtensor{s} V_2 \ceq V_1(s)\kmtensor{0} V_2\,
\]
where $V_1(s)\ceq\tau_s^*(V_1)$ with $\{\tau_s\}_{s\in \C}\subset \Aut(\Yhg)$ the one--parameter group of bialgebra automorphisms defined in Section \ref{ssec: shift-yangian}. 

\subsection{Meromorphic $R$--matrices}\label{ssec:thm1}
We are now in position to state the first main theorem of this paper.

\begin{thm}\label{thm:main}
Let $V_1,V_2\in\OY$ 
Then, there are two meromorphic
functions $\RR^\eta_{V_1,V_2}\colon\C\to\End(V_1\otimes V_2)$,
$\eta\in\{\,\uparrow\,,\, \downarrow\,\}$,
with the following properties.

\begin{enumerate}[font=\upshape]\itemsep0.25cm
\item\label{eq:intro-R-fun} $\RR^{\eta}_{V_1,V_2}(s)$ is holomorphic on
a half--plane
$\varepsilon(\eta)\Re(s/\hbar)\gg 0$ and approaches $\Id_{V_1\otimes V_2}$
as $\varepsilon(\eta)\Re(s/\hbar)\to\infty$. Here, $\varepsilon(\uparrow)=+1$
and $\varepsilon(\downarrow)=-1$.

\item\label{eq:intro-R-inter}
The following is a $\Yhg$--intertwiner which is natural in $V_1$ and $V_2$:
\begin{equation*}\label{eq:intro-mero-braid}
(1\,2)\circ\RR^\eta_{V_1,V_2}(s)\,\colon\,V_1\kmtensor{s} 
V_2\to (V_2\kmtensor{-s} V_1)(s)\,
\end{equation*}

\item\label{eq:intro-R-QYBE} For $V_1,V_2,V_3\in\OY$, the following
cabling identities hold.
\[        
\begin{aligned}
\RR_{V_1\kmtensor{s_1} V_2, V_3}^{\eta}(s_2) &= 
\RR_{V_1,V_3}^{\eta}(s_1+s_2)\cdot\RR_{V_2,V_3}^{\eta}(s_2)\,,\\
\RR_{V_1,V_2\kmtensor{s_2}V_3}(s_1+s_2) &=
\RR_{V_1,V_3}^{\eta}(s_1+s_2)\cdot\RR_{V_1,V_2}^{\eta}(s_1)\,.
\end{aligned}
\]
In particular, the quantum Yang--Baxter equation holds on $V_1\otimes V_2\otimes V_3$\,:
\begin{equation*}
\begin{aligned}
\RR_{V_1,V_2}^{\eta}(s_1)\cdot \RR_{V_1,V_3}^{\eta}(s_1+s_2)&\cdot\RR_{V_2,V_3}^{\eta}(s_2)\\
=\RR_{V_2,V_3}^{\eta}(s_2)&\cdot \RR_{V_1,V_3}^{\eta}(s_1+s_2)\cdot\RR_{V_1,V_2}^{\eta}(s_1)\,.
\end{aligned}
\end{equation*}

\item\label{eq:intro-R-translation} For any $a,b\in\C$ and $V_1,V_2\in\OY$, we have
\begin{equation*}
\RR_{V_1(a),V_2(b)}^{\eta}(s)=\RR_{V_1,V_2}^{\eta}(s+a-b)\,.
\end{equation*}

\item\label{eq:intro-R-unit} The following unitarity relation holds
\begin{equation*}
\RR_{V_1,V_2}^{\uparrow}(s)^{-1}=(1\,2)\circ\RR_{V_2,V_1}^{\downarrow}(-s)\circ(1\,2)\,.
\end{equation*}

\item\label{eq:intro-R-asym} $\RR^\eta_{V_1,V_2}(s)$ have
the same asymptotic expansion\footnote{
In this paper we use the classical notion of asymptotic
expansions \`{a} la Poincar\'{e}.
Namely, for an unbounded set $S\subset\C$, we say $f(z)\sim
\sum_{r=0}^{\infty} f_r z^{-r}$, as $z\to \infty$ in $S$,
if for every $n\geq 0$, we have
\[
\lim_{\begin{subarray}{c} z\to\infty\\ z\in S\end{subarray}}
z^n\lp f(z)-\sum_{r=0}^n f_rz^{-r}\rp = 0.
\]
}
as $\Re(s/\hbar)\to\varepsilon(\eta)\infty$.
This expansion remains
valid in a larger sector $\Sigma^\eta_\delta$, for
any $\delta>0$, where, if $\theta=\arg(\hbar)$, then
\[
\Sigma^\uparrow_\delta\ceq\{re^{\iota\phi}:r\in\R_{>0}\ \text{ and }
\phi\in (\theta-\pi+\delta,\theta+\pi-\delta)\} = -\Sigma^\downarrow_\delta\ .
\]
(see Figure \ref{afig:sector} given in Appendix \ref{app:LapDE}).

\end{enumerate}
\end{thm}

More concisely, the theorem is equivalent to
the existence of two meromorphic braidings,
related to each other by the unitary relation \eqref{eq:intro-R-unit},
on the category $\OY$ with respect to the
meromorphic (in fact, polynomial) tensor product
$\kmtensor{s}$.

\subsection{Rational $R$--matrices}\label{ssec:rat}
Our results also imply the existence of a {\em rational} $R$--matrix
on any tensor product of highest weight modules, denoted below by $\sfR(s)$.  Namely,
we prove the following:

\begin{thm}\label{thm:rat}
Let $V_1,V_2\in\OY$ be representations which are generated by highest--weight
vectors $\mathsf{v}_1,\mathsf{v}_2$, respectively.
Then, $\RR^{\uparrow/\downarrow}_{V_1,V_2}(s)$
yield, after normalizing to take value $1$ on $\mathsf{v}_1\ten\mathsf{v}_2$,
the same operator $\sfR_{V_1,V_2}(s)$, whose matrix entries are rational in $s$.
Moreover, $\ds 
\sfR_{V_1,V_2}(\infty) =
\Id_{V_1\otimes V_2}$.
\end{thm}

The theorem above is more general
than
\cite[Thm.~4]{drinfeld-qybe} for finite types. The latter
assumes that the representations are finite--dimensional
and irreducible, and its proof is
presumably\footnote{Drinfeld's original proof has not appeared in print. A proof based on generic irreducibility can be found in \cite[Thm.~3.10]{GRW-equivalences}.}
based on a ``generic irreducibility" argument.
Our proof of Theorem \ref{thm:rat}
(see Section \ref{ssec:rat-R0})
relies on the
explicit calculation of $\Ad(\RR^0(s))$ acting on
operators from $\Yhg\otimes\Yhg$.
It is valid for highest--weight representations 
and assumes neither integrability nor generic cyclicity.

\subsection{Abelianization method}\label{ssec:abel}
The strategy employed in \cite{GTLW}
is based on an interplay between the standard
tensor product and the \emph{Drinfeld tensor product}.
The latter defines a meromorphic tensor structure on $\OY$ which depends \emph{rationally} on $s$.
Namely, it gives rise to a family of actions of $\Yhg$ on $V_1\ten V_2$,
which is a rational function in $s$, denoted by  $V_1\dtensor{s} V_2$
(see Theorem \ref{thm:GTLW}).
The abelianization method 
decouples
the 
construction of $\RR(s)$
into two separate, independent problems. Assume the following datum
is given:

\begin{enumerate}[label=(M\arabic*)]\itemsep0.25cm
\item\label{prob:R0}  two meromorphic braidings
with respect to the Drinfeld tensor product on $\OY$, \ie
a natural system of invertible intertwiners
\begin{equation*}\label{eq:R0-intro}
(1\,2)\circ\RR_{V_1,V_2}^{0,\uparrow/\downarrow}(s): V_1\dtensor{s}V_2\to 
(V_2\dtensor{-s}V_1)(s)\,,
\end{equation*}
for any $V_1,V_2\in\OY$, satisfying the analogue of
Theorems~\ref{thm:main} and \ref{thm:rat} with respect to the Drinfeld tensor product;

\item\label{prob:Rm} a meromorphic twist intertwining
the standard tensor product and the Drinfeld tensor product on $\OY$,
\ie a natural system of invertible intertwiners
\begin{equation*}\label{eq:Rm-intro}
\RR^-_{V_1,V_2}(s): V_1\kmtensor{s}V_2\to V_1\dtensor{s}V_2\,,
\end{equation*}
for any $V_1,V_2\in\OY$, such that the cocycle equation
\begin{equation*}\label{eq:cocycle-intro}
\RR^-_{V_1\dtensor{s_1}V_2, V_3}(s_2)\cdot \RR^-_{V_1, V_2}(s_1) 
= \RR^-_{V_1, V_2\dtensor{s_2}V_3}(s_1+s_2)\cdot \RR^-_{V_2,V_3}(s_2)
\end{equation*}
holds for any $V_1,V_2,V_3\in\OY$. 
We further assume that the matrix entries of
$\RR^-(s)$ depend rationally on $s$, and $\RR^-(\infty) = \Id$.
These conditions
ensure that the functional nature of $\RR^0(s)$ does not change upon
multiplication by $\RR^-(s)$ on the right and by $\RR^-_{21}(-s)^{-1}$
on the left.
\end{enumerate}

Then, we consider the isomorphism of $\Yhg$--modules
\begin{equation*}\label{eq:R-intro}
(1\,2)\circ \RR_{V_1,V_2}^{\uparrow/\downarrow}(s): 
V_1\kmtensor{s}V_2\to (V_2\kmtensor{-s}V_1)(s)\,,
\end{equation*}
defined by the commutativity of the following diagram:
\begin{equation}\label{eq:intro-Gauss-diag}
\begin{tikzcd}
V_1\kmtensor{s}V_2 \arrow[d, "\RR^-_{V_1,V_2}(s)"']
\arrow[rrr, "(1\,2)\circ\RR_{V_1,V_2}^{\uparrow/\downarrow}(s)"]
& & & (V_2\kmtensor{-s}V_1)(s)\arrow[d, "\RR^-_{V_2,V_1}(-s)"] \\
V_1\dtensor{s}V_2 \arrow[rrr, "(1\,2)\circ\RR_{V_1,V_2}^{0,\uparrow/\downarrow}(s)"']&&&
(V_2\dtensor{-s}V_1)(s)
\end{tikzcd}
\end{equation}
By \cite[\S7.1]{GTLW}, this is readily seen to satisfy the properties
\ref{eq:intro-R-fun}--\ref{eq:intro-R-asym} of Theorem~\ref{thm:main}.

Finally, note that the intertwiner
\begin{equation*}\label{eq:Rp-intro}
\RR^+_{V_1,V_2}(s)\ceq(1\,2)\circ\RR^-_{V_2,V_1}(-s)^{-1}
\circ(1\,2): V_1\dtensor{s}V_2\to V_1\kmtensor{s}V_2\,
\end{equation*}
automatically provides a meromorphic twist in the opposite direction.
Thus, \eqref{eq:intro-Gauss-diag} reads
\begin{equation*}\label{eq:R-gauss-intro}
\RR^{\uparrow/\downarrow}_{V_1,V_2}(s) = 
\RR^+_{V_1,V_2}(s)\circ\RR_{V_1,V_2}^{0,\uparrow/\downarrow}(s)\circ
\RR^-_{V_1,V_2}(s)\,.
\end{equation*}

Problems \ref{prob:R0} and \ref{prob:Rm}
are solved in Theorems~\ref{thm:R0-main} and \ref{thm:R-}, respectively,
as we now spell out.

\subsection{Construction of $\RR^0(s)$}\label{ssec:intro-R0}
Let $\{t_{i,r}\}_{i\in\bfI, r\in\N}$ be the generators of the commutative subalgebra $\Ysupkm{0}\subset\Yhg$ (see Section \ref{ssec:tir}). These elements are primitive with respect to the Drinfeld coproduct. Let $\Prim(\Yhg)$ be the linear subspace spanned by them.

The operator $\RR^0(s)$ is required to satisfy the analogues of the intertwining equation and the cabling identites from Theorem~\ref{thm:main} \eqref{eq:intro-R-inter}--\eqref{eq:intro-R-QYBE} with respect to the Drinfeld coproduct. By the cabling identities, it follows that $\RR^0(s)$ is bound to have the form
\begin{equation*}\label{eq:R0-exp}
	\RR^0(s)=\exp(\Rl(s))
\end{equation*}
where\footnote{We expect that every primitive element belongs to $\Prim(\Yhg)$.} $\Rl(s)\in\Pseries{\Prim(\Yhg)^{\otimes 2}}{s^{-1}}$. Every such operator would also satisfy the intertwining equation restricted to $\Ysupkm{0}$.

Therefore, the only remaining constraint on $\RR^0(s)$ is given by the intertwining equation from raising and lowering operators in $\Yhg$. In Section~\ref{ssec:R0-thm-pf}, we prove that this forces $\Rl(s)$ to be a solution of the difference equation
\begin{equation}\label{eq:intro-L0-diff}
        (\sop-\sop^{-1})\det(\bfB(\sop)) \cdot \Rl(s) = \Rg(s),
\end{equation}
where
\begin{itemize}\itemsep0.25cm
        \item $\sop$ is the shift operator defined by $\sop\cdot f(s) = f(s-\hbar/2)$\,,

        \item $\bfB(\sop)$ is the symmetrized affine $\sop$--Cartan matrix
        (see Section \ref{ssec:T-cartan})\,,

        \item $\Rg(s) = \sum_{ij} \bfB(\sop)^{\ast}_{ji}\cdot \Rtau_{ij}(s)$, where
        $\bfB(\sop)^\ast$ is the adjoint matrix and
        $\Rtau_{ij}(s)$ is the element
        \[
        \Rtau_{ij}(s) \ceq \hbar^2
        \sum_{m\geqslant 1} m!s^{-m-1} \sum_{\begin{subarray}{c} a,b\geqslant 0 \\ a+b=m-1
        \end{subarray}} (-1)^a 
\frac{t_{i,a}}{a!}\otimes \frac{t_{j,b}}{b!}\,.
        \]
\end{itemize}

We observe that \eqref{eq:intro-L0-diff} is {\em irregular},
meaning that the difference
operator has order of vanishing $3$ at $\sop=1$ (Lemma~\ref{lem:order}) and
the right--hand side $\Rg(s)$ is $O(s^{-2})$.
Thus, there is no solution in $(\Ysupkm{0}\ten\Ysupkm{0})\fml{s^{-1}}$.
However, a direct computation  shows that the coefficient of $s^{-2}$ and $s^{-3}$ in
$\Rg(s)$ are central (see Lemma \ref{lem:Rg-terms}).

As central elements play no role in interwining equations, the last remark
allows us to regularize our difference equation and obtain
$\Rl(s)\in\Pseries{\Prim(\Yhg)^{\otimes 2}}{s^{-1}}$ as its unique formal solution
(see Proposition \ref{pr:Rl}).
We then prove that
$\ds\RR^0(s)=\exp(\Rl(s))$
is an abelian $R$--matrix for $\Yhg$, \ie it satisfies the same properties
from Section \ref{opening} with $\tau_s\ten\id\circ\Delta$ replaced by $\ddelta{s}$
(see Theorem~\ref{thm:R0-main}).

To determine the functional nature of this formal object, 
we evaluate \eqref{eq:intro-L0-diff} on a tensor product of
two representations from $\OY$. Focusing on one matrix entry
at a time, the problem becomes scalar valued.
In Theorem~\ref{appA:thm}, we prove the existence and uniqueness
of two fundamental solutions to the scalar analogue
of \eqref{eq:intro-L0-diff}.
\Omit{This part of the paper uses Watson's technique of computing
asymptotic expansions
of Laplace transforms, as well as the functional
properties of the latter. These are recalled in Section \ref{appA:watson}.}
Hence, \eqref{eq:intro-L0-diff}
admits two
fundamental solutions $\Rl_{V_1,V_2}^{\uparrow/\downarrow}(s)$
whose exponentials yield the meromorphic braiding $\RR^{0,\uparrow/\downarrow}_{V_1V_2}(s)$.
Moreover, their asymptotic expansions as
$\Re(s/\hbar)\to\pm\infty$ coincide and recover the action of $\Rl(s)$.

\begin{rem}
	The difference equation \eqref{eq:intro-L0-diff} makes sense for any type. In finite type,
	it reduces to the one considered in \cite{khoroshkin-tolstoy, sachin-valerio-III}.
	Its derivation, carried out in Proposition
	\ref{pr:tau} \eqref{tau:comm} and Corollary \ref{cor:Rg}, rests
	on the intertwining equation for $\RR^0(s)$. This is in contrast
	with the finite--type case, where
	Khoroshkin--Tolstoy \cite{khoroshkin-tolstoy} obtained the corresponding
	difference equation in order to compute the canonical tensor of a non--degenerate
	pairing. Their construction also uses the fact that for finite--type Cartan matrices,
	the determinant of the $q$--symmetrized Cartan matrix divides a $q$--number.
	The analogue of the non--degenerate pairing for affine Yangians is not known,
	and the latter statement is false (see the list of determinants
	given in Appendix \ref{app:QCM}).
\end{rem}

\subsection{Construction of $\RR^-(s)$}\label{ssec:intro-Rm} 
We now describe the solution to the problem \ref{prob:Rm}.
The rational twist $\RR^-(s)$ is obtained from the action of
a canonical element in a suitable completion of $(\Yhg\ten\Yhg)\fml{s^{-1}}$,
which is uniquely, and explicitly, determined by an intertwining equation.

To state this precisely, let $\upsigma_z:\Yhg\to\Yhg[z,z^{-1}]$
be the principal grading shift in $z$ (see Section \ref{ssec:pre-cop}), and
set
\Omit{
\[
\quad\mbox{and}\quad
\ddelta{s}:\Yhg\to(\Yhg\ten\Yhg)[s;s^{-1}]\!]\]
be, respectively, the principal grading shift in $z$ (see \ref{ssec:pre-cop})
and the \emph{deformed Drinfeld coproduct}. The latter
was introduced in \cite[Thm.~3.4]{GTLW}
and gives rise to the Drinfeld tensor product (see Section \ref{ssec:dr-tensor}).
Consider the homomorphisms}
\begin{equation*}
        \ddeltaz{s} \ceq   (\Id\otimes\upsigma_z)\circ\ddelta{s} \quad \text{ and }\quad \kmdeltaz{s}\ceq  (\tau_s\otimes \Id)\circ \Delta^z\,.
\end{equation*}
In addition, let $\Ysupkm{\pm}$ (resp. $\Ysupkm{0}$)
be the unital subalgebras
of $\Yhg$ generated by the Drinfeld
generators $\{x_{i,r}^\pm\}_{i\in\bfI, r\in\Z_{\geqslant 0}}$
(resp. by $\h$ and the commuting Drinfeld generators $\{t_{i,r}\}_{i\in\bfI, r\in\Z_{>0}}$)
equipped with the standard grading over the
root lattice $\Qaff$ (see Sections \ref{ssec: yangian}, \ref{ssec:tir}, and \ref{ssec:Q-grade}).

Our construction of $\RR^-(s)$ is based on the existence of a unique formal series
\begin{equation*}\label{eq:intro-Rm-z}
\RR^-(s, z)\ceq\sum_{\beta\in\Qaff_+}\Rmsb{\beta}{s} 
z^{\hgt{\beta}}\in(\Ysupkm{}\ten\Ysupkm{})\fml{s^{-1}}\fml{z},
\end{equation*}
where $\mathsf{ht}$ is the height function, satisfying the following three properties:
\begin{enumerate}[label=(P\arabic*)]\itemsep0.25cm
        \item\label{eq:intro-norm} {\bf Normalization:} $\Rmsb{0}{s}=1\otimes 1.$
        \item\label{eq:intro-triang} {\bf Triangularity:} for any $\beta\in\Qaff_+$,
        \begin{equation*}\label{eq:intro-Rm-norm-triang}
                \Rmsb{\beta}{s}\in(\Ysupkm{-}_{-\beta}\ten\Ysupkm{+}_{\beta})\fml{s^{-1}}\,.
        \end{equation*}
        \item\label{eq:intro-inter} {\bf Intertwining equation:} for any $x\in\Yhg$,
        \begin{equation*}
                \RR^-(s,z)\,\kmdeltaz{s}(x)= \ddelta{s}^{\!z}(x)\, \RR^-(s,z)\,.
        \end{equation*}
\end{enumerate}
In finite--type, the crucial observation in \cite[\S4]{GTLW} is that
the conditions \ref{eq:intro-norm}, \ref{eq:intro-triang} and
\ref{eq:intro-inter} for $x=t_{i,1}$, where $i$ is an arbitrary
node of the Dynkin diagram, already have a unique solution. A very
general rank $1$ reduction argument, valid for this solution,
reduces \ref{eq:intro-inter} for arbitrary $x\in\Yhg$, to
$x=x_0^{\pm}$ in $\Yhsl{2}$, which has been verified in \cite[\S4.8]{GTLW}. The rational twist alluded to in \ref{prob:Rm} is then obtained by specializing $z$ to $1$, which is possible since the matrix entries of $\RR^-(s,z)$ acting on a tensor product of
modules in $\OY$ become polynomials in $z$. 

Inspired by these results of \cite{GTLW},  we consider the analogue of
the {\em second embedding} of \S2.7 in {\em loc. cit.}:
\begin{equation*}\label{eq:intro-second-embedding}
\Topp:\h'\to\Ysupkm{0}\qquad\mbox{given by}\qquad \Topp(d_ih_i)=t_{i,1}\,,
\end{equation*}
where $\sfD\ceq (d_i)_{i\in\bfI}$ is the symmetrising diagonal matrix associated to $\g$ and
$\h'\subset\h$ is the span of the coroots $\{h_i\}_{i\in\bfI}$ (see Section
\ref{ssec:Notation}).
This is a natural choice, since the formulae of the standard coproduct are
explicit and fairly manageable in this case.
It is clear that if $h\in\h'$ and $\beta\in\Qaff_+$ are such that $\beta(h)\neq0$,
then, the intertwining equation \ref{eq:intro-inter} for $\Topp(h)$ produces an
explicit expression of the block $\Rmsb{\beta}{s}$ in terms of lower
blocks $\Rmsb{\gamma}{s}$, $\gamma<\beta$. For Yangians of finite type, this is
enough to determine $\RR^-(s,z)$ entirely. 

In our case, however, it fails almost completely.
Since the imaginary root $\delta$ vanishes on $\h'$, the resulting
system is unable to determine recursively the blocks
$\Rmsb{n\delta}{s}$, $n\in\Z_{>0}$, and, consequently, any block
$\Rmsb{\beta}{s}$ with $\beta>n\delta$ for some $n\in\Z_{>0}$.

\subsection{Extension of the second embedding} 
We overcome this issue 
by constructing an extension
of the second embedding
\begin{equation*}
        \Top:\h\to\Ysupkm{0}
\end{equation*}
satisfying the following properties (see Theorem~\ref{T:map-T}):
\begin{enumerate}[label=(T\arabic*)]\itemsep0.25cm
        \item For any $h\in\h$ and $i\in\bfI$, one has
        \[[\Top(h), x_{i,r}^\pm]=\pm\alpha_i(h)x_{i,r+1}^\pm\,.\]
        \item For any $h\in\h$ and $s\in\C$, one has
        \[\ad(\tau_s(\Top(h)))=\ad(\Top(h)+sh)\,.\]
        \item There is a family of translation--invariant elements
$\mathrm{Q}_{n\delta}\in\Ysupkm{-}_{n\delta}\ten\Ysupkm{+}_{n\delta}$, $n>0$, such that,
        for any $h\in\h$, one has
        \begin{equation*}
                \Delta^z(\Top(h))=\square\Top(h)+\hbar[h\ten 1, \Omega_z+\mathrm{Q}_z]\,,
        \end{equation*}
        where $\displaystyle \mathrm{Q}_z=\sum_{n>0}\mathrm{Q}_{n\delta}z^{n\,\hgt{\delta}}$ and $\square(y)=y\otimes 1 + 1\otimes y$ for all $y\in \Yhg$. 
\end{enumerate}
Relying on the map $\Top$, we reduce the intertwining equation \ref{eq:intro-inter}
to a system of linear equations. By choosing an element $\rho^\vee\in\h$ such that $\alpha_i(\rho^\vee)=1$ for all $i\in\bfI$, we obtain the following recursive expression
of the blocks of $\RR^-(s,z)$:
\begin{equation*}\label{eq:intro-Rm-recur}
        \begin{aligned}
                \RR^-(s)_\beta = \hbar\sum_{k\geqslant 0}
                \frac{\ad(\square\Top(\rho^\vee))^k}{(s\,\hgt{\beta})^{k+1}}
                \left(
                \sum_{\alpha\in\raff_+}\hgt{\alpha}\RR^-(s)_{\beta-\alpha}
                (\Omega_\alpha+\mathrm{Q}_\alpha)\right) \,.
        \end{aligned}
\end{equation*}
By \ref{eq:intro-norm}, the blocks are uniquely determined and satisfy \ref{eq:intro-triang}.
By a rank one reduction argument, we finally prove that $\RR^-(s,z)$ satisfies
also the intertwining equation \ref{eq:intro-inter}.

Therefore, it yields an intertwiner
\[\RR^-_{V_1,V_2}(s):V_1\kmtensor{s}V_2\to V_1\dtensor{s}V_2\,\]
which is readily seen to depend rationally on $s$, equal $\Id_{V_1\otimes V_2}$
at $s=\infty$, and
satisfy the cocycle equation \ref{prob:Rm}\footnote{Note that the cocycle equation only holds as an identity of rational functions valued in $\End(V_1\ten V_2\ten V_3)$, and does not seem to admit a natural lift to $\Yhg$, see \cite[Rmk.~4.1]{GTLW}.}.

\subsection{Outline of the paper}
We review the definition of $\Yhg$ and the basic properties of its category $\calO$ modules
in Sections~\ref{sec:Y} and \ref{sec:Y-rep}, respectively.
In Section~\ref{sec:Delta-T}, we introduce the transformation $\Top:\hkm\to\Ysupkm{0}$ and
prove its fundamental properties in Theorem~\ref{T:map-T}.
It is used in Section~\ref{sec:neg-R} to construct
the operator $\RR^-(s)$
as a rational twist relating the standard coproduct and the Drinfeld coproduct on category $\mathcal{O}$
modules (see Theorem~\ref{thm:R^--reps}).
In Section~\ref{sec:ab-R}, we construct the meromorphic braidings
$\RR^{0,\eta}$ in Theorem~\ref{thm:R0-main}
and prove that they both have the same asymptotic expansion given by the
formal abelian $R$--matrix $\RR^0(s)$.
The proof relies on well--known techniques to solve
additive difference equations, which we review in
Appendix~\ref{app:LapDE} for the reader's convenience.
Finally, in Appendices~\ref{app:Augmented} and \ref{app:QCM}, we provide the proofs of
Proposition~\ref{P:full-rank}
and Lemma~\ref{lem:order}, respectively.


\section{The affine Yangian $\Yhg$}\label{sec:Y}

\subsection{Affine Lie algebras}\label{ssec:Notation}
Throughout this paper, we fix a symmetrizable, indecomposable Cartan matrix
of {\em affine type} $\hbfA = (a_{ij})_{i,j\in\hbfI}$
and let $(d_i)_{i\in\hbfI}$ be the associated symmetrizing integers, taken to be positive and relatively prime.

As in \cite{guay-nakajima-wendlandt}, we further assume that $\hbfA$ is not of type $\mathsf{A}_1^{(1)}$  or $\mathsf{A}_2^{(2)}$. Let $\gkm$ be the Kac--Moody
Lie algebra associated to $\hbfA$. Recall that $\gkm$ is generated, as a Lie algebra,
by the following set of generators:
\begin{itemize}\itemsep0.25cm
\item A realization of $\hbfA$ \cite[\S 1.1]{kac}. That is, a $|\hbfI|+1$--dimensional
$\C$--vector space $\hkm$, together with two linearly independent subsets
$\{h_i\}_{i\in\hbfI}\subset\hkm$ and $\{\alpha_i\}_{i\in\hbfI}\subset\hkm^*$
related by $\alpha_j(h_i) = a_{ij}$ for every $i,j\in\hbfI$.

\item Raising and lowering operators $\{e_i,f_i\}_{i\in\hbfI}$.
\end{itemize}
These generators are subject to the usual Chevalley--Serre relations:
\begin{itemize}\itemsep0.25cm
\item $\hkm$ is abelian.
\item $[h,e_i]=\alpha_i(h)e_i$ and $[h,f_i]=-\alpha_i(h)f_i$
for every $h\in\hkm$ and $i\in\hbfI$.
\item $[e_i,f_j]=\delta_{ij}h_i$ for each $i,j\in \hbfI$.
\item  $\ad(e_i)^{1-a_{ij}}e_j = 0 = \ad(f_i)^{1-a_{ij}}f_j$ for each $i,j\in \hbfI$ with $i\neq j$. 
\end{itemize}
Let $\mathfrak{n}^+$ (resp. $\mathfrak{n}^-$) 
denote the Lie subalgebra
of $\g$ generated by $\{e_i\}_{i\in\hbfI}$ (resp. $\{f_i\}_{i\in\hbfI}$).

We recall some of the important ingredients in the theory of
affine Lie algebras. We mostly follow \cite[Ch. 6]{kac}. For the list
of tables of affine Dynkin diagrams, see \cite[\S 4.8]{kac}.
\begin{itemize}\itemsep0.25cm
\item 
By \cite[Prop. 4.7]{kac}, there is a unique $\hbfI$-tuple $(a_i)_{i\in \hbfI}$ of positive, relatively prime integers satisfying $\sum_{i\in\hbfI} a_{ji}a_i = 0$ for all $j\in\hbfI$. These numbers are listed explicitly in \cite[Ch. 4, Tables Aff 1-2-3]{kac}.
Following \cite[Thm.~5.6]{kac}, we define $\delta\in \hkm^\ast$ by 
$\delta=\sum_{i\in\hbfI} a_i\alpha_i \in\sum_i \Z_{>0}\alpha_i$. In particular, this element satisfies $\delta(h_j)=0$ for all $j\in \hbfI$. 

\item Using $A^T=DAD^{-1}$, we get that
$(a_i^{\vee}=d_ia_i)_{i\in\hbfI}$ are the coefficients of a
linear dependence relation among the rows of $\hbfA$. 
We record the linear relations:
$\sum_{i\in\bfI} a_i^{\vee}a_{ij}=0,\ \forall\ j\in\hbfI$.
Let $\caff \ceq   \sum_i a_i^{\vee} h_i$. Note that $\caff$ is central
in $\gkm$; as in \cite[\S6.2]{kac}, we call it the \textit{canonical central element} of $\gkm$.

\item Let $\haff\subset\hkm$ denote the span of $\{h_i:i\in\hbfI\}$.
Let $\raff$ denote the set of roots of $(\gkm,\hkm)$,
which comes naturally with the polarization $\raff=\raff_+\sqcup
\raff_-$. Let 
$\Qaff=\Z\raff$ be the root lattice, and $\Qaff_+
=\sum_{\alpha\in\raff_+} \N\alpha$. Set $\gaff = [\gkm,\gkm]$.

\end{itemize}

\subsection{The Yangian $\Ykm$ \cite{drinfeld-yangian-qaffine}}\label{ssec: yangian}

Let $\hbar\in\nC$. The Yangian $\Ykm$ is the unital, associative $\C$--algebra
generated by elements $\hkm \cup \{x^{\pm}_{i,r},\xi_{i,r}\}_{i\in\hbfI,
r\in\N}$, subject to the following relations:
\begin{enumerate}[font=\upshape, label=(Y\arabic*)]\itemsep0.25cm
\item\label{Y1} $\xi_{i,0}=d_ih_i\in\hkm$ for every $i\in\hbfI$.
$\hkm$ is abelian, and $\ds [h,\xi_{i,r}] = 0 = [\xi_{i,r}, \xi_{j,s}]$
for every $h\in\hkm$, $i,j\in\hbfI$ and $r,s\in\N$.

\item\label{Y2} For $i,j\in\hbfI$ and $s\in \N$: $
[\xi_{i,0}, x_{j,s}^{\pm}] = \pm d_ia_{ij} x_{j,s}^{\pm}$.

\item\label{Y3} For $i,j\in\hbfI$ and $r,s\in\N$:
\[[\xi_{i,r+1}, x^{\pm}_{j,s}] - [\xi_{i,r},x^{\pm}_{j,s+1}] =
\pm\hbar\frac{d_ia_{ij}}{2}(\xi_{i,r}x^{\pm}_{j,s} + x^{\pm}_{j,s}\xi_{i,r})\ .
\]

\item\label{Y4} For $i,j\in\hbfI$ and $r,s\in \N$:
\[
[x^{\pm}_{i,r+1}, x^{\pm}_{j,s}] - [x^{\pm}_{i,r},x^{\pm}_{j,s+1}]=
\pm\hbar\frac{d_ia_{ij}}{2}(x^{\pm}_{i,r}x^{\pm}_{j,s} 
+ x^{\pm}_{j,s}x^{\pm}_{i,r})\ .
\]
\item\label{Y5} For $i,j\in\hbfI$ and $r,s\in \N$:
$[x^+_{i,r}, x^-_{j,s}] = \delta_{ij} \xi_{i,r+s}$.

\item\label{Y6} Let $i\not= j\in\bfI$ and set $m = 1-a_{ij}$. For any
$r_1,\cdots, r_m, s\in \N$:
\[
\sum_{\pi\in\Sym_m}
\left[x^{\pm}_{i,r_{\pi(1)}},\left[x^{\pm}_{i,r_{\pi(2)}},\left[\cdots,
\left[x^{\pm}_{i,r_{\pi(m)}},x^{\pm}_{j,s}\right]\cdots\right]\right]\right]=0.
\]
\end{enumerate}

Note that it follows from this definition that the assignment 
\begin{equation*}
\sqrt{d_i} e_i\mapsto x_{i,0}^+, \quad \sqrt{d_i} f_i\mapsto  x_{i,0}^-, \quad 
d_i h_i \mapsto \xi_{i,0}, \quad h\mapsto h,
\end{equation*}
for all $i\in \hbfI$ and $h\in \hkm$, extends to an algebra homomorphism $U(\gkm)\to \Ykm$. We shall frequently work through this homomorphism without further comment.

We denote by $\Ysupkm{0}$ and $\Ysupkm{\pm}$ the unital subalgebras of
$\Ykm$ generated by $\hkm\cup\{\xi_{i,r}\}_{i\in\hbfI, r\in\Z_{\geqslant 0}}$
and $\{x_{i,r}^{\pm}\}_{i\in\hbfI, r\in\Z_{\geqslant 0}}$, respectively.
Let $\Ysupkm{\geqslant}$ (resp. $\Ysupkm{\leqslant}$) denote the
subalgebras of $\Ykm$ generated by $\Ysupkm{0}$ and $\Ysupkm{+}$
(resp. $\Ysupkm{0}$ and $\Ysupkm{-}$).

\subsection{Filtration}\label{ssec:filtration}
The Yangian $\Ykm$ is a filtered algebra once given the {\em loop
filtration}. That is, let $\deg(y_r)=r$, where $y$ is one of
$\xi_i$ or $x_i^{\pm}$. For each $k\geqslant 0$, set
\[
\mathbf{F}_k(\Ykm) \ceq   \text{Linear subspace spanned by monomials
of degree } \leqslant k\ .
\]

The induced filtration on $\Ykm^{\otimes 2}$ will again be
denoted by $\mathbf{F}_\bullet(\Ykm^{\otimes 2})$.

\subsection{Formal currents}\label{ssec:currents}

For each $i\in\hbfI$, we 
define $\xi_i(u),x_i^{\pm}(u)\in\Pseries{\Ykm}{u^{-1}}$ by
\[
\xi_i(u)=1+\hbar\sum_{r\in\N} \xi_{i,r}u^{-r-1} \quad \text{ and }\quad
x_i^{\pm}(u) = \hbar\sum_{r\in\N} x_{i,r}^{\pm} u^{-r-1}.
\]
The relations \ref{Y1}--\ref{Y6} can be written in terms
of these formal series, see \cite[Prop. 3.3]{sachin-valerio-2}.
Namely, it is proven that
these relations are equivalent to the following
identities:

\begin{enumerate}[font=\upshape, label=($\mathcal{Y}$\arabic*)]\itemsep0.25cm
\item\label{cY1} For any $i\in\hbfI$, $\xi_{i,0}=d_ih_i\in\hkm$.
For any $i,j\in\bfI$ and $h,h'\in\hkm$,
\[[\xi_i(u), \xi_j(v)]=0\qquad\qquad [\xi_i(u),h]=0\qquad\qquad [h,h']=0.\] 

\item\label{cY2} For any $i\in\hbfI$, and $h\in\hkm$,
$\ds[h,x^\pm_i(u)]=\pm\alpha_i(h)x^\pm_i(u)$.

\item\label{cY3} For any $i,j\in \hbfI$, and $a = \hbar d_ia_{ij}/2$
\[(u-v\mp a)\xi_i(u)x_j^{\pm}(v)=
(u-v\pm a)x_j^{\pm}(v)\xi_i(u)\mp 2a x_j^{\pm}(u\mp a)\xi_i(u).\]

\item\label{cY4}For any $i,j\in \hbfI$, and $a = \hbar d_ia_{ij}/2$
\begin{multline*}
(u-v\mp a) x_i^{\pm}(u)x_j^{\pm}(v)\\
= (u-v\pm a)x_j^{\pm}(v)x_i^{\pm}(u)
+\hbar\lp [x_{i,0}^{\pm},x_j^{\pm}(v)] - [x_i^{\pm}(u),x_{j,0}^{\pm}]\rp.
\end{multline*}

\item\label{cY5}For any $i,j\in \hbfI$
\[(u-v)[x_i^+(u),x_j^-(v)]=-\delta_{ij}\hbar\left(\xi_i(u)-\xi_i(v)\right).\]

\item\label{cY6}For any $i\neq j\in\hbfI$, $m=1-a_{ij}$, $r_1,\cdots, r_m\in
\N$, and $s\in \N$
\[\sum_{\pi\in\Sym_m}
\left[x^{\pm}_i(u_{\pi_1}),\left[x^{\pm}_i(u_{\pi(2)}),\left[\cdots,
\left[x^{\pm}_i(u_{\pi(m)}),x^{\pm}_j(v)\right]\cdots\right]\right]\right]=0.\]
\end{enumerate}
Here the relations \ref{cY1}--\ref{cY5} are identities in $\Ykm[u,v;u^{-1},v^{-1}]\!]$, while the Serre relation \ref{cY6} (for a fixed pair $i,j\in \hbfI$ with $i\neq j$) is an equality in the formal series space $\Ykm[\![u_1^{-1},\ldots,u_m^{-1},v^{-1}]\!]$, where $m=1-a_{ij}$.

\subsection{Alternate set of generators for $\Ysupkm{0}$}
\label{ssec:tir}
For each $i\in\hbfI$, let $t_i(u)$ and $B_i(z)$ be the formal series defined by
\begin{gather*}
t_i(u) = \hbar\sum_{r\geqslant 0} t_{i,r}u^{-r-1}
\ceq   \log(\xi_i(u)) = \sum_{n=1}^{\infty}
\frac{(-1)^{n-1}}{n} \lp\hbar \sum_{r\in\N}\xi_{i,r}u^{-r-1}\rp^n,\\
B_i(z) \ceq   \hbar\sum_{r\geqslant 0} t_{i,r}\frac{z^r}{r!}.
\end{gather*}
That is, $t_i(u)$ is the formal series logarithm of $\xi_i(u)$ and $B_i(z)$ is the formal Borel transform of $t_i(u)$. 
The elements $\{t_{i,r}\}_{i\in\hbfI,r\in\N}\subset\Ysupkm{0}$
are polynomials in $\{\xi_{i,r}\}_{i\in\hbfI,r\in\Z}$, and
together with $\hkm$, generated $\Ysupkm{0}$. We record the formulae
for the first few terms, for future use:
\begin{align}
t_{i,0} &= \xi_{i,0}, \label{eq:ti0} \\
t_{i,1} &= \xi_{i,1} - \frac{\hbar}{2}\xi_{i,0}^2, \label{eq:ti1}\\
t_{i,2} &= \xi_{i,2} - \hbar \xi_{i,1}\xi_{i,0} + 
\frac{\hbar^2}{3}\xi_{i,0}^3, \label{eq:ti2}\\
t_{i,3} &= \xi_{i,3} - \hbar\xi_{i,2}\xi_{i,0} - \frac{\hbar}{2}\xi_{i,1}^2
+\hbar^2\xi_{i,1}\xi_{i,0}^2 - \frac{\hbar^3}{4}\xi_{i,0}^4\ \label{eq:ti3} .
\end{align}

The following commutation relation was obtained in
\cite[\S 2.9]{sachin-valerio-1}:
\begin{equation}\label{eq:comm-Bi}
[B_i(z),x^{\pm}_{k,n}] = \pm \frac{e^{\frac{d_ia_{ik}\hbar}{2}z} - 
e^{-\frac{d_ia_{ik}\hbar}{2}z}}{z} \lp\sum_{r\geqslant 0} x^{\pm}_{k,n+r} \frac{z^r}{r!}\rp.
\end{equation}
Comparing coefficients of $z$, one obtains (see \cite[Remark 2.9]{sachin-valerio-1}):
\begin{equation*}
[t_{i,m},x_{k,n}^{\pm}] = 
\pm d_ia_{ik}
\sum_{\ell=0}^{\lfloor \frac{m}{2}\rfloor}
\cbin{m}{2\ell} \frac{(\hbar d_ia_{ik}/2)^{2\ell}}{2\ell+1}
x_{k,m+n-2\ell}^{\pm}\ .
\end{equation*}
A few special cases of this relation which will be particularly relevant to us are
\begin{align}
[t_{i,1},x_{j,n}^{\pm}] &= \pm d_ia_{ij} x_{j,n+1}^{\pm}\ ,\label{eq:commt1} \\
[t_{i,2},x_{j,n}^{\pm}] &= \pm d_ia_{ij} x_{j,n+2}^{\pm}
\pm \frac{\hbar^2}{12} (d_ia_{ij})^3 x_{j,n}^{\pm}, \label{eq:commt2}\ ,\\
[t_{i,3},x_{j,n}^{\pm}] &= \pm d_ia_{ij} x_{j,n+3}^{\pm}
\pm \frac{\hbar^2}{4} (d_ia_{ij})^3 x_{j,n+1}^{\pm}\label{eq:commt3}\ .
\end{align}

\subsection{Shift automorphism}\label{ssec: shift-yangian}

The group of translations of the complex plane acts on
$\Ykm$ as follows. For $s\in\C$, 
$\tau_s(h)=h,\ \forall\ h\in\hkm$, and 
$\tau_s(y(u)) = y(u-s)$, where $y$ is one of $\xi_i,x_i^{\pm}$.
In terms of modes, we have 
\[\tau_s(y_r) = \sum_{k=0}^r
\left(\begin{array}{c}r\\k\end{array}\right)
s^{r-k}y_k\ .
\]
Note that, since $t_i(u)=\log(\xi_i(u))$, the formula above also holds for $y=t_i$.\\

If $s$ is instead viewed as a formal variable, then we obtain an algebra homomorphism $\Ykm\to\Ykm[s]$, still denoted $\tau_s$. This is a filtered algebra homomorphism, provided $\Ykm[s]\cong \Ykm\otimes \C[s]$ is equipped with the standard tensor product filtration in which $\deg s=1$.

\subsection{$\Qaff$--grading}\label{ssec:Q-grade}
Viewed as a module over $\hkm$, we have
$\Ykm = \bigoplus_{\beta\in \Qaff}\Ykm_{\beta}$, where\
\[
\Ykm_\beta = \{y\in\Ykm : [h,y]=\beta(h)y,\ \forall\ 
h\in\hkm\}.
\]
This gives rise to a $\Qaff$-graded algebra structure on $\Ykm$ for which $\Ysupkm{\pm}$, $\Ysupkm{\scriptscriptstyle{\geqslant}}$ and $\Ysupkm{\scriptscriptstyle{\leqslant}}$ are all $\Qaff$-graded subalgebras.  In particular, we have
\[
\Ysupkm{\scriptscriptstyle{\geqslant}} = \bigoplus_{\beta\in\Qaff_+} \Ysupkm{\scriptscriptstyle{\geqslant}}_{\beta}
\quad \text{ and }\quad \Ysupkm{\scriptscriptstyle{\leqslant}} = \bigoplus_{\beta\in\Qaff_+} \Ysupkm{\scriptscriptstyle{\leqslant}}_{-\beta}.
\]

\subsection{The elements \texorpdfstring{$\caff_r$}{c_r}}\label{ssec:cr}
Recall from Section \ref{ssec:Notation} that $\caff = \sum_{i\in\hbfI} a_id_ih_i\in\hkm$ is the canonical central element of $\gkm$. We define higher order analogues of $\caff$ in $\Ykm$ by setting 
\begin{equation}\label{eq:Cr}
\caff_r\ceq   \sum_{i\in \hbfI} a_it_{i,r}  \quad \forall r\geqslant 0. 
\end{equation}
In particular, $\caff_0$ is the image of $\caff$ in $\Ykm$. In general, these elements do not belong to the center of $\Ykm$. However, we do have the following corollary of the relations \eqref{eq:commt1}--\eqref{eq:commt3}.
\begin{cor}\label{cor:c123}
The elements $\caff_0,\caff_1$ are central, while $\caff_2$ and $\caff_3$ satisfy
\begin{equation*}
[\caff_2,x_{j,s}^\pm]=\pm \frac{\hbar^2}{12} \mu_j x_{j,s}^\pm \quad 
\text{ and }\quad [\caff_3,x_{j,s}^\pm]=\pm \frac{\hbar^2}{4} \mu_j x_{j,s+1}^\pm,
\end{equation*}
for all $j\in \hbfI$ and $s\geqslant 0$, where $\mu_j$ is defined by 
\begin{equation}\label{eq:mu}
\mu_j=\sum_{i\in \hbfI}   a_i (d_i a_{ij})^3 \quad \forall j\in \hbfI.
\end{equation}
\end{cor}
That $\mathrm{ad}(\caff_j)$ is a non-trivial derivation for $j=2,3$ will play a crucial role in the constructions of this article.
\subsection{The coproduct $\Delta^{z}$}\label{ssec:pre-cop}

Finally, we recall the definition of the twisted standard coproduct $\Delta^z$ on $\Ykm$, which was introduced in \cite[\S6.1]{guay-nakajima-wendlandt}. 
By \cite[\S6.1]{guay-nakajima-wendlandt}, the assignment
\begin{equation*}
	\upsigma_z(x_{i,r}^\pm)=z^{\pm 1} x_{i,r}^\pm, \quad \upsigma_z(\xi_{i,r})=\xi_{i,r} 
	\quad \text{ and } \quad \upsigma_z(h)=h,
\end{equation*}
for all $i\in \hbfI$, $r\geqslant 0$ and $h\in \hkm$, extends to an injective  algebra homomorphism  
\begin{equation*}
	\upsigma_z:\Ykm\to \Ykm[z^{\pm 1}].
\end{equation*}
Note that if $\rho^{\vee}\in\hkm$ is chosen so that $\alpha_i(\rho^{\vee})=1$,
for every $i\in\hbfI$, then $\upsigma_z = \Ad(z^{\rho^{\vee}})$.
Next, for each $\beta\in \raff_+\cup\{0\}$, let $\Omega_\beta\in \gkm_{-\beta}\otimes \gkm_\beta$ be the canonical element defined by the restriction of the standard invariant form on $\gkm$ to $\gkm_{-\beta}\times \gkm_\beta$. We then define $\Omega_z$ and $\Omega_z^-$ in $(\gkm\otimes \gkm)[\![z]\!]$ by 
\begin{equation*}
	\Omega_z^-\ceq  \sum_{\beta\in \raff_+}\Omega_\beta z^{\mathsf{ht}(\beta)} \quad \text{ and }\quad \Omega_z\ceq  \Omega_0+\Omega_z^-,
\end{equation*}
where $\mathsf{ht}:\Qaff_+\to \Z_{\geqslant 0}$ is the additive \textit{height} function, defined on $\beta=\sum_{j}n_j \alpha_j$ by $\mathsf{ht}(\beta)=\sum_{j} n_j$. Note in particular that $\Omega_\beta z^{\mathsf{ht}(\beta)}=(\Id\otimes \upsigma_z)(\Omega_\beta)$ for each $\beta\in \Qaff_+$. 

By  Theorem 6.2 of \cite{guay-nakajima-wendlandt}, there is an algebra homomorphism 
\begin{equation*}
	\Delta^z:\Ykm\to \Ykm^{\otimes 2}[z^{-1};z]\!]
\end{equation*}
uniquely determined by the formulae
\begin{equation}\label{Delta^z:def}
	\begin{gathered}
		\Delta^z(y)=y\otimes 1 + 1\otimes \upsigma_z(y),\\
		\Delta^z(t_{i,1})=t_{i,1}\otimes 1 + 1\otimes t_{i,1}+\hbar[\xi_{i,0}\otimes 1, \Omega_z], 
	\end{gathered}
\end{equation}
for each $y\in \hkm\cup\{x_{j,0}^\pm\}_{j\in \hbfI}$ and $i\in \hbfI$. We will call $\Delta^z$ the (twisted) \textit{standard coproduct} on $\Ykm$. It is not coassociative or counital, but satisfies the twisted coalgebra relations
\begin{equation}\label{twisted-coalg}
	\begin{gathered}
		(\Delta^z \otimes \Id)\circ \Delta^{zw}=(\Id\otimes \Delta^w)\circ \Delta^z,\\
		(\veps\otimes \Id)\Delta^z(x)=1\otimes \upsigma_z(x)\quad \text{ and }\quad  (\Id\otimes \veps)\Delta^z(x)=x\otimes 1 \quad \forall \; x\in \Ykm,
	\end{gathered}
\end{equation}
where $\veps: \Ykm\to \C$ is the counit, defined by $\veps(y)=0$ for all $y\in \hkm\cup\{\xi_{i,r},x_{i,r}^\pm\}_{i\in \hbfI,r\geqslant 0}$.

It also follows easily from the formulae given above that $\Delta^z$ preserves
the loop filtration on $\Ykm$ introduced in Section \ref{ssec:filtration}
above. That is, for every $k\in\N$ we have 
\begin{equation}\label{eq:Delta-filt}
	\Delta^z(y)\in \mathbf{F}_k(\Ykm^{\otimes 2})[z^{-1},z]\!],\ 
	\forall\ y\in\mathbf{F}_k(\Ykm)\ .
\end{equation} 

We shall make use of the linear map $\square:\Ykm\to \Ykm\otimes \Ykm$ defined by 
\begin{equation*}
	\square(y)=y\otimes 1 + 1\otimes y \quad \forall \; y\in \Ykm.
\end{equation*}
Though it is not an algebra homomorphism, it satisfies $\square([x,y])=[\square(x),\square(y)]$ for all $x,y\in \Ykm$. In particular, by (4.13) of \cite{guay-nakajima-wendlandt}, we have  
\begin{equation}\label{Delta-xi1}
	\begin{aligned}
		\Delta^z(x_{i,1}^+)&=\square^z(x_{i,1}^+)-z\hbar[1\otimes x_{i,0}^+,\Omega_z], \\
		\Delta^z(x_{i,1}^-)&=\square^z(x_{i,1}^-)+\hbar[x_{i,0}^-\otimes 1,\Omega_z],
	\end{aligned}
\end{equation}
for each $i\in \hbfI$, where  $\square^z\ceq   (\Id\otimes \upsigma_z)\circ \square$.
Using the fact that $\xi_{i,r}-t_{i,r}\in\mathbf{F}_{r-1}(\Ykm)$,
it follows from the proof of \cite[Prop.~2.6]{GW-Poles} that
\begin{equation}\label{Delta-filt2}
	(\Delta^z-\square^z)(t_{i,r}) \in \Pseries{\mathbf{F}_{r-1}(\Ykm^{\otimes 2})}{z},
	\ \forall\ i\in\bfI,r\in\N.
\end{equation} 
Note that, since $t_{i,r}$ is a weight zero element, $\square^z(t_{i,r})=\square(t_{i,r})$.



\section{Representations of affine Yangians}
\label{sec:Y-rep}

In this section we recall the definitions of two different parameter--dependent 
tensor structures on category $\mathcal{O}$ representations
of $\Ykm$: the standard tensor product
$\kmtensor{s}$ and the Drinfeld tensor product $\dtensor{s}$.

\subsection{The category $\OY$}\label{ssec:cat-O}
A representation $V$ of $\Ykm$ is said to be in category
$\OY$ if its restriction to $U(\gkm)$ is in category
$\calO$. That is,

\begin{enumerate}[font=\upshape]\itemsep0.25cm
\item $V$ is a direct sum of finite--dimensional weight spaces.
\[
V = \bigoplus_{\mu\in\hkm^*} V[\mu],\ \dim(V[\mu])<\infty,
\]
where $V[\mu] \ceq   \{v\in V : h\cdot v = \mu(h)v,\ \forall h\in\hkm\}$.

\item Let $P(V) \ceq   \{\mu : V[\mu]\neq\{0\}\}$ be the set of
weight of $V$. Then, there exist $\lambda_1,\ldots,\lambda_r\in\hkm^*$
such that 
\begin{equation*}
P(V) \subset \bigcup_{j=1}^r \lambda_j - \Qaff_+.
\end{equation*}
\end{enumerate}
%
%
Given $V\in\OY$ and $s\in \C$, let $V(s)$ denote the pull--back
representation $\tau_s^*(V)$.

\begin{rem}
This is the Yangian analogue of the corresponding category
for quantum affinizations studied by Hernandez in 
\cite{hernandez-affinizations} (see also \cite[\S 3]{sachin-valerio-2}).
\end{rem}

The analogue of the following {\em rationality property} for quantum affine algebras was obtained in \cite[\S 6]{beck-kac} and \cite[Prop. 3.8]{hernandez-drinfeld}. The Yangian version, stated below,
can be found in \cite[Prop. 3.6]{sachin-valerio-2}.

\begin{prop}\label{pr:rationality}
Let $V$ be a representation of $\Ykm$ which is $\hkm$--diagonalizable,
with finite--dimensional weight spaces. Then, for every weight
$\mu\in\hkm^*$ of $V$, the generating series
\[
\xi_i(u)\in \Pseries{\End(V[\mu])}{u^{-1}},\qquad
x_i^{\pm}(u) \in \Pseries{\Hom(V[\mu],V[\mu\pm\alpha_i])}{u^{-1}},
\]
defined in \ref{ssec:currents} above, are the Taylor series expansions at $\infty$
of rational functions of $u$.
\end{prop}

\subsection{Matrix logarithms}\label{ssec:log}

Let $i\in\hbfI$, $V\in\OY$ and $\mu\in P(V)$. By the previous proposition, the operator $\xi_i(u)$
acting on the finite--dimensional weight space $V[\mu]$ becomes
a rational, abelian, function of $u\in\C$, taking value $1$
at $u=\infty$. Let $A\subset \nC$ be the set of poles
of $\xi_i(u)^{\pm 1}$. As shown in
\cite[Prop. 5.4]{sachin-valerio-III}, we can view $t_i(u)$ as
Taylor series near $\infty$ of a single--valued function defined
on the cut plane:
\[
t_i(u) = \log(\xi_i(u)) : \C\setminus \bigcup_{a\in A} [0,a] \to \End(V[\mu]).
\]

\subsection{Standard tensor product}\label{ssec:km-tensor}
Let $\Delta^z$ be the twisted standard coproduct introduced in Section \ref{ssec:pre-cop}.
Given $V_1,V_2\in\OY$, with action homomorphisms $\pi_\ell:\Ykm\to\End(V_\ell)$,
the composition
\[
(\pi_1\otimes \pi_2) \circ \Delta^z : \Ykm\to \End(V_1\otimes V_2)[z^{\pm 1}]
\]
 can be evaluated
at $z=1$ to yield an action of $\Ykm$ on $V_1\otimes V_2$
(see \cite[Cor.~6.9]{guay-nakajima-wendlandt}). The resulting representation of
$\Ykm$ is denoted $V_1\kmtensor{0}V_2$.
More generally, we set 
\begin{equation*}
V_1\kmtensor{s} V_2 \ceq   V_1(s)\kmtensor{0} V_2\quad \forall \; s\in \C.
\end{equation*}
The properties of this tensor product are summarized in the following theorem, which is a consequence of the results of \cite{guay-nakajima-wendlandt}; see also
\cite[Prop.~8.2]{GTLW}.
\begin{thm}\label{thm:guay-nakajima-wendlandt}
The category $\OY$ together with the tensor product $\kmtensor{s}$
is a (polynomial) tensor category. In more detail, we have
the following properties:
\begin{enumerate}[font=\upshape]\itemsep0.25cm
\item For $V_1,V_2\in\OY$, $V_1\kmtensor{s}V_2$ depends polynomially
in $s$.

\item The tensor product is compatible with the shift automorphism.
That is, for every $V_1,V_2\in\OY$, we have: 
\[
(V_1\kmtensor{s_1}V_2)(s_2) = V_1(s_2)\kmtensor{s_1} V_2(s_2),\qquad
V_1(s_1)\kmtensor{s_2} V_2 = V_1\kmtensor{s_1+s_2} V_2.
\]

\item Let $\Triv$ denote the $1$--dimensional, trivial representation
of $\Ykm$. Then the following natural identifications of vector spaces
are $\Ykm$--intertwiners.
\[
\Triv\kmtensor{s} V_2 \cong V_2,\qquad
V_1\kmtensor{s} \Triv \cong V_1(s).
\]

\item The tensor product is asssociative in the following sense.
For any $V_1,V_2,V_3\in\OY$, the natural identification of vector
spaces is a $\Ykm$--intertwiner:
\[
(V_1\kmtensor{s_1}V_2)\kmtensor{s_2} V_3 \cong
V_1\kmtensor{s_1+s_2} (V_2\kmtensor{s_2} V_3).
\]
\end{enumerate}
\end{thm}

\subsection{Drinfeld tensor product}\label{ssec:dr-tensor}
For quantum affine algebras, Drinfeld coproduct was introduced
in \cite{drinfeld-yangian-qaffine} as a formal algebra homomorphism. Its regularization
using the shift automorphisms was obtained by Hernandez \cite{hernandez-drinfeld}
and it was shown to depend rationally on the shift parameter.

The following version of the Drinfeld tensor product on
representation of Yangians was introduced in
\cite[\S 4]{sachin-valerio-III} (see also \cite[\S 3.2--3.4, and Prop.~8.1]{GTLW}).

Let $\ddelta{s}$ be the assignment on the generating set $\{h,\xi_{i,r},x_{i,r}^\pm\}_{h\in \hkm,i\in \hbfI,r\geq 0}$ of $\Ykm$, with values in $\Ykm\otimes\Ykm[s;s^{-1}]\negthinspace]$, defined as follows:

\begin{itemize}\itemsep0.25cm
\item For any $h\in\hkm$, $\ddelta{s}(h) = \square(h)=h\otimes 1 + 1\otimes h$.
\item For any $i\in\hbfI$, $\ddelta{s}(\xi_i(u)) = \xi_i(u-s)\otimes\xi_i(u)$.
Thus,
\[
\ddelta{s}(\xi_{i,r}) = \tau_s(\xi_{i,r})\otimes 1 + 1\otimes\xi_{i,r}
+\hbar\sum_{p=0}^{r-1} \tau_s(\xi_{i,p})\otimes \xi_{i,r-1-p}\ .
\]
Note that the elements $\{t_{i,r}\}$ introduced in Section \ref{ssec:tir}
are primitive with respect to the Drinfeld coproduct. That is,
\[
\ddelta{s}(t_{i,r}) = \tau_s(t_{i,r})\otimes 1 + 1\otimes t_{i,r}\ .
\]

\item For each $i\in\hbfI$, we have:
\begin{align*}
\ddelta{s}(x^+_{i,r}) &= \tau_s(x_{i,r}^+)\otimes 1 + 1\otimes x_{i,r}^+ \\
&\phantom{=}
+ \hbar \sum_{N\geqslant 0} s^{-N-1} \lp\sum_{n=0}^N (-1)^{n+1}
\cbin{N}{n} \xi_{i,n}\otimes x_{i,r+N-n}^+ \rp ,\\
\ddelta{s}(x^-_{i,r}) &= \tau_s(x^-_{i,r})\otimes 1 + 1\otimes x^-_{i,r}
+ \hbar \sum_{p=0}^{r-1} \tau_s(x^-_{i,r-1-p})\otimes \xi_{i,p} \\
&\phantom{=}+ \hbar \sum_{N=0}^{\infty} s^{-N-1} \left( \sum_{n=0}^N (-1)^n 
\left(\begin{array}{c} N \\ n \end{array}\right)
x^-_{i,n}\otimes \xi_{i,r+N-n} \right) .
\end{align*}

\end{itemize}
It was shown in \cite[Thm.~3.4]{GTLW} that the analogue of this assignment for the Yangian $\Yhgfin$ of a  finite-dimensional simple Lie algebra $\gfin$ defines an algebra homomorphism  $\Yhgfin\to\Yhgfin\otimes\Yhgfin[s;s^{-1}]\negthinspace]$. This result extends naturally
to any Kac--Moody algebra satisfying the condition that each rank $2$ diagram subalgebra is
of finite type, since the defining relations of the Yangian are inherently
of rank $2$ type. In particular, this applies to the affine Kac--Moody algebra $\gkm$: the above assignment extends to an algebra homomorphism
\begin{equation*}
\ddelta{s}:\Ykm\to \Ykm\otimes\Ykm[s;s^{-1}]\negthinspace].
\end{equation*}

In contrast with the infinite sums encountered in the definition
of $\Delta$ from Section \ref{ssec:pre-cop}, the infinite sums written above do not
truncate. It is a consequence of the rationality property provided by Proposition
\ref{pr:rationality}, that these formal Laurent series
in $s^{-1}$ become rational functions of $s$, once
evaluated on $V_1\otimes V_2$. More precisely, given $V_1,V_2\in\OY$,
the image of the composition
\[
(\pi_1\otimes \pi_2)\circ \ddelta{s} : \Ykm
\to
\End(V_1\otimes V_2)[s;s^{-1}]\negthinspace]
\]
lies in the subspace of $\End(V_1\ten V_2)$-valued rational functions in $s$.
\footnote{Since $V_1\ten V_2$ is in general infinite-dimensional, this phrase is to be understood at the level of weight spaces. By a slight abuse of terminology, for every $V\in\OY$,  we say that
	$\varphi\in \End(V)[s,s^{-1}]\negthinspace]$ is rational if, for every
	$\lambda,\mu\in\h^*$, its restriction to the corresponding finite-dimensional weight spaces
	$\varphi_{\lambda,\mu}\in \Hom(V[\lambda],V[\mu])[s,s^{-1}]\negthinspace]$
	is a rational function of $s$.}
We let $V_1\dtensor{s}V_2$ denote the resulting representation
of $\Ykm$. The following theorem summarizes the results of
\cite[\S 3.2--3.4]{GTLW}.

\begin{thm}\label{thm:GTLW}
The category $\OY$ together with the tensor product $\dtensor{s}$
is a (rational) tensor category. In more detail, we have
the following properties:
\begin{enumerate}[font=\upshape]\itemsep0.25cm
\item For $V_1,V_2\in\OY$, $V_1\dtensor{s}V_2$ depends rationally
on $s$.

\item The tensor product is compatible with the shift automorphism.
That is, for every $V_1,V_2\in\OY$, we have: 
\[
(V_1\dtensor{s_1}V_2)(s_2) = V_1(s_2)\dtensor{s_1} V_2(s_2),\qquad
V_1(s_1)\dtensor{s_2} V_2 = V_1\dtensor{s_1+s_2} V_2.
\]

\item Let $\Triv$ denote the $1$--dimensional, trivial representation
of $\Ykm$. Then the following natural identifications of vector spaces
are $\Ykm$--intertwiners.
\[
\Triv\dtensor{s} V_2 \cong V_2,\qquad
V_1\dtensor{s} \Triv \cong V_1(s).
\]

\item The tensor product is associative in the following sense.
For any $V_1,V_2,V_3\in\OY$, the natural identification of vector
spaces is a $\Ykm$--intertwiner:
\[
(V_1\dtensor{s_1}V_2)\dtensor{s_2} V_3 \cong
V_1\dtensor{s_1+s_2} (V_2\dtensor{s_2} V_3).
\]
\end{enumerate}

\end{thm}

\begin{rem}
We wish to stress the point that $\ddelta{s}$ is {\em not}
coassociative. The associativity of $\dtensor{s}$ rests on
an identity among rational functions of two variables,
which does not seem to have a lift at the level of algebra
(see \cite[Remark 3.1]{GTLW}).
\end{rem}


\section{The transformation \texorpdfstring{$\Top$}{T}}
\label{sec:Delta-T}

In this section, we describe a linear map
$\Top:\hkm\to\Ysupkm{0}$, which plays a crucial role in the construction
of $\RR^-(s)$ carried out in Section \ref{sec:neg-R} below.

\subsection{The transformation $\Top$}\label{ssec:T-main}

Let $\anode\in\bfI$ denote the {\em extending vertex}, as in
\cite[\S6.1]{kac} and let $d\in\hkm$ denote the {\em scaling element},
chosen so as to have $\alpha_i(d) = \delta_{i,\anode}$ (see \cite[\S6.2]{kac}).
We need the following result, which is proven in detail in Appendix \ref{app:Augmented}.
\begin{prop}\label{P:full-rank} Let $\bfB=(d_ia_{ij})_{i,j\in \hbfI}$ be the symmetrized Cartan matrix of $\gkm$, and let $\underline{\mu}=(\mu_j)_{i\in \hbfI}\in \mathbb{Z}^{\hbfI}$,
	where $\mu_j$ is defined by \eqref{eq:mu}. 
	Then the augmented matrix $(\bfB\,|\, \underline{\mu})$ has rank $|\hbfI|$. 
\end{prop}
Since the augmented matrix $(\mathbf{B}\,|\, \underline{\mu})$ has full rank, the unit vector $(\delta_{j,\anode})_{j\in \hbfI}\in \mathbb{Q}^{\hbfI}$ lies in its range. In particular, there exists a tuple $(\zeta_i)_{i\in \hbfI}\in \mathbb{Q}^{\hbfI}$ and $\zeta\in \mathbb{Q}^\times$ such that 
\begin{equation*}
	\frac{1}{4}\zeta\mu_j +\sum_{i\in \hbfI}d_ja_{ji}\zeta_i =\delta_{j,\anode} \quad \forall \; j\in \hbfI.
\end{equation*} 
%
Since $\{d\}\cup \{d_i h_i\}_{i\in \hbfI}$ is a basis of $\hkm$, the assignment $\Top:\{d\}\cup \{d_i h_i\}_{i\in \hbfI}\to \Ysupkm{0}$ defined by 
\begin{equation*}
	\Top(d)=\sum_{i\in \hbfI}\zeta_i t_{i,1} + \frac{1}{\hbar^2}\zeta \caff_3\quad \text{ and }\quad \Top(d_ih_i)=t_{i,1} \quad \forall\; i\in \hbfI
\end{equation*}
uniquely extends to a linear map $\Top:\hkm\to \Ysupkm{0}$. The following theorem  provides the main result of this section. 
\begin{thm}\label{T:map-T}
	The linear map $\Top$ has the following properties: 
	\begin{enumerate}[font=\upshape]\setlength{\parskip}{3pt}\itemsep0.25cm
		\item\label{map-T:1} For each $h\in \hkm$, $j\in \hbfI$ and $r\geqslant 0$, one has 
		\begin{equation*}
			[\Top(h),x_{j,r}^\pm]=\pm \alpha_j(h) x_{j,r+1}^\pm.
		\end{equation*}
		\item\label{map-T:2} For each $h\in \hkm$ and $s\in \C$, one has 
		\begin{equation*}
			\mathrm{ad}(\tau_s(\Top(h)))=\mathrm{ad}(\Top(h)+sh).
		\end{equation*}
		\item\label{map-T:3} For each $h\in \hkm$, $\Delta^z\circ \Top$ satisfies 
		\begin{equation*}
			\Delta^z(\Top(h))=\Top(h)\otimes 1 + 1\otimes \Top(h) + \hbar[h\otimes 1, \Omega_z+\mathrm{Q}_z] 
		\end{equation*}
		where  $\mathrm{Q}_z$ is a formal series in $z$ with the following properties:
		\vspace{0.25cm}
		\begin{enumerate}[font=\normalshape, label=\emph{\roman*})]\itemsep0.25cm
			\item\label{Q_z:i} $\mathrm{Q}_z=\sum_{n>0} \mathrm{Q}_{n\delta} z^{n\mathsf{ht}(\delta)}$, where each coefficient $\mathrm{Q}_{n\delta}$ satisfies
			\begin{equation*}
				\mathrm{Q}_{n\delta}\in \Ysupkm{-}_{-n\delta}\otimes \Ysupkm{+}_{n\delta}.
			\end{equation*}
			\item\label{Q_z:ii} For each $a,b\in \C$, one has 
			\begin{equation*}
				(\tau_a\otimes \tau_b)(\mathrm{Q}_z)=\mathrm{Q}_{z}.
			\end{equation*}

			\item \label{Q_z:iii}The tensor factors of $\mathrm{Q}_z$ are primitive:
			\begin{align*}
				(\Delta^w\otimes \Id)(\mathrm{Q}_{z})&=\mathrm{Q}_{z}^{13}+\mathrm{Q}_{z/w}^{23} \\
				(\Id\otimes \Delta^w)(\mathrm{Q}_z)&=\mathrm{Q}_z^{12}+\mathrm{Q}_{zw}^{13}.
			\end{align*}
			\item\label{Q_z:iv}  $\mathrm{Q}_{z}$ belongs to the centralizer of $U(\gaff)^{\otimes 2}$ in $\Ykm^{\otimes 2}[\![z]\!]$.
		\end{enumerate}
	\end{enumerate}
\end{thm}
\begin{pf} Let us begin with some preliminary observations that will be useful throughout the proof. 
	Since the constant $\zeta$ is nonzero, we can set
	$
	\mathsf{h}\ceq  \frac{\hbar^2}{\zeta}(d-\sum_{i\in \hbfI}\zeta_id_i h_i )\in \hkm
	$.
	Then $\{\mathsf{h}\}\cup\{d_ih_i\}_{i\in \hbfI}$ is a basis of $\hkm$, and we have 
	\begin{equation}\label{3C_2-preimage}
		\begin{gathered}
			\Top(\mathsf{h})=\frac{\hbar^2}{\zeta}\bigg(\Top(d)-\sum_{i\in \hbfI}\zeta_it_{i,1} \bigg)=\caff_3,\\ 
			[\mathsf{h},x_{j,r}^\pm]=\pm\frac{\hbar^2}{\zeta}\bigg(\delta_{j,\anode}-\sum_{i\in \hbfI}\zeta_id_ia_{ij} \bigg)x_{j,r}^\pm	=\pm \frac{\hbar^2}{4} \mu_jx_{j,r}^\pm,
		\end{gathered}
	\end{equation}
	for all $j\in \hbfI$ and $r\geqslant 0$. It thus follows by Corollary \ref{cor:c123} that we have the equality of operators $\mathrm{ad}(\mathsf{h})=\mathrm{ad}(3\caff_2)$ on $\Ykm$. 
	
	Consider now 
	Parts \eqref{map-T:1} and \eqref{map-T:2} of the theorem. They clearly hold when $h\in \{d_ih_i\}_{i\in \hbfI}$, so it is sufficient to verify them for $h=\mathsf{h}$. For Part \eqref{map-T:1}, this is immediate from \eqref{3C_2-preimage} and Corollary \ref{cor:c123}. As for Part \eqref{map-T:2}, since $\Top(\mathsf{h})=\caff_3$, we have 
	\begin{align*}
		\mathrm{ad}(\tau_s(\Top(\mathsf{h})))&=\mathrm{ad}(\caff_3+3s\caff_2+3s^2\caff_1+s^3 \caff_0)\\
		&=\mathrm{ad}(\caff_3+3s\caff_2)\\
		&=\mathrm{ad}(\Top(\mathsf{h})+s\mathsf{h}),
	\end{align*}
	where the second and third equalities follow since $\caff_0$ and $\caff_1$ are central and  $\mathrm{ad}(\mathsf{h})=\mathrm{ad}(3\caff_2)$.
	
	Let us now turn to Part \eqref{map-T:3}. We will construct the unique series $\mathrm{Q}_z$ satisfying
	\begin{equation*}
		\Delta^z(\Top(h))=\Top(h)\otimes 1 + 1\otimes \Top(h) + \hbar[h\otimes 1, \Omega_z+\mathrm{Q}_z] 
	\end{equation*}
	for all $h\in \hkm$, in addition to the conditions \ref{Q_z:i}--\ref{Q_z:iv}, in several steps which will be carried out in Sections \ref{ssec:T-beta-nperp}--\ref{ssec:Q_z-exp}. In order to explain these steps, we introduce $\Thz{z}{}{}\in \Pseries{\Ykm^{\otimes 2}}{z}$ via the following
	equation:
	\begin{equation}\label{Theta_z-def}
		\Delta^z(\caff_3) = \caff_3\otimes 1 + 1\otimes \caff_3 + \hbar\Thz{z}{}{}.
	\end{equation} 
	By \eqref{eq:Cr}, $\caff_3=\sum_{i\in\bfI}a_it_{i,3}$. Thus, it follows from \cite[Prop. 2.9]{GW-Poles} that $\Thz{z}{}{}=\sum_{\beta>0}\Theta_\beta z^{\mathsf{ht}(\beta)}$, where the summation is taken over all nonzero $\beta \in \Qaff_+$ and 
	\begin{equation*}
		\Theta_\beta \in \Ysupkm{{\scriptscriptstyle\leqslant}}_{-\beta}\otimes \Ysupkm{{\scriptscriptstyle\geqslant}}_{\beta}.
	\end{equation*} 
	Next, let $\EuScript{K}_\beta\ceq  \Theta_\beta+\beta(\mathsf{h})\Omega_\beta$ for each nonzero $\beta\in \Phi_+$, so that 
	\begin{equation}\label{K_z}
		\EuScript{K}_z\ceq  \sum_{\beta>0}\EuScript{K}_\beta z^{\mathsf{ht}(\beta)}=\Thz{z}{}{}-[\mathsf{h}\otimes 1, \Omega_z].
	\end{equation}
	These elements are of interest to us because of the equation
	\[
	\Delta^z(\caff_3) = \caff_3\otimes 1 + 1\otimes \caff_3 + \hbar
	\left(\EuScript{K}_z + [\mathsf{h}\otimes 1, \Omega_z]\right).
	\]
	Thus, for the equation in \eqref{map-T:3} to hold, we must find
	$\mathrm{Q}_z$ so that $[\mathsf{h}\otimes 1, \mathrm{Q}_z]=\EuScript{K}_z$
	and establish its stated properties. We proceed as follows: 
	
	\begin{enumerate}[label=\emph{\arabic*})]\setlength{\parskip}{5pt}\itemsep0.25cm
		\item In Section \ref{ssec:T-beta-nperp}, we will prove that $\EuScript{K}_\beta=0$ for all $\beta \in \Qaff_+\!\setminus \Z_{\geqslant 0}\delta$, and that $\EuScript{K}_z$ belongs to the centralizer of $U(\gaff)^{\otimes 2}$ in $\Ykm^{\otimes 2}[\![z]\!]$. 
		This will be achieved in Proposition \ref{P:Theta} with the help of Lemma \ref{L:Theta_z}, which computes the commutation relations between $\Thz{z}{}{}$ and $\Delta^z(y)$, for all $y\in \{x_{i,0}^\pm, t_{i,1}\}_{i\in \hbfI}$.
		
		\item By the previous step, $\EuScript{K}_z$ takes the form $\EuScript{K}_z=\sum_{n>0}\EuScript{K}_{n\delta}z^{n\mathsf{ht}(\delta)}$. In Sections \ref{ssec:T-Gamma_i} and \ref{ssec:T-Gamma_i-2}, we will explicitly compute each coefficient $\EuScript{K}_{n\delta}$ and study some of their properties; see Lemma \ref{L:Gamma-i} and Proposition \ref{P:K-exp}.
		
		\item In Section \ref{ssec:Q_z-exp}, we finally define $\mathrm{Q}_z=\sum_{n>0}\mathrm{Q}_{n\delta}z^{n\mathsf{ht}(\delta)}$ by setting 
		$
		\mathrm{Q}_{n\delta}\ceq  -\frac{1}{n\delta(\mathsf{h})}\EuScript{K}_{n\delta}
		$
		for each positive integer $n$, so that $[\mathsf{h}\otimes 1,\mathrm{Q}_z]=\EuScript{K}_z$. In Proposition \ref{P:Q_z}, we show that $\mathrm{Q}_z$ has all the properties stated in Part \eqref{map-T:3} of the theorem.  \qedhere
	\end{enumerate}
\end{pf}

\subsection{Proof that $\EuScript{K}_\beta=0$ for $\beta\in \Qaff_+\!\setminus\! \Z_{\geqslant 0}\delta$}\label{ssec:T-beta-nperp}
Our goal in this section is to show that the series $\EuScript{K}_z\in \Ykm^{\otimes 2}[\![z]\!]$ defined in \eqref{K_z} belongs to the centralizer of $U(\gaff)^{\otimes 2}$ in $\Ykm^{\otimes 2}[\![z]\!]$ and, consequently, that its component $\EuScript{K}_\beta$ is zero for $\beta\in \Qaff_+\!\setminus\! \Z_{\geqslant 0}\delta$.
We begin with the following lemma, which spells out some of the commutation relations satisfied by $\Thz{z}{}{}$.
\begin{lem}\label{L:Theta_z}
	For each $i\in \hbfI$, one has 
	\begin{gather*}
		[\Thz{z}{}{},\Delta^z(t_{i,1})]=-[\xi_{i,0}\otimes 1,[\square(\caff_3),\Omega_z]],\\
		[\Thz{z}{}{},\Delta^z(x_{i,0}^+)]=\frac{\hbar^2}{4}z\mu_i[\Omega_z,1\otimes x_{i,0}^+] \quad \text{ and }\quad [\Thz{z}{}{},\Delta^z(x_{i,0}^-)]=\frac{\hbar^2}{4}\mu_i[\Omega_z,x_{i,0}^-\otimes 1].
	\end{gather*}
\end{lem}
\begin{pf}
	Since  $\caff_3$ and $t_{i,1}$ commute and $\Delta^z$ is an algebra homomorphism, we have 
	\begin{equation*}
		[\hbar\Thz{z}{}{},\Delta^z(t_{i,1})]=[\Delta^z(\caff_3),\Delta^z(t_{i,1})]-[\square(\caff_3),\Delta^z(t_{i,1})]=-[\square(\caff_3),\Delta^z(t_{i,1})].
	\end{equation*}
	Since $\square(\caff_3)$ and $\square(t_{i,1})$ commute, we obtain 
	\begin{equation*}
		[\hbar\Thz{z}{}{},\Delta^z(t_{i,1})]=-\hbar[\square(\caff_3),[\xi_{i,0}\otimes 1, \Omega_z]]=-\hbar[\xi_{i,0}\otimes 1,[\square(\caff_3), \Omega_z]],
	\end{equation*}
	which yields the first identity of the lemma. Similarly, we have 
	\begin{equation*}
		[\hbar\Thz{z}{}{}, \Delta^z(x_{i,0}^\pm)]=\Delta^z([\caff_3,x_{i,0}^\pm])-\square^z[\caff_3,x_{i,0}^\pm]=\pm\frac{\hbar^2}{4}\mu_i (\Delta^z-\square^z)(x_{i,1}^\pm),
	\end{equation*}
	where we recall that $\square^z=(1\otimes \upsigma_z)\circ \square$. 
	The second and third relations of the lemma now follow from \eqref{Delta-xi1}.
\end{pf}
\begin{prop}\label{P:Theta}
	The element $\EuScript{K}_z$ has the following properties:
	\begin{enumerate}[font=\upshape]\itemsep0.25cm
		\item\label{Theta:2} For each $x\in \gaff$, one has 
		\begin{equation*}
			[\EuScript{K}_z, x\otimes 1]=0=[1\otimes x,\EuScript{K}_z].
		\end{equation*}
		\item\label{Theta:3} $\EuScript{K}_\beta=0$ for all $\beta \in \Qaff_+\!\setminus \Z_{\geqslant 0}\delta$.

	\end{enumerate}
\end{prop}

\begin{pf}
	If Part \eqref{Theta:2} holds, then by taking $x=\xi_{i,0}$ in the relation $[\EuScript{K}_z, x\otimes 1]=0$ we find that 
	\begin{equation*}
		0=[\EuScript{K}_z, \xi_{i,0}\otimes 1]=\sum_{\beta>0} (\alpha_i,\beta) \EuScript{K}_\beta z^{\mathsf{ht}(\beta)} \quad \forall \; i\in \hbfI.
	\end{equation*}
	Projecting onto the $\Ykm_{-\beta}\otimes \Ykm_{\beta}$-component, we obtain $(\alpha_i,\beta) \EuScript{K}_\beta=0$ for all $i\in \hbfI$, and thus $\EuScript{K}_\beta=0$ provided $(\beta,\alpha_i)\neq 0$ for some $i\in \hbfI$. This shows that Part \eqref{Theta:3} follows from Part \eqref{Theta:2}.

	Let us now turn to proving Part \eqref{Theta:2}. 
	Define 
	\begin{equation*}
		\Thz{z;s}{}{}\ceq  
		(\tau_s\ten\Id)(\Thz{z}{}{})-\Thz{z}{}{}
		=\sum_{\beta>0}\Thz{1;s}{}{\beta}z^{\mathsf{ht}(\beta)}
	\end{equation*}	
	where
	\begin{equation*}
		\Thz{1;s}{}{\beta}\ceq  (\tau_s\otimes \Id)(\Theta_\beta)-\Theta_\beta \in s(\Ysupkm{{\scriptscriptstyle\leqslant}}_{-\beta}\otimes \Ysupkm{{\scriptscriptstyle\geqslant}}_{\beta})[s]
	\end{equation*}
	for each nonzero $\beta \in \Qaff_+$.
	Note that each of these is a polynomial of degree at most $2$. Indeed, 
	by \eqref{Delta-filt2}, $(\Delta^z-\square)(t_{i,3})$ belongs to the subspace $\mathbf{F}_2(\Ykm^{\otimes 2})[\![z]\!]$ of $\Ykm^{\otimes 2}[\![z]\!]$.
	So, the same is true for $(\Delta^z-\square)(\caff_3)$ and thus $\Thz{z}{}{}$.  
	As $\tau_s\otimes \Id$ sends any element of $\mathbf{F}_k(\Ykm^{\otimes 2})$ to a polynomial in $s$ of degree at most $k$ (see Section \ref{ssec: shift-yangian} above), the assertion follows. This allows us to write
	\begin{equation*}
		\Thz{z;s}{}{}\ceq \Thz{z}{1}{}s+ \Thz{z}{2}{}s^2,
	\end{equation*}
	where $\Thz{z}{1}{}, \Thz{z}{2}{}\in \Ykm^{\otimes 2}[\![z]\!]$. 
	Furthermore,  $\Thz{z;s}{}{}$ is $\gkm$-invariant, in the sense that it satisfies
	\begin{equation}\label{Theta_z(s)-inv}
		[\Delta^z(x),\Thz{z;s}{}{}]=0 \quad \forall  \;x\in \gkm. 
	\end{equation}
	For $x\in \hkm$ this is clear from the form of $\Thz{z;s}{}{}$, and for $x=x_{i}^\pm$ this follows by applying  $\tau_s\otimes \Id$ to the second line of identities from Lemma \ref{L:Theta_z}. 
	
	Applying $\tau_s\otimes \Id$ instead to the first equation of Lemma \ref{L:Theta_z} while using that $\caff_0$ and $\caff_1$ are central, we obtain 
	\begin{equation*}
		[\Thz{z}{}{}+\Thz{z;s}{}{},\Delta^z(t_{i,1})+s\xi_{i,0}\otimes 1]=-[\xi_{i,0}\otimes 1,[\square(\caff_3)+3s\caff_2\otimes 1,\Omega_z]].
	\end{equation*}
	Expanding both sides and subtracting the first relation of Lemma \ref{L:Theta_z} yields
	\begin{equation}\label{Theta_z}
		s[\Thz{z}{}{},\xi_{i,0}\otimes 1]+s[\Thz{z;s}{}{},\xi_{i,0}\otimes 1]+[\Thz{z;s}{}{},\Delta^z(t_{i,1})]=-3s[\xi_{i,0}\otimes 1,[\caff_2\otimes 1,\Omega_z]]
	\end{equation}
	Comparing powers of $s$, we immediately deduce that 
	\begin{equation}\label{Theta_z'}
		[\Thz{z}{2}{},\xi_{i,0}\otimes 1]=0 \quad \text{ and }\quad [\Thz{z}{1}{},\xi_{i,0}\otimes 1]=-[\Thz{z}{2}{},\Delta^z(t_{i,1})].
	\end{equation}
	It follows from the first identity that the $\beta$-component $\Thz{z}{2}{\beta}$ of $\Thz{z}{2}{}$ is zero unless $(\beta,\alpha_i)=0$ for all $i\in \hbfI$. Furthermore, by using this identity, and applying $\ad(\xi_{i,0}\otimes 1)$ to
	\eqref{Theta_z(s)-inv} with $x=x^{\pm}_{i,0}$, we deduce that 
	\begin{equation}\label{Theta_z^2}
		[\Thz{z}{2}{},x\otimes 1]=0=[\Thz{z}{2}{},1\otimes x] \quad \forall\; x\in \gaff.
	\end{equation}
	Consequently, the bracket $[\Thz{z}{2}{},\Delta^z(t_{i,1})]$ coincides with $[\Thz{z}{2}{},\square(t_{i,1})]$, and hence the second equation of \eqref{Theta_z'} yields
	\begin{equation*}
		[\Thz{z}{1}{\beta},\xi_{i,0}\otimes 1]=-[\Thz{z}{2}{\beta},\square(t_{i,1})] \quad \forall \; \beta\in \Qaff_+\setminus\{0\}, \, i\in \hbfI.
	\end{equation*}
	Note that the left--hand side is zero if $\beta\in \alpha_i^\perp$, and the right--hand side is zero if $\beta\not\in \alpha_i^\perp$. Hence both sides are identically zero for all nonzero $\beta \in \Qaff_+$.
	Using this observation together with \eqref{Theta_z(s)-inv},  \eqref{Theta_z}, and \eqref{Theta_z^2},  and the fact that $\ad(\mathsf{h})=3\ad(\caff_2)$, we conclude that  
	\begin{equation}
		\begin{gathered}\label{Theta-die}
			[\Thz{z;s}{}{},x\otimes 1]=0=[\Thz{z;s}{}{},1\otimes x]\\
			s[\EuScript{K}_z,\xi_{i,0}\otimes 1]=-[\Thz{z;s}{}{},\Delta^z(t_{i,1})]=-[\Thz{z;s}{}{},\square(t_{i,1})]
		\end{gathered}
	\end{equation}
	for all $x\in \gaff$ and $i\in \hbfI$. From the first equality with $x\in \haff$, we find that the component $\Thz{s}{}{\beta}$ is zero unless $\beta\in \cap_{i\in \hbfI}\alpha_i^\perp$. As $s[\EuScript{K}_z,\xi_{i,0}\otimes 1]$ has no such component and $\mathrm{ad}(\square (t_{i,1}))$ is weight zero, we deduce from the second line that 
	\begin{equation*}\label{Theta-die'}
		[\EuScript{K}_z,\xi_{i,0}\otimes 1]=0 \quad \forall\; i \in \hbfI.
	\end{equation*}

	To prove that $[\EuScript{K}_z,x\otimes 1]=0=[\EuScript{K}_z,1\otimes x]$ for general $x\in \gaff$, it now suffices to show that $\EuScript{K}_z$ commutes with $\Delta^z(x_{i,0}^\pm)$ for each $i\in \hbfI$. We have 
	\begin{equation*}
		3[[\caff_2\otimes 1, \Omega_z],\Delta^z(x_{i,0}^{\pm})]=\pm \frac{\hbar^2}{4}\mu_i [x_{i,0}^\pm \otimes 1, \Omega_z]+3[\caff_2\otimes 1, [\Omega_z,\square^z(x_{i,0}^\pm)]].
	\end{equation*}
	By Lemma 4.2 and \S6.3 of \cite{guay-nakajima-wendlandt}, we have $[\Omega_z,\square^z(x_{i,0}^+)]=x_{i,0}^+\otimes \xi_{i,0}$ and $[\Omega_z,\square^z(x_{i,0}^-)]=-z^{-1}\xi_{i,0}\otimes x_{i,0}^-$. Hence, we obtain 
	\begin{align*}
		3[[\caff_2\otimes 1, \Omega_z],\Delta^z(x_{i,0}^{-})]&=-\frac{\hbar^2}{4}\mu_i [x_{i,0}^- \otimes 1, \Omega_z],\\
		3[[\caff_2\otimes 1, \Omega_z],\Delta^z(x_{i,0}^{+})]&=\frac{\hbar^2}{4}\mu_i [x_{i,0}^+ \otimes 1, \Omega_z]+\frac{\hbar^2}{4}\mu_i x_{i,0}^+\otimes \xi_{i,0}=-\frac{\hbar^2}{4}z\mu_i [1\otimes x_{i,0}^+, \Omega_z].
	\end{align*}
	Combining these calculations with Lemma \ref{L:Theta_z} yields $[\EuScript{K}_z,\Delta^z(x_{i,0}^{\pm})]=0$ for all $i\in \hbfI$, as desired. This completes the proof of Part \eqref{Theta:2}. \qedhere

\end{pf}

\subsection{Computing $\EuScript{K}_{n\delta}$ explicitly, I}
\label{ssec:T-Gamma_i}

By the results of the previous section, $\EuScript{K}_z$ is a series of the form $\sum_{n>0}\EuScript{K}_{n\delta}z^{n\mathsf{ht}(\delta)}$ with each coefficient $\EuScript{K}_{n\delta}$ belonging to the centralizer of $U(\gaff)^{\otimes 2}$. The goal of this section and Section \ref{ssec:T-Gamma_i-2} is to  compute these coefficients explicitly. Our starting point is the following lemma, which will also play a crucial role in Section \ref{ssec:Q_z-exp}.
\begin{lem}\label{L:Gamma-i}
	Let $n$ be a positive integer and set $\beta=n\delta$. Then, for each $i\in \hbfI$, one has 
	\begin{equation*}
		[\square(t_{i,1}),\Theta_\beta]=d_i\hbar \sum_{\alpha+\gamma=\beta}\alpha(\mathsf{h})\gamma(h_i)\left[\Omega_\alpha,\Omega_\gamma\right],
	\end{equation*}
	where the summation runs over all positive roots $\alpha,\gamma\in \raff_+$ such that $\alpha+\gamma=\beta$. Moreover, the element 
	$
	\Gamma_{i,\beta}\ceq  [\square(t_{i,1}),\EuScript{K}_\beta]
	$
	has the following properties:
	\begin{enumerate}[font=\upshape]\itemsep0.25cm
		\item \label{Gamma:1}For each $a,b\in \C$, one has $(\tau_a\otimes \tau_b)(\Gamma_{i,\beta})=\Gamma_{i,\beta}$.
		\item\label{Gamma:2} One has 
		\begin{align*}
			(\Delta^z\otimes \Id)(\Gamma_{i,\beta})&=\Gamma_{i,\beta}^{13}+\Gamma_{i,\beta}^{23}z^{-\mathsf{ht}(\beta)}+\Gamma_{i,\beta}^-(z),\\
			(\Id\otimes \Delta^z)(\Gamma_{i,\beta})&=\Gamma_{i,\beta}^{12}+\Gamma_{i,\beta}^{13}z^{\mathsf{ht}(\beta)}
			+\Gamma_{i,\beta}^+(z),
		\end{align*}
		where $\Gamma^-_{i,\beta}(z)$ and $\Gamma^+_{i,\beta}(z)$ are given explicitly by 
		\begin{equation*}
			\begin{aligned}
				\Gamma^-_{i,\beta}(z)&=d_i\hbar \sum_{\alpha+\gamma=\beta}(\alpha(\mathsf{h})\gamma(h_i)-\gamma(\mathsf{h})\alpha(h_i))[\Omega_\alpha^{13},\Omega_\gamma^{23}]z^{-\mathsf{ht}(\gamma)}\\[-1em]
				&\hspace{7em} + \hbar \beta(\mathsf{h}) \mathrm{ad}(\xi_{i,0}^{(1)})\left( [\Omega_z^{12},\Omega_{\beta}^{13}+\Omega_\beta^{23}z^{-\mathsf{ht}(\beta)}]\right), \\[3pt]
				\Gamma^+_{i,\beta}(z)&=d_i\hbar \sum_{\alpha+\gamma=\beta}\left(\alpha(\mathsf{h})\gamma(h_i)-\gamma(\mathsf{h})\alpha(h_i)\right) [\Omega_\alpha^{12},\Omega_\gamma^{13}]z^{\mathsf{ht}(\gamma)}\\[-1em]
				&\hspace{7em} +\hbar \beta(\mathsf{h})\mathrm{ad}(\xi_{i,0}^{(2)})\left([\Omega_z^{23},\Omega_\beta^{12}+\Omega_\beta^{13}z^{\mathsf{ht}(\beta)}]\right),
			\end{aligned}
		\end{equation*}
		where both summations are taken over all $\alpha,\gamma\in \raff_+$ such that $\alpha+\gamma=\beta$,
		$\xi_{i,0}^{(1)}=\xi_{i,0}\ten 1\ten 1$ and $\xi_{i,0}^{(2)}=1\ten \xi_{i,0}\ten 1$.
		
		\item\label{Gamma:3} The series $\Gamma_{i,\beta}^+(z)$ and $\Gamma_{i,\beta}^-(z)$ satisfy
		\begin{equation*}
			\Gamma_{i,\beta}^\pm(z)\in (\naff^-\otimes \naff^\pm \otimes \naff^+)[z^{\pm 1}]
		\end{equation*}
	\end{enumerate}
\end{lem}
\begin{pf}
	
	From the first identity of Lemma \ref{L:Theta_z}, we obtain 
	\begin{equation*}
		[\Thz{z}{}{},\square(t_{i,1})]+\hbar[\Thz{z}{}{},[\xi_{i,0}\otimes 1,\Omega_z]]=-[\xi_{i,0}\otimes 1,[\square(\caff_3),\Omega_z]].
	\end{equation*}
	As $\mathrm{ad}(\caff_3)$ is weight zero, the right-hand side has no $\Ykm_{-\beta}\otimes \Ykm_{\beta}$ component. Since $\mathrm{ad}(t_{i,1})$ is weight zero, this implies that 
	\begin{equation}\label{sqt_i,theta_b:1}
		[\square(t_{i,1}),\Theta_\beta]=\hbar\sum_{\alpha+\gamma=\beta}[\Theta_\alpha,[\xi_{i,0}\otimes 1,\Omega_\gamma]]
		=-d_i\hbar\sum_{\alpha+\gamma=\beta}\gamma(h_i)[\Theta_\alpha,\Omega_\gamma],
	\end{equation}
	where the summations are taken over all $\gamma\in\raff_+$ and $\alpha\in \Qaff_+\setminus\{0\}$ such that $\alpha+\gamma=\beta$. However, for the term $\gamma(h_i)[\Theta_\alpha,\Omega_\gamma]$ to provide a nonzero contribution, we must have $\gamma\notin \alpha_i^\perp$. As $\alpha+\gamma=\beta$ and $\beta\in \alpha_i^\perp$, this implies that $\alpha\notin \alpha_i^\perp$ as well. By Part \eqref{Theta:3} of Proposition \ref{P:Theta}, for any such $\alpha$ we have 
	\begin{equation*}
		\EuScript{K}_\alpha=\Theta_\alpha+\alpha(\mathsf{h})\Omega_\alpha=0.
	\end{equation*}
	The formula for $[\square(t_{i,1}),\Theta_\beta]$ stated in the lemma now follows from \eqref{sqt_i,theta_b:1} after replacing $\Theta_\alpha$ by $-\alpha(\mathsf{h})\Omega_\alpha$.
	
	We now turn to establishing Parts \eqref{Gamma:1}--\eqref{Gamma:3} of the lemma. 
	
	\begin{pf}[Proof of Part \eqref{Gamma:1}]\let\qed\relax
		From what has been proven above, we have 
		\begin{equation}\label{eq:Gammandelta}
			\Gamma_{i,\beta}=d_i\hbar \sum_{\alpha+\gamma=\beta}\alpha(\mathsf{h})\gamma(h_i)\left[\Omega_\alpha,\Omega_\gamma\right]+\beta(\mathsf{h})[\square(t_{i,1}),\Omega_\beta].
		\end{equation}
		It is therefore sufficient to prove that $[\square(t_{i1}),\Omega_\beta]$ is fixed by $\tau_a\otimes \tau_b$ for any $a,b\in \C$. We have 
		\begin{align*}
			(\tau_a\otimes \tau_b)([\square(t_{i,1}),\Omega_\beta])&=[\square(t_{i,1})+a\xi_{i,0}\otimes 1 +1\otimes b\xi_{i,0},\Omega_\beta]\\
			&=[\square(t_{i,1}),\Omega_\beta]+(\beta,\alpha_i)(b-a)\Omega_\beta.
		\end{align*}
		Since $(\beta,\alpha_i)=n(\delta,\alpha_i)=0$, this coincides with $[\square(t_{i,1}),\Omega_\beta]$. 
	\end{pf}
	
	\begin{pf}[Proof of Part \eqref{Gamma:2}]\let\qed\relax
		By the first assertion of the lemma, we have 
		\begin{align*}
			&(\Delta^z\otimes \Id)([\square(t_{i,1}),\Theta_\beta])\\
			&
			=d_i\hbar \!\sum_{\alpha+\gamma=\beta}\!\alpha(\mathsf{h})\gamma(h_i)\left[\Omega_\alpha^{13}+\Omega_\alpha^{23}z^{-\mathsf{ht}(\alpha)},\Omega_\gamma^{13}+\Omega_\gamma^{23}z^{-\mathsf{ht}(\gamma)}\right]\\
			&
			=(\square^z\otimes \Id)[\square(t_{i,1}),\Theta_\beta]\\
			&\hspace{5em}+d_i\hbar \!\sum_{\alpha+\gamma=\beta}\!\alpha(\mathsf{h})\gamma(h_i)
			([\Omega_\alpha^{13},\Omega_\gamma^{23}]z^{-\mathsf{ht}(\gamma)}
			+[\Omega_\alpha^{23},\Omega_\gamma^{13}]z^{-\mathsf{ht}(\alpha)})\\
			&
			=(\square^z\otimes \Id)[\square(t_{i,1}),\Theta_\beta]+d_i\hbar \!\sum_{\alpha+\gamma=\beta}\!(\alpha(\mathsf{h})\gamma(h_i)-\gamma(\mathsf{h})\alpha(h_i))[\Omega_\alpha^{13},\Omega_\gamma^{23}]z^{-\mathsf{ht}(\gamma)}.
		\end{align*}
		Similarly, we have 
		\begin{align*}
			&(\Delta^z\otimes \Id)[\square(t_{i,1}),\Omega_\beta]\\
			&=
			[t_{i,1}^{(3)}+\square(t_{i,1})\otimes 1+\hbar[\xi_{i,0}\otimes 1,\Omega_z]\otimes 1, \Omega_{\beta}^{13}+\Omega_\beta^{23}z^{-\mathsf{ht}(\beta)}]\\
			&
			=(\square^z\otimes \Id)[\square(t_{i,1}),\Omega_\beta]+\hbar [\mathrm{ad}(\xi_{i,0}^{(1)})(\Omega_z^{12}),\Omega_{\beta}^{13}+\Omega_\beta^{23}z^{-\mathsf{ht}(\beta)}]\\
			&
			=(\square^z\otimes \Id)[\square(t_{i,1}),\Omega_\beta]+\hbar \mathrm{ad}(\xi_{i,0}^{(1)})( [\Omega_z^{12},\Omega_{\beta}^{13}+\Omega_\beta^{23}z^{-\mathsf{ht}(\beta)}]),
		\end{align*}
		where $t_{i,1}^{(3)}=1\ten 1\ten t_{i,1}$ and in the last line we have used that $\mathrm{ad}(\xi_{i,0}^{(1)})$ commutes with $\Omega_{\beta}^{13}$ and $\Omega_{\beta}^{23}$ as $\beta=n\delta$. As $\Gamma_{i,\beta}=[\square(t_{i,1}),\Theta_\beta+\beta(\mathsf{h})\Omega_\beta]$, this implies that 
		\begin{equation*}
			(\Delta^z\otimes \Id)(\Gamma_{i,\beta})=\Gamma_{i,\beta}^{13}+\Gamma_{i,\beta}^{23}z^{-\mathsf{ht}(\beta)}+\Gamma_{i,\beta}^-(z)
		\end{equation*}
		with $\Gamma_{i,\beta}^-(z)$ as in the statement of the lemma. The computation of $(\Id\otimes \Delta^z)(\Gamma_{i,\beta})$ is nearly identical, and hence omitted.  \qedhere
	\end{pf}
	
	\begin{pf}[Proof of Part \eqref{Gamma:3}]
		As the proofs for $\Gamma_{i,\beta}^-(z)$ and $\Gamma_{i,\beta}^+(z)$ are the same, we will focus on the former. By the formulas of Part \eqref{Gamma:2}, it is sufficient to prove that 
		\begin{gather*}
			\mathrm{ad}(\xi_{i,0}^{(1)})\left( [\Omega_z^{12},\Omega_{\beta}^{13}+\Omega_\beta^{23}z^{-\mathsf{ht}(\beta)}]\right)
			\in (\naff^-\otimes \naff^-\otimes \naff^+)[z^{-1}].
		\end{gather*}
		It is clear that $\mathrm{ad}(\xi_{i,0}^{(1)})\left( [\Omega_z^{12},\Omega_{\beta}^{13}+\Omega_\beta^{23}z^{-\mathsf{ht}(\beta)}]\right)$ is a Laurent series in $z$ with coefficients in $\naff^-\otimes \gaff \otimes \naff^+$. Let $\widetilde{\Omega}_z\ceq  \Omega_z+(\Omega^{-}_{z^{-1}})^{21} \in \gkm^{\otimes 2}[\![z^{\pm 1}]\!]$. Note that this is just $(\Id\otimes \upsigma_z)(\widetilde{\Omega})$ where $\widetilde{\Omega}$ is the \textit{full} Casimir tensor of $\gkm$. It satisfies
		\begin{equation*}
			[\widetilde{\Omega}_z,\square^z(x)]=0\quad \forall\; x\in \gkm.
		\end{equation*}
		Hence, we have 
		\begin{equation*}
			[\Omega_z^{12},\Omega_{\beta}^{13}+\Omega_\beta^{23}z^{-\mathsf{ht}(\beta)}]
			=
			-[(\Omega^{-}_{z^{-1}})^{21}, \Omega_{\beta}^{13}+\Omega_\beta^{23}z^{-\mathsf{ht}(\beta)}],
		\end{equation*}
		whose right-hand side is a formal series in $z^{-1}$ with coefficients in $\gaff\otimes \naff^-\otimes \naff^+$. 
		Hence, we have 
		\begin{equation*}
			\mathrm{ad}(\xi_{i,0}^{(1)})\left( [\Omega_z^{12},\Omega_{\beta}^{13}+\Omega_\beta^{23}z^{-\mathsf{ht}(\beta)}]\right)\in (\naff^-\otimes \gaff \otimes \naff^+)[z^{-1};z]\!]\cap (\gaff\otimes \naff^-\otimes \naff^+)[\![z^{-1}]\!]
		\end{equation*}
		This shows that  $\mathrm{ad}(\xi_{i,0}^{(1)})\left( [\Omega_z^{12},\Omega_{\beta}^{13}+\Omega_\beta^{23}z^{-\mathsf{ht}(\beta)}]\right)$ is a polynomial in $z^{-1}$ with coefficients in $\naff^-\otimes \naff^-\otimes \naff^+$, and thus that $\Gamma_{i,\beta}^-(z)\in (\naff^-\otimes \naff^-\otimes \naff^+)[z^{-1}]$. \qedhere
	\end{pf}\let\qed\relax
\end{pf}

\subsection{Computing $\EuScript{K}_{n\delta}$ explicitly, II}
\label{ssec:T-Gamma_i-2}

We now turn towards computing each coefficient $\EuScript{K}_{n\delta}$ of the series $\EuScript{K}_z$. 
For the sake of brevity, let us define a linear operator $\mathrm{ad}_{x,y}$ on $\Ykm$, for each $x,y\in \Ykm$, by setting 
\begin{equation*} 
	\mathrm{ad}_{x,y}\ceq  \mathrm{ad}(x)\circ \mathrm{ad}(y):\Ykm\to \Ykm.
\end{equation*} 
For each $n>0$ and $\ell\in \{1,2\}$, we then define 
\begin{gather*}
	\EuScript{K}_{i,n\delta;\ell}
	\ceq  
	(\mathrm{ad}^{(\ell)}_{x_{i,1}^-,x_{i,0}^+}+\mathrm{ad}^{(\ell)}_{x_{i,1}^+,x_{i,0}^-})(\Gamma_{i,n\delta})
	-\hbar \mathrm{ad}^{(\ell)}_{x_{i,0}^-,x_{i,0}^+}(\Gamma_{i,n\delta}) \cdot \xi_{i,0}^{(\ell)},
\end{gather*}
where the superscript in $\mathrm{ad}^{(\ell)}_{x,y}$ indicates that it operates in the $\ell$-th tensor factor. 
\begin{prop}\label{P:K-exp}
	Let $n$ be a positive integer and fix $\ell\in \{1,2\}$. Then $\EuScript{K}_{n\delta}$ is given explicitly by
	\begin{equation*}
		\EuScript{K}_{n\delta}=(-1)^\ell\frac{3}{2{n\delta}(\mathsf{h})}\sum_{i\in \hbfI} \frac{a_i}{d_i}\EuScript{K}_{i,{n\delta};\ell}.
	\end{equation*}
\end{prop}
\begin{pf}
	We shall establish the claimed formula for $\EuScript{K}_{n\delta}$ in the case  where $\ell=1$. The proof in the $\ell=2$ case follows by a simple modification of the same argument. To begin, note that since $[3\caff_2\otimes 1,\EuScript{K}_{n\delta}]=-n\delta(\mathsf{h})\EuScript{K}_{n\delta}$ (and $\delta(\mathsf{h})\neq 0$), it is sufficient to establish that 
	\begin{equation}\label{[C_2,K_beta]}
		[\caff_2\otimes 1,\EuScript{K}_{n\delta}]=\sum_{i\in \hbfI} \frac{a_i}{2d_i}\EuScript{K}_{i,{n\delta};1}.
	\end{equation}
	By Part \eqref{Theta:2} of Proposition \ref{P:Theta}, we have 
	\begin{equation*}
		[x_{i,0}^\pm \otimes 1,\Gamma_{i,{n\delta}}]=[x_{i,0}^\pm \otimes 1,[\square(t_{i,1}),\EuScript{K}_{n\delta}]]=\mp 2d_i[x_{i,1}^\pm \otimes 1, \EuScript{K}_{n\delta}],
	\end{equation*}
	and therefore 
	\begin{equation*}
		[x_{i,1}^\pm \otimes 1, \EuScript{K}_{n\delta}]=\mp\frac{1}{2d_i}[x_{i,0}^\pm \otimes 1,\Gamma_{i,{n\delta}}]
		\quad \forall \; i\in \hbfI.
	\end{equation*}
	It follows that, for each $i\in \hbfI$, we have 
	\begin{equation*}
		\begin{aligned}
			[\xi_{i,1}\otimes 1,\EuScript{K}_{n\delta}]
			&=
			[[x_{i,1}^+\otimes 1,x_{i,0}^-\otimes 1],\EuScript{K}_{n\delta}]
			\\
			&=
			-\frac{1}{2d_i}[[x_{i,0}^+\otimes 1,\Gamma_{i,{n\delta}}], x_{i,0}^-\otimes 1]
			=\frac{1}{2d_i}\mathrm{ad}_{x_{i,0}^-,x_{i,0}^+}^{(1)}(\Gamma_{i,{n\delta}})
		\end{aligned}
	\end{equation*}
	where in the second equality we have used that, by Proposition \ref{P:Theta}, $[x_{i,0}^-\otimes 1, \EuScript{K}_{n\delta}]=0$. Similarly, we obtain 
	\begin{equation*}
		\begin{aligned}
			[\xi_{i,2}\otimes 1,\EuScript{K}_{n\delta}]
			&=
			[[x_{i,1}^+\otimes 1,x_{i,1}^-\otimes 1],\EuScript{K}_{n\delta}]
			\\
			&
			=-[x_{i,1}^-\otimes 1,[x_{i,1}^+\otimes 1,\EuScript{K}_{n\delta}]]+[x_{i,1}^+\otimes 1,[x_{i,1}^-\otimes 1,\EuScript{K}_{n\delta}]]\\
			&
			=\frac{1}{2d_i}(\mathrm{ad}_{x_{i,1}^-,x_{i,0}^+}^{(1)}+\mathrm{ad}_{x_{i,1}^+,x_{i,0}^-}^{(1)})(\Gamma_{i,{n\delta}}).
		\end{aligned}
	\end{equation*}
	Since $t_{i,2}=\xi_{i,2}-\hbar \xi_{i,1}\xi_{i,0}+\frac{\hbar^2}{3}\xi_{i,0}^3$ (see \eqref{eq:ti2}) and $[\xi_{i,0}\otimes 1,\EuScript{K}_{n\delta}]=0$, the above computations yield 
	\begin{equation*}
		\begin{aligned}
			2d_i[t_{i,2}\otimes 1,\EuScript{K}_{n\delta}] &= (\mathrm{ad}_{x_{i,1}^-,x_{i,0}^+}^{(1)}+\mathrm{ad}_{x_{i,1}^+,x_{i,0}^-}^{(1)})(\Gamma_{i,{n\delta}})-\hbar\mathrm{ad}_{x_{i,0}^-,x_{i,0}^+}^{(1)}(\Gamma_{i,{n\delta}})\cdot \xi_{i,0}^{(1)} \\
			&=\EuScript{K}_{i,{n\delta};1}.
		\end{aligned}
	\end{equation*}
	As $\caff_2=\sum_{i\in \hbfI} a_i t_{i,2}$, this implies the formula \eqref{[C_2,K_beta]}. \qedhere
\end{pf}

\subsection{The series $\mathrm{Q}_z$}\label{ssec:Q_z-exp}

In this section, we use the results of the previous subsections to construct a series $\mathrm{Q}_z$ satisfying all the properties listed in Part \eqref{map-T:3} of Theorem \ref{T:map-T}. To begin, recall from Proposition \ref{P:K-exp} that for each positive integer $n$, index $i\in \hbfI$, and number $\ell\in \{1,2\}$, we defined 
\begin{equation*}
	\EuScript{K}_{i,n\delta;\ell}
	\ceq  
	(\mathrm{ad}^{(\ell)}_{x_{i,1}^-,x_{i,0}^+}+\mathrm{ad}^{(\ell)}_{x_{i,1}^+,x_{i,0}^-})(\Gamma_{i,n\delta})
	-\hbar \mathrm{ad}^{(\ell)}_{x_{i,0}^-,x_{i,0}^+}(\Gamma_{i,n\delta}) \cdot \xi_{i,0}^{(\ell)},
\end{equation*}
where $\xi_{i,0}^{(\ell)}=1^{\otimes (\ell-1)}\otimes \xi_{i,0}\otimes 1^{\otimes (2-\ell)}$ and $\Gamma_{i,n\delta}$ and  $\mathrm{ad}^{(\ell)}_{x,y}$ are given as follows:
\begin{enumerate}[label=\emph{\roman*})]\itemsep0.25cm
	\item $\Gamma_{i,n\delta}=[\square(t_{i,1}),\EuScript{K}_{n\delta}]$. By \eqref{eq:Gammandelta} in the proof of Lemma \ref{L:Gamma-i}, this  may be written equivalently as 
	\begin{equation*}
		\Gamma_{i,n\delta}=d_i\hbar \sum_{\alpha+\gamma=n\delta}\alpha(\mathsf{h})\gamma(h_i)\left[\Omega_\alpha,\Omega_\gamma\right]+n\delta(\mathsf{h})[\square(t_{i,1}),\Omega_{n\delta}],
	\end{equation*}
	with 
	$\mathsf{h}=\frac{\hbar^2}{\zeta}(d-\sum_{i\in \hbfI}\zeta_id_i h_i )\in \hkm$, as in \eqref{3C_2-preimage}. 
	
	\item $\mathrm{ad}_{x,y}\ceq  \mathrm{ad}(x)\circ \mathrm{ad}(y)$ for all $x,y\in \Yhg$, while 
	\begin{equation*}
		\mathrm{ad}_{x,y}^{(1)}=\mathrm{ad}_{x,y}\otimes \Id \quad \text{ and }\quad \mathrm{ad}_{x,y}^{(2)}=\Id\otimes \mathrm{ad}_{x,y}.
	\end{equation*}
\end{enumerate}

We now introduce the distinguished series $\mathrm{Q}_z\in \Ykm^{\otimes 2}[\![z]\!]$ by setting
\begin{equation*}
	\mathrm{Q}_z\ceq  -\sum_{n>0} \frac{1}{n\delta(\mathsf{h})}\EuScript{K}_{n\delta}z^{n\mathsf{ht}(\delta)}.
\end{equation*}
By Proposition \ref{P:K-exp}, the coefficient $\mathrm{Q}_{n\delta}\ceq  \frac{1}{n\delta(\mathsf{h})}\EuScript{K}_{n\delta}$ is given by 
\begin{equation}\label{Q:2-defs}
	\mathrm{Q}_{n\delta}=\frac{3}{2n^2\delta(\mathsf{h})^2}\sum_{i\in \hbfI} \frac{a_i}{d_i}\EuScript{K}_{i,n\delta;1}
	=
	-\frac{3}{2n^2\delta(\mathsf{h})^2}\sum_{i\in \hbfI} \frac{a_i}{d_i}\EuScript{K}_{i,n\delta;2}.
\end{equation}
The following proposition completes the proof of Theorem \ref{T:map-T} by showing that $\mathrm{Q}_z$ satisfies all the conditions spelled out in Part \eqref{map-T:3} therein.
\begin{prop}\label{P:Q_z}
	The series $\mathrm{Q}_z$ has the following properties:
	\begin{enumerate}[font=\upshape]\itemsep0.25cm
		\item\label{Qz:1} For each $h\in \hkm$, one has 
		\begin{equation*}
			\Delta^z(\Top(h))=\Top(h)\otimes 1 + 1\otimes \Top(h) + \hbar[h\otimes 1, \Omega_z+\mathrm{Q}_z].
		\end{equation*}
		
		\item\label{Qz:2} For each $n>0$, one has 
		\begin{equation*}
			\mathrm{Q}_{n\delta}\in \Ysupkm{-}_{-n\delta}\otimes \Ysupkm{+}_{n\delta}. 
		\end{equation*}
		\item\label{Qz:3} For each $a,b\in \C$, one has 
		\begin{equation*}
			(\tau_a\otimes \tau_b)(\mathrm{Q}_z)=\mathrm{Q}_{z}.
		\end{equation*}
		\item\label{Qz:4} The tensor factors of $\mathrm{Q}_z$ are primitive:
		\begin{equation*}
			(\Delta^w\otimes \Id)(\mathrm{Q}_{z})=\mathrm{Q}_{z}^{13}+\mathrm{Q}_{z/w}^{23} \quad \text{ and }\quad 
			(\Id\otimes \Delta^w)(\mathrm{Q}_z)=\mathrm{Q}_z^{12}+\mathrm{Q}_{zw}^{13}.
		\end{equation*}
		\item\label{Qz:5}  $\mathrm{Q}_{z}$ belongs to the centralizer of $U(\gaff)^{\otimes 2}$ in $\Ykm^{\otimes 2}[\![z]\!]$.
		
	\end{enumerate}
\end{prop}
\begin{pf} Note that Part \eqref{Qz:5} follows immediately from the definition of $\mathrm{Q}_z$ and Part \eqref{Theta:2} of Proposition \ref{P:Theta}. We shall prove the remaining four statements in order. 
	\begin{pf}[Proof of \eqref{Qz:1}]\let\qed\relax
		If $h\in \haff$, then $[h\otimes 1,\mathrm{Q}_z]=0$ and so the claimed formula for $\Delta^z(\Top(h))$ holds by definition of $\Delta^z$ (see \eqref{Delta^z:def}). It thus suffices to verify the formula in the case where $h=\mathsf{h}$. Since $\Top(\mathsf{h})=\caff_3$, this amounts to checking
		\begin{equation*}
			\Delta^z(\caff_3)=\caff_3\otimes 1 + 1\otimes \caff_3+\hbar[\mathsf{h}\otimes 1,\Omega_z+\mathrm{Q}_z].
		\end{equation*}
		This is a consequence of the definition of $\mathrm{Q}_z$. Indeed, we have 
		\begin{equation*}
			[\mathsf{h}\otimes 1,\mathrm{Q}_z]=-\sum_{n>0}n\delta(\mathsf{h})\mathrm{Q}_{n\delta}z^{n\mathsf{ht}(\delta)}=\EuScript{K}_z=\Thz{z}{}{}-[\mathsf{h}\otimes 1,\Omega_z]. \qedhere
		\end{equation*}
	\end{pf}
	\begin{pf}[Proof of \eqref{Qz:2}]
		By Lemma \ref{L:Gamma-i}, we have 
		\begin{equation*}
			\Gamma_{i,n\delta}=d_i\hbar \sum_{\alpha+\gamma=n\delta}\alpha(\mathsf{h})\gamma(h_i)\left[\Omega_\alpha,\Omega_\gamma\right]+n\delta(\mathsf{h})[\square(t_{i,1}),\Omega_{n\delta}]
			\in\Ysupkm{-}_{-n\delta}\otimes \Ysupkm{+}_{n\delta}. 
		\end{equation*}
		It follows by definition of $\EuScript{K}_{i,n\delta;1}$ and $\EuScript{K}_{i,n\delta;2}$ that 
		\begin{equation*}
			\EuScript{K}_{i,n\delta;1}\in \Ykm_{-n\delta}\otimes \Ysupkm{+}_{n\delta}
			\quad \text{ and }\quad 
			\EuScript{K}_{i,n\delta;2}\in \Ysupkm{-}_{-n\delta}\otimes \Ykm_{n\delta}. 
		\end{equation*}
		Thus, we conclude from \eqref{Q:2-defs} that $\mathrm{Q}_{n\delta}$
		belongs to $\Ysupkm{-}_{-n\delta}\otimes \Ysupkm{+}_{n\delta}$. \let\qed\relax
	\end{pf}
	
	\begin{pf}[Proof of \eqref{Qz:3}]\let\qed\relax
		Fix a positive integer $n$ and complex numbers $a,b\in \C$. Then, by definition of $\EuScript{K}_{i,n\delta;1}$ and $\EuScript{K}_{i,n\delta;2}$, we have 
		\begin{equation*}
			\begin{aligned}
				(\Id&\otimes \tau_b)(\EuScript{K}_{i,n\delta;1})\\
				&=(\mathrm{ad}^{(1)}_{x_{i,1}^-,x_{i,0}^+}+\mathrm{ad}^{(1)}_{x_{i,1}^+,x_{i,0}^-})((\Id\otimes \tau_b)\Gamma_{i,{n\delta}})
				-\hbar \mathrm{ad}^{(1)}_{x_{i,0}^-,x_{i,0}^+}((\Id\otimes \tau_b)\Gamma_{i,{n\delta}}) \cdot \xi_{i,0}^{(1)},
				\\[5pt]
				(\tau_a&\otimes \Id)(\EuScript{K}_{i,n\delta;2})\\
				&
				=(\mathrm{ad}^{(2)}_{x_{i,1}^-,x_{i,0}^+}+\mathrm{ad}^{(2)}_{x_{i,1}^+,x_{i,0}^-})((\tau_a\otimes \Id)\Gamma_{i,{n\delta}})
				-\hbar \mathrm{ad}^{(2)}_{x_{i,0}^-,x_{i,0}^+}((\tau_a\otimes \Id)\Gamma_{i,{n\delta}}) \cdot \xi_{i,0}^{(2)}.
			\end{aligned}
		\end{equation*}
		Since $(\tau_a\otimes \tau_b)(\Gamma_{i,{n\delta}})=\Gamma_{i,{n\delta}}$ by 
		Part \eqref{Gamma:1} of Lemma \ref{L:Gamma-i}, we get that $(\Id\otimes \tau_b)(\EuScript{K}_{i,n\delta;1})=\EuScript{K}_{i,n\delta;1}$ and $(\tau_a\otimes \Id)(\EuScript{K}_{i,n\delta;2})=\EuScript{K}_{i,n\delta;2}$. By \eqref{Q:2-defs}, this implies that 
		\begin{equation*}
			(\tau_a\otimes \tau_b)(\mathrm{Q}_z)=\mathrm{Q}_z \quad \forall\; a,b\in \C. \qedhere
		\end{equation*}
	\end{pf}
	
	\begin{pf}[Proof of \eqref{Qz:4}]
		By Part \eqref{Gamma:2} of Lemma \ref{L:Gamma-i}, we have 
		\begin{align*}
			&((\Delta^z-\square^z)\otimes \Id)(\EuScript{K}_{i,{n\delta};2})\\
			&=(\mathrm{ad}^{(3)}_{x_{i,1}^-,x_{i,0}^+}+\mathrm{ad}^{(3)}_{x_{i,1}^+,x_{i,0}^-})(\Gamma_{i,{n\delta}}^-(z))
			-\hbar \mathrm{ad}^{(3)}_{x_{i,0}^-,x_{i,0}^+}(\Gamma_{i,{n\delta}}^-(z)) \cdot \xi_{i,0}^{(3)}.
		\end{align*}
		By Part \eqref{Gamma:3} of Lemma \ref{L:Gamma-i} and Proposition \ref{P:K-exp}, this yields that 
		\begin{equation}\label{K^-:lives}
			\EuScript{K}_{n\delta}^-(z)\ceq  ((\Delta^z-\square^z)\otimes \Id)(\EuScript{K}_{n\delta})\in (\naff^-\otimes \naff^-\otimes \Yhg)[z^{-1}]
		\end{equation}
		Similarly, Lemma \ref{L:Gamma-i} and Proposition \ref{P:K-exp} imply that 
		\begin{equation}\label{K^+:lives}
			\EuScript{K}_{n\delta}^+(z)\ceq  (\Id \otimes(\Delta^z-\square^z))(\EuScript{K}_{n\delta})\in (\Ykm\otimes \naff^+\otimes \naff^+)[z]. 
		\end{equation}
		Note that, by definition of $\mathrm{Q}_z$, the statement of Part \eqref{Qz:4} is equivalent to the assertion that $\EuScript{K}_{n\delta}^+(z)=0=\EuScript{K}_{n\delta}^-(z)$. This is established in the following claim. 
		
		\noindent \textit{Claim}. $\EuScript{K}_{n\delta}^+(z)$ and $\EuScript{K}_{n\delta}^-(z)$ are both equal to $0$. 
		\begin{pf}[Proof of Claim]
			Define $\Xi^{\pm}(z,w)$ by the equations 
			\begin{align*}
				(\Delta^z\otimes \Id)(\Thz{zw}{}{})&=\Thz{zw}{}{}^{13}+\Thz{w}{}{}^{23}+\Xi^-(z,w),\\  
				(\Id\otimes \Delta^w)(\Theta_z)&=\Thz{z}{}{}^{12}+\Thz{zw}{}{}^{13}+\Xi^+(z,w).
			\end{align*}
			By the twisted coassociativity of $\Delta^z$ (see \eqref{twisted-coalg}), we have 
			\begin{equation*}
				\begin{aligned}
					0=&\hbar^{-1}((\Delta^z\otimes \Id)\circ \Delta^{zw}-(\Id\otimes \Delta^w)\circ \Delta^z)(\caff_3)\\
					= &\Thz{z}{}{}^{12}-\Thz{w}{}{}^{23}+ (\Delta^z\otimes \Id)(\Thz{zw}{}{})-(\Id\otimes \Delta^w)(\Thz{z}{}{})\\
					=&\Xi^-(z,w)-\Xi^+(z,w).
				\end{aligned}
			\end{equation*}
			Hence, we have $\Xi^+(z,w)=\Xi^-(z,w)$. Therefore, we drop the superscript and write $\Xi(z,w)$ for $\Xi^{\pm}(z,w)$. Next, observe that, since the coefficients of $[\mathsf{h}\otimes 1,\Omega_z]$ belong to $\gkm\otimes \gkm$, the element $\EuScript{K}_z=\Thz{z}{}{}-[\mathsf{h}\otimes 1,\Omega_z]$ satisfies
			\begin{align*}
				(\Delta^z\otimes \Id)(\EuScript{K}_{zw})&=\EuScript{K}_{zw}^{13}+\EuScript{K}_{w}^{23}+\Xi(z,w)\,,\\
				(\Id\otimes \Delta^w)(\EuScript{K}_z)&=\EuScript{K}_z^{12}+\EuScript{K}_{zw}^{13}+\Xi(z,w).
			\end{align*}
			In particular, we have 
			\begin{equation*}
				\sum_{n>0}\EuScript{K}_{n\delta}^-(z)(zw)^{n\mathsf{ht}(\delta)}=\Xi(z,w)=
				\sum_{n>0}\EuScript{K}_{n\delta}^+(w)z^{n\mathsf{ht}(\delta)}.
			\end{equation*}
			However, by \eqref{K^-:lives} and \eqref{K^+:lives}, this is only possible if both the left-hand side and right-hand side are identically zero. Therefore, $\EuScript{K}_{n\delta}^-(z)=0=\EuScript{K}_{n\delta}^+(w)$ for each positive integer $n$, as desired. This completes the proof of Part \eqref{Qz:4}, and thus the proof of the proposition. \qedhere
		\end{pf}
		\let\qed\relax
	\end{pf}
	\let\qed\relax
\end{pf}

\section{Construction of $\RR^-(s)$}\label{sec:neg-R}

In this section, we show that the standard and Drinfeld tensor products on $\OY$ are
conjugate to each other by a triangular element $\RR^-(s)$.
In other words, $\RR^-(s)$ is a rational tensor structure on
the identity functor 
\[
(\Id,\RR^-(s)) : \lp\OY, \dtensor{s}\rp \to \lp\OY, \kmtensor{s}\rp.
\]
We refer the reader to Theorem \ref{thm:R^--reps} below for
the precise meaning of this statement.
Our proof is almost verbatim to the one for finite type
\cite[Thm.~4.1]{GTLW}, with some subtle differences highlighted
in this section. Additionally, the argument in {\em loc. cit.} rests on the
linear map $\Top:\hkm\to\Ysupkm{0}$ and the formulae for
$\Delta_s(\Top(h))$. For affine Yangians, the map only
exists, {\em a priori}, from $\haff\to\Ysupkm{0}$. That is
the reason we had to extend it, using $\caff_3$ as in
Theorem \ref{T:map-T}.

\subsection{The formal series $\RR^-(s,z)$}\label{ssec:R^--formal} 
Recall from Section \ref{ssec:T-main} that there is a series $\mathrm{Q}_z=\sum_{n>0}\mathrm{Q}_{n\delta}z^{n\mathsf{ht}(\delta)}$, with $\mathrm{Q}_{n\delta}$ defined explicitly in \eqref{Q:2-defs}, which satisfies
\begin{equation*}
\Delta^z(\Top(h))=\Top(h)\otimes 1 +1\otimes \Top(h)+\hbar[h\otimes 1,\Omega_z+\mathrm{Q}_z] \quad \forall\; h\in \hkm
\end{equation*}
in addition to the conditions of Part \eqref{map-T:3} in Theorem \ref{T:map-T}. Here the transformation $\Top$ itself is defined in Theorem \ref{T:map-T} above, while $\Omega_z\in (\gkm\otimes \gkm)[\![z]\!]$ is as in Section \ref{ssec:pre-cop}. Furthermore, we introduce the algebra homomorphisms
\begin{equation*}
\ddelta{s}^z \ceq   (\Id\otimes\upsigma_z)\circ\ddelta{s} \quad \text{ and }\quad \Delta_s^z\ceq  (\tau_s\otimes \Id)\circ \Delta^z
\end{equation*}
in addition to the function $\nu:\Qaff_+\to \Z_{\geqslant 0}$ defined by 
\begin{equation*}
\nu(\beta)=\text{min}\{k\in\N : \beta=\alpha^{(1)}+\cdots+\alpha^{(k)} \text{ for } \alpha^{(1)},\ldots,\alpha^{(k)}\in \raff_+\}.
\end{equation*}
The following theorem provides the formal version of the main
result of this section.
\begin{thm}\label{thm:R-}
There exists a unique family of elements $\{\RR^-(s)_\beta\}_{\beta\in \Qaff_+}$, with 
 $\RR^-(s)_\beta\in (\Ykm_{-\beta}\otimes \Ykm_\beta)[\![s^{-1}]\!]$ for each $\beta\in \Qaff_+$, 
 satisfying $\RR^-(s)_0=1\otimes 1$ in addition to the intertwiner equation 
 \begin{equation}\label{eq:IER-}
 \RR^-(s,z)\Delta_s^z(\Top(h))= \ddelta{s}^{\!z}(\Top(h)) \RR^-(s,z) \quad \forall\; h\in \hkm
 \end{equation}
 in $\Ykm^{\otimes 2}[s;s^{-1}]\!] [\![z]\!]$, where $\RR^-(s,z)=\sum_{\beta\in \Qaff_+}\RR^-(s)_\beta z^{\mathsf{ht}(\beta)}$. Moreover: 
\begin{enumerate}[font=\upshape]\itemsep0.25cm
\item\label{R-(s,z):1} For each $x\in \Ykm$, one has 
\begin{equation*}
\RR^-(s,z)\Delta_s^z(x)= \ddelta{s}^{\!z}(x) \RR^-(s,z).
\end{equation*}
\item\label{R-(s,z):2}  For each $\beta\in \Qaff_+$, $\RR^-(s)_\beta$ satisfies 
\begin{equation*}
\RR^-(s)_\beta\in s^{-\nu(\beta)}\left( \Ysupkm{-}\otimes \Ysupkm{+}\right)[\![s^{-1}]\!].
\end{equation*}
\item\label{R-(s,z):3}  For each $\beta\in \Qaff_+$ and $a,b\in \C$, $\RR^-(s)_\beta$ satisfies 
\begin{equation*}
\tau_a\otimes \tau_b(\RR^-(s)_\beta)=\RR^-(s+a-b)_\beta.
\end{equation*}
\item\label{R-(s,z):4}  As an element of $\Ykm^{\otimes 2}[\![z]\!][\![s^{-1}]\!]$, $\RR^-(s,z)$ is of the form
\begin{equation*}
\RR^-(s,z)=1+\hbar s^{-1}(\Omega_z^-+\mathrm{Q}_z)+O(s^{-2}).
\end{equation*}
\end{enumerate}

\end{thm}
\begin{pf}
With Theorem \ref{T:map-T} at our disposal, the proof of the finite-type counterpart to this result, given in \cite[\S4]{GTLW},
carries over to the present setting with only minor adjustments. 
In what follows we summarize the key steps while making transparent the role played by $z$ (which does not appear in \cite{GTLW}) and Theorem \ref{T:map-T}. 

\begin{pf}[Proof of Existence and Uniqueness] By Theorem \ref{T:map-T},  $(\tau_s\otimes \Id)(\mathrm{Q}_z)=\mathrm{Q}_z$ and $\mathrm{ad}(\tau_s(\Top(h)))=\mathrm{ad}(\Top(h)+sh)$ for all $h\in \hkm$. It follows that the intertwiner equation \eqref{eq:IER-} is equivalent to 
\begin{equation*}
(\adT(h)+\mathrm{ad}(sh\otimes \Id))\cdot \RR^-(s,z)=\hbar \RR^-(s,z)[h\otimes 1,\Omega_z+\mathrm{Q}_z]
\end{equation*}
for all $h\in \hkm$, where we have set $\adT(h)=\mathrm{ad}(\square(\Top(h)))$. 
Taking the $\Ykm_{-\beta}\otimes \Ykm_\beta$ component of this equation, for any fixed $\beta\in \Qaff_+$, and dividing both sides of the resulting identity by $z^{\mathsf{ht}(\beta)}$ yields 
\begin{equation}\label{R^-:recur-rat}
(\adT(h)-s\beta(h)) \cdot \RR^-(s)_\beta
= -\hbar\sum_{\alpha\in\Phi_+} \alpha(h) \RR^-(s)_{\beta-\alpha}
(\Omega_\alpha+\mathrm{Q}_\alpha),
\end{equation}
where it is understood that $\RR^-(s)_\gamma=0=\mathrm{Q}_\alpha$ for $\gamma\not\in\Qaff_+$ and $\alpha \notin \Z\delta$. 
In particular, the right--hand side of the equation above is a finite sum. The proof that there is at most one family $\{\RR^-(s)_\beta\}_{\beta\in \Qaff_+}$ satisfying the conditions of the theorem now proceeds identically to the proof of the analogous assertion for the Yangian of a finite--dimensional simple Lie algebra given in \cite[\S4.2]{GTLW}. Indeed, if $h\in \hkm$ is chosen so that $\beta(h)\neq 0$, then the operator on the left--hand side of the above equation can be inverted to yield 
\begin{equation}\label{R^-:recur}
\begin{aligned}
\RR^-(s)_\beta &= \frac{\hbar}{s\beta(h)} \lp 1 - \frac{\adT(h)}{s\beta(h)}\rp^{-1}
\sum_{\alpha\in\raff_+}\alpha(h)\RR^-(s)_{\beta-\alpha}(\Omega_\alpha+\mathrm{Q}_\alpha) \\
&= \hbar\sum_{k\geqslant 0} \frac{\adT(h)^k}{(s\beta(h))^{k+1}}
\sum_{\alpha\in\raff_+}\alpha(h)\RR^-(s)_{\beta-\alpha}(\Omega_\alpha+\mathrm{Q}_\alpha). 
\end{aligned}
\end{equation}
This equation determines $\RR^-(s)_\beta$ recursively in terms of $\RR^-(s)_{\gamma}$ with $\mathsf{ht}(\gamma)<\mathsf{ht}(\beta)$, and thus establishes the uniqueness of $\{\RR^-(s)_\beta\}_{\beta\in \Qaff_+}$ satisfying \eqref{eq:IER-} and the initial condition $\RR^-(s)_0=1\otimes 1$.  

To prove existence, we proceed as in \cite[\S4.2]{GTLW}. 
Choose $h\in\hkm$ such that $\beta(h)\neq 0$ for every nonzero $\beta\in\Qaff_+$. 
Corresponding to this choice of $h$, we may define a family of elements $\{\RR^-(s)_\beta\}_{\beta\in \Qaff_+}$ recursively on the height of $\beta$ using \eqref{R^-:recur} and the condition $\RR^-(s)_0=1$. This produces elements $\RR^-(s)_\beta\in (\Ykm_{-\beta}\otimes \Ykm_\beta)[\![s^{-1}]\!]$ which by definition satisfy
\begin{equation}\label{eq:IER-:fixed}
\RR^-(s,z)\Delta_s^z(\Top(h))= \ddelta{s}^{\!z}(\Top(h)) \RR^-(s,z)
\end{equation}
for our fixed choice of $h$, where $\RR^-(s,z)=\sum_{\beta\in \Qaff_+}\RR^-(s)_\beta z^{\mathsf{ht}(\beta)}$.  To see that $\RR^-(s,z)$ in fact satisfies \eqref{eq:IER-} in general, let $h'\in \hkm$ and note that, since $\Delta_s^z$
and $\ddelta{s}^{\!z}$ are algebra homomorphisms, the two series 
\begin{equation*}
\RR^-(s,z)\Delta_s^z(\Top(h')) \quad \text{ and }\quad  \ddelta{s}^{\!z}(\Top(h'))\RR^-(s,z)
\end{equation*}
both solve \eqref{eq:IER-:fixed}, are of the form $\sum_{\beta\in \Qaff_+}\EuScript{R}(s)_\beta z^{\mathsf{ht}(\beta)}$ with $\EuScript{R}(s)_\beta\in (\Ykm_{-\beta}\otimes \Ykm_{\beta})[s;s^{-1}]\!]$ for each $\beta\in \Qaff_+$,   and have $\Ykm_0\otimes \Ykm_0$ components equal to $(\tau_s\otimes \Id)\square(\Top(h'))$. They therefore coincide by the uniqueness argument given above (see \eqref{R^-:recur}).
\let\qed\relax
\end{pf}

\begin{pf}[Proof of Parts \eqref{R-(s,z):2}--\eqref{R-(s,z):4}] By Part \eqref{map-T:3} of Theorem \ref{T:map-T}, $\Omega_\alpha+\mathrm{Q}_\alpha$ lies in $\Ysupkm{-}_{-\alpha}\otimes \Ysupkm{+}_{\alpha}$ for each $\alpha\in \raff_+$. Therefore, the recursion \eqref{R^-:recur} implies that 
\[
\RR^-(s)_\beta\in s^{-\nu(\beta)} 
\Pseries{\lp\Ysupkm{-}_{-\beta}\otimes\Ysupkm{+}_{\beta}\rp}{s^{-1}}
\]
for all $\beta\in \Qaff_+$, with $\RR^-(s)_{\alpha}=(\Omega_\alpha+\mathrm{Q}_\alpha)s^{-1} +O(s^{-2})$ for all $\alpha\in \raff_+$. This proves Parts \eqref{R-(s,z):2} and \eqref{R-(s,z):4} of the theorem. Similarly, given $a,b\in \C$,  Parts \eqref{map-T:2} and \eqref{map-T:3} of Theorem \ref{T:map-T} imply that the elements $(\tau_a\otimes \tau_b)(\RR^-(s)_\beta)$ and $\RR^-(s+a-b)_\beta$ both satisfy \eqref{R^-:recur-rat} with $s$ replaced by $s+a-b$, and thus coincide by uniqueness (see also \cite[\S4.3]{GTLW}). Note that the argument given in {\em loc. cit.} doesn't work as stated since $\Top(\hkm)\oplus\hkm$ is not stable under $\tau_b$. Still $\mathrm{ad}$ of this subspace is, which is all that is needed here. \let\qed\relax 
\end{pf}

\begin{pf}[Proof of Part \eqref{R-(s,z):1}] Since the desired identity is satisfied for $x\in \hkm\oplus \Top(\hkm)$, it  suffices to prove that 
\begin{equation*}
\RR^-(s,z)\Delta_s^z(x_{i,0}^\pm)= \ddelta{s}^{\!z}(x_{i,0}^\pm) \RR^-(s,z)
\end{equation*}
 for each $i\in \hbfI$. This is proven using a rank $1$ reduction argument, exactly
as in \cite[\S 4.7]{GTLW}, with the $\Lsl_2$ case having been verified
in \cite[\S 4.8]{GTLW}.
\end{pf} \let\qed\relax
\end{pf}
\begin{rem}
We emphasize that the use of Theorem \ref{T:map-T} is essential in the above proof, and in particular in the existence and uniqueness argument for $\{\RR^-(s)_\beta\}_{\beta\in \Qaff_+}$. Indeed, unlike in the setting of \cite[Thm~4.1]{GTLW}, it is not sufficient to replace $\{\Top(h)\}_{h\in \hkm}$ by the family $\{t_{i,1}\}_{i,\in \hbfI}$ in \eqref{eq:IER-}, as it is not true that for any nonzero $\beta\in \Qaff_+$ there is $h\in \haff$ such that $\beta(h)\neq 0$. The existence of such an $h$ is necessary to pass from \eqref{R^-:recur-rat} to \eqref{R^-:recur}.

 In addition, unlike in the finite case considered in \cite{GTLW}, the infinite sum $\sum_\beta \RR^-(s)_\beta$ does
not converge to an element of $\Ykm^{\otimes 2}[\![s^{-1}]\!]$. This is the reason behind the introduction of the 
auxiliary parameter $z$ in the statement of the theorem.  
\end{rem}
\subsection{The operators $\RR^-_{V_1,V_2}(s)$}\label{ssec:R^--rep}

Let $V_1,V_2\in\OY$, with $\pi_\ell:\Ykm\to\End(V_\ell)$ being
the action homomorphisms. By the category $\mathcal{O}$ condition,
the matrix entries of $\pi_1\otimes\pi_2 (\Omega_z)$ and $\pi_1\otimes\pi_2 (\mathrm{Q}_z)$
are polynomials in $z$ and therefore can be evaluated at
$z=1$. Let
\[
\Omega_{V_1,V_2} \ceq   \left.\pi_1\otimes\pi_2 (\Omega_z)\right|_{z=1}
\quad \text{ and }\quad 
\mathrm{Q}_{V_1,V_2} \ceq   \left.\pi_1\otimes\pi_2 (\mathrm{Q}_z)\right|_{z=1}.
\]

We have the following representation--theoretic version
of Theorem \ref{thm:R-}.
\begin{thm}\label{thm:R^--reps}
Let $V_1,V_2\in\OY$. Then, there exist a unique $\End(V_1\otimes V_2)$--valued,
rational function of $s$, denoted by $\RR^-_{V_1,V_2}(s)$,
satisfying the following properties.
\begin{enumerate}[font=\upshape]\itemsep0.25cm
\item\label{R_V^-:1} (Normalization) $\RR^-_{V_1,V_2}(\infty)=\Id_{V_1\otimes V_2}$.

\item\label{R_V^-:2} (Triangularity) $\ds\RR^-_{V_1,V_2}(s) = \sum_{\beta\in\Qaff_+}
\RR^-_{V_1,V_2}(s)_\beta$, where, for any $\mu_1,\mu_2\in\hkm^*$
\[
\RR^-_{V_1,V_2}(s)_\beta : V_1[\mu_1]\otimes V_2[\mu_2]
\to V_1[\mu_1-\beta]\otimes V_2[\mu_2+\beta].
\]

\item\label{R_V^-:3} (Intertwiner) For every $h\in\hkm$, the following
intertwining relation holds:
\[
\left[\square(\Top(h)) + s h\otimes 1, \RR^-_{V_1,V_2}(s)\right]
= \hbar \RR_{V_1,V_2}^-(s) [h\otimes 1, \Omega_{V_1,V_2}+\mathrm{Q}_{V_1,V_2}].
\]
\end{enumerate}

Moreover, this operator has the following properties:

\begin{enumerate}[font=\upshape, start=4]\itemsep0.25cm
\item\label{R_V^-:4} If $\RR^-(s,z)$ is as in Theorem \ref{thm:R-}, then 
\[
\RR^-_{V_1,V_2}(s) = \left.\pi_1\otimes \pi_2 (\RR^-(s,z))\right|_{z=1}.
\]

\item\label{R_V^-:5} $\ds \RR^-_{V_1,V_2}(s) : V_1\kmtensor{s} V_2 \to V_1\dtensor{s} V_2$
is a $\Yhg$--intertwiner, which is natural in $V_1,V_2$ and compatible
with the shift automorphism:
\[
\RR^-_{V_1(a),V_2(b)}(s) = \RR^-_{V_1,V_2}(s+a-b).
\]

\item\label{R_V^-:6} (Cocycle equation) The following diagram commutes, for every $V_1,V_2,V_3\in\OY$: 
		\[
\begin{tikzcd}[column sep=huge, row sep=large]
	(V_1\kmtensor{s_1} V_2)\kmtensor{s_2} V_3
	\arrow[r, "="]
	\arrow[d, "\RR^-_{V_1,V_2}(s_1)\otimes \Id_{V_3}"']
	&
	V_1\kmtensor{s_1+s_2}(V_2\kmtensor{s_2}\otimes V_3)
	\arrow[d, "\Id_{V_1}\otimes \RR^-_{V_2,V_3}(s_2)"]
	\\
	(V_1\dtensor{s_1} V_2)\kmtensor{s_2} V_3
	\arrow[d, "\RR^-_{V_1\dtensor{s_1}V_2,V_3}(s_2)"']
	&
	V_1\kmtensor{s_1+s_2}(V_2\dtensor{s_2}\otimes V_3)
	\arrow[d, "\RR^-_{V_1,V_2\dtensor{s_2}V_3}(s_1+s_2)"]
	\\
	(V_1\dtensor{s_1} V_2)\dtensor{s_2} V_3
	\arrow[r, "="]
	&
	V_1\dtensor{s_1+s_2}(V_2\dtensor{s_2}\otimes V_3)
\end{tikzcd}
\]

\end{enumerate}

\end{thm}

\vspace{0.25cm}
\begin{rem}
We wish to highlight that the cocycle equation only holds for the
rational intertwiners, not for the formal $\RR^-(s,z)$; see 
\cite[Remark 4.1]{GTLW}.
\end{rem}

\begin{pf}
For notational simplicity, we drop the subscript $V_1,V_2$.
Let us fix a $\mu\in\hkm^*$ and consider the intertwining
equation valued in the finite--dimensional
vector space $\mathcal{E}\ceq  \End((V_1\otimes V_2)[\mu])$.
For $h\in\hkm$,
let $\mathsf{r}(h) \ceq   [h\otimes 1, \Omega+\mathrm{Q}]$,
viewed as a strictly lower triangular operator on
$\mathcal{E}$.
Now, we have:
\[
A(h,s)\ceq  \ad(\square(\Top(h))+sh\otimes 1) - \rho(\hbar\mathsf{r}(h)) 
\in \End(\mathcal{E})[s],
\]
where $\rho(X)$ denotes the operator of right multiplication by $X$.
This operator is triangular in the weight grading on $\mathcal{E}$,
with diagonal blocks being linear in $s$, for generic $h$.
The same argument given in the proof of Theorem \ref{thm:R-},
carried out in $\End(\mathcal{E})$, shows that there is
a unique unipotent solution to \eqref{eq:IER-}. Moreover, each
$\RR^-(s)_\beta$ is a rational function of $s$, vanishing
at $s=\infty$ for $\beta>0$. This proves the first part of the
theorem.

The uniqueness of the solution, together with Theorem \ref{thm:R-}, imply that
the Taylor series expansion of this rational function
at $s=\infty$
is the same as the evaluation of $\RR^-(s,z)$ on $V_1\otimes V_2$,
specialized at $z=1$. This proves \eqref{R_V^-:4}. Part \eqref{R_V^-:5} is now a consequence
of Theorem \ref{thm:R-}.

The cocycle equation also rests on uniqueness, and
the coassociativity of $\Delta_s$ and $\ddelta{s}$ on representations
(see \cite[\S 4.7]{GTLW}). Namely, consider
the two sides of the cocycle equation:
\[
\begin{aligned}
L(s_1,s_2) &\ceq   \RR^-_{V_1\dtensor{s_1}V_2,V_3}(s_2)\circ 
\left(\RR^-_{V_1,V_2}(s_1)\otimes \Id_{V_3}\right), \\
R(s_1,s_2) &\ceq   \RR^-_{V_1,V_2\dtensor{s_2}V_3}(s_1+s_2)\circ 
\left(\Id_{V_1}\otimes \RR^-_{V_2,V_3}(s_2)\right).
\end{aligned}
\]
It is straightforward to see that both of these operators
intertwine the two actions of $\Top(h)$. That is, 
using the formulae for $\Delta_s$ and Theorem
\ref{T:map-T} \eqref{map-T:3}, these operators solve
the following equation:
\begin{gather*}
\ad\left(\sum_{a=1}^3 \Top(h)^{(a)} + (s_1+s_2)h^{(1)}
+s_2 h^{(2)}\right)\cdot X(s_1,s_2) = \\
\hbar X(s_1,s_2)(\mathsf{r}_{12}(h)+\mathsf{r}_{23}(h)
+\mathsf{r}_{13}(h)).
\end{gather*}
As in \cite[\S4.6]{GTLW}, there is at most one solution
to this equation which is triangular in the sense
that $X = \sum_{\beta,\gamma\in\Qaff_+} X_{\beta,\gamma}$, where, for each $\mu_1,\mu_2,\mu_3\in\h^*$,
\[
X_{\beta,\gamma} : V_1[\mu_1]\otimes V_2[\mu_2]\otimes V_3[\mu_3]
\to V_1[\mu_1-\beta]\otimes V_2[\mu_2+\beta-\gamma]\otimes V_3[\mu_3+\gamma]
\]
and $X_{0,0}=\Id$. That $L(s_1,s_2)$ and $R(s_1,s_2)$ take this form, and thus coincide by uniqueness, is readily deduced from Part \eqref{R_V^-:2} of the theorem together with the fact that the Drinfeld coproduct $\ddelta{s}$ satisfies 
\begin{equation*}
\ddelta{s}(\Ysupkm{+})\subset (\Ykm\otimes\Ysupkm{+})[s;s^{-1}]\negthinspace] \quad \text{ and }\quad \ddelta{s}(\Ysupkm{-})\subset (\Ysupkm{-}\otimes\Ykm)[s;s^{-1}]\negthinspace].
\end{equation*}
This completes the proof of Part \eqref{R_V^-:6}, and thus the proof of the theorem. \qedhere
\end{pf}

\section{Construction of $\RR^0(s)$}\label{sec:ab-R}

In this section, we introduce an additive, regular difference equation,
whose coefficients come
from $\Prim(\Ykm)^{\otimes 2}$, where $\Prim(\Ykm)$ is the linear span of
$\{t_{i,r}\}_{i\in\hbfI, r\in\N}$. 
We show that the
exponential of the two fundamental solutions of this difference equation 
give rise to two meromorphic
braidings on $(\OY,\dtensor{s})$, related by a unitarity condition. 
They have the same asymptotic expansion,
which is the unique formal solution of our difference equation, \ie
a formal abelian $R$--matrix.

\subsection{Main construction and result}\label{ssec:R0-main}

We begin by stating the main result of this section.

\begin{thm}\label{thm:R0-main}
Given $V_1,V_2\in\OY$, there exist two meromorphic
 $\End(V_1\otimes V_2)$--valued
functions, $\RR^{0,\eta}_{V_1,V_2}(s)$, $\eta\in\{\uparrow,\downarrow\}$,
which are natural in $V_1$ and $V_2$
and have the following properties:
\begin{enumerate}[font=\upshape]\itemsep0.25cm
\item\label{R0:fn} Let $\mathbb{H}^{\uparrow} = \{s\in\C : \Re(s/\hbar)\gg 0\} 
= -\mathbb{H}^{\downarrow}$.
Then, for every $V_1,V_2\in\OY$, $\RR^{0,\eta}_{V_1,V_2}(s)$ is holomorphic
and invertible in $\mathbb{H}^{\eta}$ and approaches $\Id_{V_1\otimes V_2}$
as $|s|\to\infty$.

\item\label{R0:int} For every $V_1,V_2$, the following is a $\Ykm$--intertwiner:
\[
(1\ 2) \circ \RR^{0,\eta}_{V_1,V_2}(s) : V_1\dtensor{s} V_2 \to
(V_2\dtensor{-s}V_1)(s)\ .
\]

\item\label{R0:shift} The following holds, for every $V_1,V_2\in\OY$
and $a,b\in\C$:
\[
\RR^{0,\eta}_{V_1(a),V_2(b)}(s) = \RR^{0,\eta}_{V_1,V_2}(s+a-b)\ .
\]

\item\label{R0:cabling} For every $V_1,V_2,V_3\in\OY$, we have:
\[
\begin{aligned}
\RR^{0,\eta}_{V_1\dtensor{s_1}V_2,V_3}(s_2) &= \RR^{0,\eta}_{V_1,V_3}(s_1+s_2)
\RR^{0,\eta}_{V_2,V_3}(s_2)\ , \\
\RR^{0,\eta}_{V_1,V_2\dtensor{s_2}V_3}(s_1+s_2) &= \RR^{0,\eta}_{V_1,V_3}(s_1+s_2)
\RR^{0,\eta}_{V_1,V_2}(s_1)\ .
\end{aligned}
\]

\item\label{R0:unitary} For every $V_1,V_2\in\OY$, we have:
\[
\RR^{0,\uparrow}_{V_1,V_2}(s)^{-1} = (1\ 2)\circ
\RR^{0,\downarrow}_{V_2,V_1}(-s)\circ (1\ 2)\ .
\]

\item\label{R0:asym} There exists a unique formal series
$\RR^0(s)\in\Pseries{\Ysupkm{0}^{\otimes 2}}{s^{-1}}$ such that,
for every $V_1,V_2\in\OY$, with $\pi_\ell:\Ykm\to\End(V_\ell)$
the corresponding action homomorphism, we have:
\[
\RR^{0,\eta}_{V_1,V_2}(s)\sim \pi_1\otimes\pi_2 (\RR^0(s)),\ 
\text{as } \pm\Re(s/\hbar) \to\infty .
\]
This asymptotic expansion remains valid in a larger sector $\Sigma^\eta_\delta$,
for any $\delta>0$, where, if $\theta = \arg(\hbar)$ then:
\[
\Sigma^\uparrow_\delta \ceq   \{re^{\iota\phi} : r\in\R_{>0}, 
\phi\in (\theta-\pi+\delta,\theta+\pi-\delta)\} = -\Sigma^\downarrow_\delta.
\]
(see Figure \ref{afig:sector} in Appendix \ref{app:LapDE} where $\chi=\hbar/2$).

\item\label{R0:leading} The first order term of $\RR^0(s)$ is given by:
\[
\RR^0(s) = \exp\lp s^{-1}\lp \frac{1}{\hbar}\ctsing + 
\frac{\hbar}{4\qdzero} \ctzero \rp + O(s^{-2})\rp\,.
\]
Here, the constant $\sczero$ and the tensor $\ctzero$
are given in Lemma \ref{lem:Rg-terms}, $\qdzero$ in Lemma \ref{lem:order}
below, and
\[
\ctsing = \frac{\sczero}{\qdzero} \lp \frac{\caff_2\otimes\caff_0}{2}
-\caff_1\otimes\caff_1 + \frac{\caff_0\otimes \caff_2}{2}\rp.
\]
\end{enumerate}
\end{thm}


The proof of this theorem is carried out in Sections \ref{ssec:R0-thm-pf}--\ref{ssec:borel-growth}. 

\subsection{Rationality}\label{ssec:rat-R0}
In this section, we prove the 
rationality property of (a normalization of) the abelian $R$--matrix. The proof relies on
Theorem~\ref{thm:R0-main} \eqref{R0:int} and the fact that
\begin{equation}\label{eq:R0-zero}
	[\RR^{0,\eta}_{V_1,V_2}(s), \pi_1(y_1)\ten\pi_2(y_2)]=0
\end{equation}
for every $V_1,V_2\in\OY$, with $\pi_\ell:\Ykm\to\End(V_\ell)$
the corresponding action homomorphism, and $y_1,y_2\in\Ysupkm{0}$.
This identity, though not stated in the theorem, is evident from the explicit construction of $\RR^{0}(s)$, see Section~\ref{ssec:R0-thm-pf}.

Fix $\eta\in\{\uparrow,\downarrow\}$ and assume that $V_1,V_2\in\OY$ are two highest--weight representations.
Let $\lambda_\ell\in\hkm^*$ be their highest weights, and fix
$\mathsf{v}_\ell\in V_\ell[\lambda_\ell]$ highest--weight vectors
($\ell=1,2$). 
By \eqref{eq:R0-zero}, $\RR^{0,\eta}_{V_1,V_2}(s)$  is a weight zero operator, hence it preserves the one--dimensional vector space $V_1[\lambda_1]\otimes
V_2[\lambda_2]$.
Let $\mathfrak{f}^\eta(s)$ be the corresponding eigenvalue:
\[
\RR^{0,\eta}(s)\cdot \mathsf{v}_1\otimes\mathsf{v}_2
= \mathfrak{f}^\eta(s) \mathsf{v}_1\otimes\mathsf{v}_2\ .
\]
\begin{thm}\label{thm:rat-R0}
	With the notational set up as above, the normalized operator
	\[
	\mathsf{R}^0_{V_1,V_2}(s) \ceq   \mathfrak{f}^\eta(s)^{-1}\RR^{0,\eta}_{V_1,V_2}(s)
	\]
	is independent of $\eta$, and is rational in $s$.
\end{thm}

\begin{pf}
	Let $\mu_\ell\in P(V_\ell)$ be two weights ($\ell=1,2$). We argue by
	induction on $\hit(\lambda_1-\mu_1) + \hit(\lambda_2-\mu_2)$.
	For notational convenience we will drop $V_1,V_2$ from the subscripts
	of our operators.
	The base case being clear, we focus on the induction step.
	Note that, by the highest--weight property, we have
	\[
	V_1[\mu_1]\otimes V_2[\mu_2] = \sum_{\begin{subarray}{c}
			k_1,k_2\in\hbfI\\ r_1,r_2\in\N \end{subarray}}
	x^-_{k_1,r_1}\lp V_1[\mu_1+\alpha_{k_1}]\rp 
	\otimes
	x^-_{k_2,r_2}\lp V_2[\mu_2+\alpha_{k_2}]\rp. 
	\]
	Moreover, by the commutation relation $[t_{k,1},x^-_{k,r}] = -2d_kx^-_{k,r+1}$,
	for any $V\in\OY$ and $\mu\in P(V)$, we have
	\[
	\sum_{r\in\N} x^-_{k,r}\lp V[\mu]\rp = 
	\Ysupkm{0}\cdot x^-_{k,0}\lp V[\mu]\rp.
	\]
	By \eqref{eq:R0-zero}, $\RR^{0,\eta}(s)$ commutes with the action of
	$\Ysupkm{0}\otimes \Ysupkm{0}$ on $V_1\otimes V_2$. Thus, it suffices
	to show that the
	action of $\mathsf{R}^0(s)$ on the image of  
	$x^-_{k,0}\otimes 1$ (and $1\otimes x^-_{k,0}$) is rational, 
	and independent of $\eta$. We focus of the former, as the latter will follow
	either with a similar proof, or using unitarity.
	We will use the following commutation relation between $\RR^{0,\eta}(s)$
	and $x^-_{k,0}\otimes 1$, which will be established in the proof of Part \eqref{R0:int} of Theorem \ref{thm:R0-main}, outlined in Section \ref{ssec:R0-thm-pf}:
	\begin{equation}\label{eq:pf-ratR0-1}
		\Ad(\RR^{0,\eta}(s)^{-1})\cdot (x^-_{k,0}\otimes 1)
		= x^-_{k,0}\otimes 1 + \mathfrak{Y}_k(s),
	\end{equation}
	where
	\[
	\mathfrak{Y}_k(s) = 
	\hbar\sum_{N\geqslant 0} s^{-N-1} \lp
	\sum_{n=0}^N (-1)^{n} \cbin{N}{n} x^-_{k,n}\otimes \xi_{k,N-n}\rp .
	\]

	The right--hand side of \eqref{eq:pf-ratR0-1} has the following
	contour integral representation, which immediately shows that
	it is rational in $s$ (see \cite[\S4.5]{sachin-valerio-III}):
	\[
	x^-_{k,0}\otimes 1 + \mathfrak{Y}_k(s) = 
	\frac{1}{\hbar} \oint_{C_1} x^-_k(u)\otimes \xi_k(u+s)\, du\ .
	\]
	Here, $C_1$ is a contour enclosing all the poles of $x^-_k(u)$
	acting on $V_1[\mu_1+\alpha_k]$, and $s$ is so large that $\xi_k(u+s)$
	acting on $V_2[\mu_2]$ is holomorphic on and within $C_1$.
	
	Now, let $w_1\in V_1[\mu_1+\alpha_k]$ and $w_2\in V_2[\mu_2]$. Then, since $\Ad(\RR^{0,\eta}(s)) = \Ad(\mathsf{R}^0(s))$, we have 
	\[
	\begin{aligned}
		\mathsf{R}^0(s)^{-1}\circ &(x^-_{k,0}\otimes 1)\cdot (w_1\otimes w_2)\\
		&=
		\lp\Ad(\RR^{0,\eta}(s)^{-1})\cdot (x^-_{k,0}\otimes 1)\rp
		\lp \mathsf{R}^0(s)^{-1}\cdot (w_1\otimes w_2)\rp\\
		&=
		\lp x^-_{k,0}\otimes 1+ \mathfrak{Y}_k(s)\rp \cdot
		\lp \mathsf{R}^0(s)^{-1}\cdot (w_1\otimes w_2)\rp.
	\end{aligned}
	\]
	The last line gives a vector in $V_1[\mu_1]\otimes V_2[\mu_2]$
	depending rationally on $s$, by the induction hypothesis combined
	with the rationality of $x^-_{k,0}\otimes 1 + \mathfrak{Y}_k(s)$.
	The theorem is proved.
\end{pf}

\subsection{Proof of Theorem \ref{thm:R0-main}}\label{ssec:R0-thm-pf}
Our proof is based on an explicit construction of $\RR^0(s)$, both
as a formal series, and as an $\End(V_1\otimes V_2)$-valued meromorphic function of $s$. This $\RR^0(s)$ has the form
\[
\RR^0(s) = \exp(\Rl(s)),\ \text{ where } \Rl(s)\in 
s^{-1}\Pseries{\Prim(\Ykm)^{\otimes 2}}{s^{-1}}.
\]
The $\Rl(s)$, in turn, is obtained in a few steps. We outline
the steps here, which are carried out in Sections \ref{ssec:T-cartan}--\ref{ssec:Rl-props}.
To begin, consider the  difference equation
\begin{equation}\label{eq:L0-diff}
(\sop-\sop^{-1})\det(\bfB(\sop)) \cdot \Rl(s) = \Rg(s),
\end{equation}
where $\sop$, $\bfB(\sop)$ and $\Rg(s)$ are defined as follows:
\begin{itemize}\itemsep0.25cm
\item $\sop$ is the shift operator: $\sop\cdot f(s) = f(s-\hbar/2)$.

\item $\bfB(\sop)$ is the symmetrized, affine $\sop$--Cartan matrix (see
Section \ref{ssec:T-cartan} below). 

\item $\Rg(s) = \sum_{ij} \bfB(\sop)^{\ast}_{ji}\cdot \Rtau_{ij}(s)$, where $\bfB(\sop)^\ast$ is the
adjoint matrix to $\bfB(\sop)$ and 
\[
\Rtau_{ij}(s) \ceq   
\hbar^2
\sum_{m\geqslant 1} m!s^{-m-1} \sum_{\begin{subarray}{c} a,b\geqslant 0 \\ a+b=m-1
\end{subarray}} (-1)^a 
\frac{t_{i,a}}{a!}\otimes \frac{t_{j,b}}{b!}
\]
\end{itemize}

In Section \ref{ssec:tau}, we prove important properties of the operators
$\Rtau_{ij}(s)$. These are used to establish nearly identical properties for the operator
$\Rg(s)$  in Section \ref{ssec:Rg}. We compute the first
few terms of $\Rg(s)$ in Section \ref{ssec:Rg-reg} and show
that the coefficients of $s^{-2}$ and $s^{-3}$ are central.
This fact is used in Section \ref{ssec:R0-diff} to regularize the difference equation
\eqref{eq:L0-diff} to
\begin{equation}\label{eq:Rl-diff}
(\sop-\sop^{-1})\det(\bfB(\sop)) \cdot \Rl(s) = \Rg_{\reg}(s),
\end{equation}
where $\Rg_{\reg}(s)$ is obtained by removing the aforementioned central terms from $\Rg(s)$, and is  defined explicitly in \eqref{eq:Rg-reg}.
We show in  Corollary \ref{cor:Rl-formal} that this equation
has a unique formal solution. The properties of $\Rl(s)$,
analogous to the assertions in this theorem, are obtained
in Proposition \ref{pr:Rl}. 

Now, given $V_1,V_2\in\OY$, the evaluation of the difference equation \eqref{eq:Rl-diff} on $V_1\otimes V_2$
has coefficients from the following subalgebra of
$\End(V_1\otimes V_2)$:
\[
\bigoplus_{\mu_1\in P(V_1),\mu_2\in P(V_2)} \End_{\Ysupkm{0}}(V_1[\mu_1])
\otimes \End_{\Ysupkm{0}}(V_2[\mu_2]).
\]
Thus, fixing $\mu_1,\mu_2$, we can view \eqref{eq:Rl-diff} as an
equation for a
finite size matrix--valued function of $s$. In Appendix \ref{app:LapDE}
(see Theorem \ref{appA:thm}),
we establish the existence and uniqueness of the solutions of such equations with the prescribed asymptotics, for
any right--hand side of $O(s^{-4})$ (this $4$ is $1$ plus the
order of vanishing of the polynomial on the left--hand side at $\sop=1$;
see Lemma \ref{lem:order} below). 
We verify that our difference equation satisfies the hypotheses
of Theorem \ref{appA:thm} in Section \ref{ssec:borel-growth} (see
Lemma \ref{lem:borel-growth}).

Given these preparatory results, let us prove the theorem. Note that
\eqref{R0:fn} and \eqref{R0:asym} follow from definitions. 
Parts \eqref{R0:shift} and \eqref{R0:unitary} of the theorem are a direct consequence of 
Parts \eqref{Rl:shift} and \eqref{Rl:unitary} of Proposition \ref{pr:Rl}, respectively. Part \eqref{R0:cabling}
follows from the fact that the coefficients of $\Rl(s)$ are tensors of primitive
 elements with respect to the Drinfeld coproduct.

Let us prove \eqref{R0:int}. It is clear that $\RR^{0,\eta}_{V_1,V_2}(s)$
commutes with operators from $\Ysupkm{0}$. Therefore, it is enough to
show that $(1\ 2)\circ\RR^{0,\eta}_{V_1,V_2}(s)$ commutes with
$\ddelta{s}(x_{k,0}^{\pm})$, for all $k\in\hbfI$. We focus on the
$+$ case, the other one being similar. Recall that
\[
\ddelta{s}(x_{k,0}^+) = \square(x_{k,0}^+)
+\hbar\sum_{N\geqslant 0} s^{-N-1} \lp
\sum_{a+b=n} (-1)^{a+1} \cbin{N}{a} \xi_{k,a}\otimes x^+_{k,b}
\rp.
\]
Let $\mathfrak{X}_k(s)$ denote the second term on the right--hand side.
The desired commutation relation decouples to the following two:
\[
\begin{aligned}
\Ad(\RR^{0,\eta}_{V_1,V_2}(s))\cdot x^+_{k,0}\otimes 1 &=
x^+_{k,0}\otimes 1 + \mathfrak{X}_k^{\op}(-s), \\
\Ad(\RR^{0,\eta}_{V_1,V_2}(s)^{-1})\cdot 1\otimes x^+_{k,0} &=
1\otimes x^+_{k,0} + \mathfrak{X}_k(s).
\end{aligned}
\]
Note that the second follows from the first, given Part \eqref{R0:unitary} of the theorem (unitarity). For the first,
we write
\[
x^+_{k,0}\otimes 1 + \mathfrak{X}^{\op}_k(-s) = \sum_{a=0}^{\infty}
x^+_{k,a}\otimes \partial_s^{(a)}(\xi_k(s)).
\]
The desired relation is then a consequence of the following claim. 

\noindent {\bf Claim:} For every $k\in\hbfI, n\in\N$ and $y\in\Ysupkm{0}$,
we have
\[
\ad(\Rl_{V_1,V_2}^{\eta}(s))\cdot (x^+_{k,n}\otimes y) = 
\sum_{a=0}^{\infty} x^+_{k,n+a}\otimes \partial_s^{(a)}(t_k(s))y.
\]

We show that this identity is true for the formal $\Rl(s)$ in
Proposition \ref{pr:Rl} \eqref{Rl:comm}. To deduce it for the
meromorphic functions $\Rl_{V_1,V_2}^\eta(s)$, we apply the
difference operator $\Rd(\sop)=(\sop-\sop^{-1})\det(\bfB(\sop))$ to both
sides to conclude that they are solutions to the same difference
equation, with the same asymptotic expansion as $s\to\infty$.
Hence, by the
uniqueness result of Theorem \ref{appA:thm}, they are equal.

\noindent{\em A word on set up and proofs.}
Our construction of $\RR^0(s)$
is weight preserving. So, let $V_1,V_2\in\OY$ and $\mu_\ell\in P(V_\ell)$
($\ell=1,2$) be two fixed weights. All the operators considered
in Sections \ref{ssec:tau}--\ref{ssec:Rg-reg} can be viewed as
meromorphic $\End(V_1[\mu_1]\otimes V_2[\mu_2])$--valued functions
of a complex variable $s$, which are regular near $s=\infty$.
We will be interested in both their functional nature, and naturality
with regards to $V_1,V_2$. Therefore, their Taylor series expansions
will be shown to be the evaluation on $V_1\otimes V_2$ of
an element of $\Pseries{\Ysupkm{0}^{\otimes 2}}{s^{-1}}$.

An identity among such operators can be shown either functionally, or formally.
The functional proofs, by which we mean proofs involving the usual
``contour deformation" trick, can be found in \cite{sachin-valerio-III}
(we will make precise citation when needed). Here, we also give the formal
proofs --- both for completeness and for their aesthetic beauty.

\subsection{$\sop$--Cartan matrix}\label{ssec:T-cartan}
Let $\sop$ be an indeterminate and let $\bfB(\sop) \ceq   ([d_ia_{ij}]_\sop)
\in \mathrm{M}_{\hbfI\times\hbfI}(\Z[\sop,\sop^{-1}])$. Here,
we use the standard notation of Gaussian numbers:
\[
[n]_q \ceq   \frac{q^n-q^{-n}}{q-q^{-1}}.
\]
Let $\bfB(\sop)^\ast$ be the adjoint matrix of $\bfB(\sop)$, so that
\[
\bfB(\sop)^\ast \bfB(\sop) = \bfB(\sop)\bfB(\sop)^\ast = \det(\bfB(\sop))\Id.
\]
The following lemma is crucial, and is proved by direct
inspection (see Appendix \ref{app:QCM} for the table of determinants
of symmetrized, affine $\sop$--Cartan matrices).

\begin{lem}\label{lem:order}
The order of vanishing of $\det(\bfB(\sop))$ at $\sop=1$ is $2$. That is,
\[
\qdzero\ceq   \left.
\frac{\det(\bfB(\sop))}{(\sop-\sop^{-1})^2}
\right|_{\sop=1} \text{ exists and } \neq 0\ .
\]
Moreover, all the zeroes of $\det(\bfB(\sop))$ have modulus $1$.
\end{lem}

Below we will view $\sop$ as an operator on either rational, matrix--valued
functions of $s$, regular near $s=\infty$; or formal series in $s^{-1}$,
via:
\begin{equation}\label{eq:T-action}
\sop\cdot f(s) = f(s-\hbar/2)\ .
\end{equation}

\subsection{The $\Rtau$ operators}\label{ssec:tau}

Recall that for
$i,j\in\hbfI$, we defined
\begin{equation}\label{eq:tau-formal}
\Rtau_{ij}(s) = \hbar^2
\sum_{m\geqslant 1} m!s^{-m-1} \sum_{\begin{subarray}{c} a,b\geqslant 0 \\ a+b=m-1
\end{subarray}} (-1)^a 
\frac{t_{i,a}}{a!}\otimes \frac{t_{j,b}}{b!}\ .
\end{equation}
In Lemma 6.5 of \cite{GTLW},  it was shown that, when evaluated on 
$V_1[\mu_1]\otimes V_2[\mu_2]$, $\Rtau_{ij}(s)$ is the
Taylor series expansion near $s=\infty$ of the 
contour integral
\begin{equation}\label{eq:tau}
\Rtau_{ij}(s) \ceq  \oint_{C_1} t_i'(u)\otimes t_j(u+s)\, du,
\end{equation}
where $C_1$ is a contour enclosing zeroes and poles of $\xi_i(u)$
acting on $V_1[\mu_1]$, $t_i'(u)$ is the derivative of $t_i(u)$, and $s$ is large enough so that
$t_j(u+s)$ acting on $V_2[\mu_2]$
is holomorphic on and within $C_1$. As usual, we suppress
$2\pi\iota$ factor in our notations $\oint = \frac{1}{2\pi\iota} \int$.
Here, $t_j(w)$ is viewed as a single--valued function
on a cut plane as in Section \ref{ssec:log}.
This expression is used to show that $\exp(\Rtau_{ij}(s))$
becomes a rational function of $s$ on $V_1[\mu_1]\otimes V_2[\mu_2]$.
(see the proof of Proposition \ref{pr:tau} \eqref{tau:fun} below).
We can also rewrite \eqref{eq:tau-formal} 
above using $B_i(z)$ (see Section \ref{ssec:tir} above),
and the fact that $m!s^{-m-1} = (-\partial_s)^m\cdot s^{-1}$, as
\begin{equation}\label{eq:tau-borel}
\Rtau_{ij}(s) = \left. zB_i(-z)\otimes B_j(z)\right|_{z=-\partial_s} \cdot s^{-1}
 = \left. B_i(-z)\otimes B_j(z)\right|_{z=-\partial_s} \cdot s^{-2}.
\end{equation}

\begin{prop}\label{pr:tau}
The elements $\{\Rtau_{ij}(s)\}_{i,j\in\hbfI}$ have the following properties:
\begin{enumerate}[font=\upshape]\itemsep0.25cm
\item\label{tau:fun} As an operator on $V_1[\mu_1]\otimes V_2[\mu_2]$,  $\exp(\Rtau_{ij}(s))$ is a rational function of $s$, regular
at $s=\infty$ with value $\Id_{V_1[\mu_1]\otimes V_2[\mu_2]}$ at $s=\infty$.

\item\label{tau:unitary} $\Rtau_{ij}^{\op}(s) = \Rtau_{ji}(-s)$.

\item\label{tau:shift} $\tau_a\otimes \tau_b (\Rtau_{ij}(s)) = \Rtau_{ij}(s+a-b)$ for each  $a,b\in \C$. 

\item\label{tau:comm} Let $k\in\hbfI$, $n\in\N$ and let $y\in\Ysupkm{0}$.
Then, the following commutation relations hold, where $\sop$ is the shift
operator \eqref{eq:T-action}:
\begin{align*}
[\Rtau_{ij}(s), x_{k,n}^{\pm}\otimes y] &= 
\pm (\sop^{d_ia_{ik}} - \sop^{-d_ia_{ik}}) \cdot
\lp \sum_{a=0}^{\infty} x^{\pm}_{k,n+a} \otimes \partial_s^{(a)}(t_j(s))y\rp \\
&=\pm (\sop^{d_ia_{ik}} - \sop^{-d_ia_{ik}}) \cdot \\
&\hspace*{0.5in} \hbar\lp
\sum_{N\geqslant 0} N! s^{-N-1} \sum_{a+b=N} (-1)^a \frac{x_{k,n+a}^{\pm}}{a!}
\otimes \frac{t_{j,b}y}{b!} \rp, \\
[\Rtau_{ij}(s), y\otimes x_{k,n}^{\pm}] &=
\mp (\sop^{d_ja_{jk}} - \sop^{-d_ja_{jk}}) \cdot
\lp \sum_{b=0}^{\infty} y(-\partial_s)^{(b)}(t_i(-s))\otimes  x^{\pm}_{k,n+b}\rp\\
&=\pm (\sop^{d_ja_{jk}} - \sop^{-d_ja_{jk}}) \cdot \\
&\hspace*{0.5in}\hbar\lp
\sum_{N\geqslant 0} N! s^{-N-1} \sum_{a+b=N} (-1)^a \frac{t_{i,a}y}{a!}
\otimes \frac{x^{\pm}_{k,n+b}}{b!} \rp.
\end{align*}

\end{enumerate}
\end{prop}

\begin{pf}
We remind the reader of the following general fact, proved in Claims 1 and 2 of the
proof of \cite[Thm. 5.5]{sachin-valerio-III}:

Let $V,W$ be two finite--dimensional $\C$--vector spaces, $A:\C\to\End(V)$
and $B:\C\to\End(W)$ two rational functions, taking value $\Id$ at $\infty$,
such that $[A(s),A(s')]=0=[B(s),B(s')]$ for all $s,s'\in\C$. Let $\sigma(A)$
and $\sigma(B)$ be the set of poles of $A(s)^{\pm 1}$ and $B(s)^{\pm 1}$
respectively. Then the following is a rational $\End(V\otimes W)$--valued
function of $s$, taking value $\Id$ at $s=\infty$:
\[
X(s) = \exp\lp \oint_{C_1} A(u)^{-1}A'(u)\otimes \log(B(u+s))\, du\rp.
\]
Here $C_1$ is a contour enclosing
$\sigma(A)$, and $s$ is large enough so that $\log(B(u+s))$ is analytic
within and on $C_1$.
Moreover, we have 
\begin{equation*}
 X(s) = \exp\lp \oint_{C_2} \log(A(u-s))\otimes B(u)^{-1}B'(u)\, du\rp.
\end{equation*}
This proves \eqref{tau:fun} and \eqref{tau:unitary}. Note that \eqref{tau:unitary}
can also be easily deduced from the formal expansion \eqref{eq:tau-formal}.

To prove \eqref{tau:shift}, note that $\tau_a B_i(z) = e^{az}B_i(z)$. Therefore,
using \eqref{eq:tau-borel}, we have:
\[
\begin{aligned}
(\tau_a\otimes \tau_b)(\Rtau_{ij}(s)) &= \left.\tau_a(B_i(-z))\otimes \tau_b(B_j(z))\right|_{z=-\partial_s}
\cdot s^{-2} \\
&= \left.e^{-(a-b)z} B_i(-z)\otimes B_j(z)\right|_{z=-\partial_s}\cdot s^{-2} \\
&= e^{(a-b)\partial_s}\cdot \Rtau_{ij}(s) = \Rtau_{ij}(s+a-b).
\end{aligned}
\]
In the last line, we used Taylor's theorem $e^{c\partial_s}\cdot f(s) = f(s+c)$.

We remark that a ``contour deformation" style proof of \eqref{tau:comm}
is given in \cite[Prop. 5.10]{sachin-valerio-III}. Here,
we give a different ``formal" proof, using the expression \eqref{eq:tau-borel} of $\Rtau_{ij}(s)$
in terms of $B_i(z)$ and the commutation relation
\eqref{eq:comm-Bi}. Let us focus on the first relation (the second
one can be deduced easily from the first, using the unitarity relation
\eqref{tau:unitary}). Setting $c_{ik} = d_ia_{ik}\hbar/2$, we have
\begin{align*}
[\Rtau_{ij}(s),x_{k,n}^\pm\otimes y] &= \left. z[B_i(-z),x^\pm_{k,n}]\otimes B_j(z)y\right|
_{z=-\partial_s} \cdot s^{-1} \\
&= \left. \pm (e^{c_{ik}z}-e^{-c_{ik}z}) \left(\hbar\sum_{a,b\geqslant 0} (-1)^a \frac{x^\pm_{k,n+a}}{a!}
\otimes \frac{t_{j,b}y}{b!} z^{a+b}\right) \right|_{z=-\partial_s} \cdot s^{-1}\\
&= \pm (e^{-c_{ik}\partial_s} - e^{c_{ik}\partial_s}) \cdot X^\pm_{k,j;n}(s) (1\otimes y)\ ,
\end{align*}
where
\begin{equation}\label{eq:X-notation}
\begin{aligned}
X^\pm_{k,j;n}(s) &= \hbar\sum_{N\geqslant 0} \left(\sum_{a+b=N} (-1)^a \frac{x^\pm_{k,n+a}}{a!}
\otimes \frac{t_{j,b}}{b!} \right) (-\partial_s)^N \cdot s^{-1} \\
&= \hbar\sum_{N\geqslant 0} N!s^{-N-1} \left(\sum_{a+b=N} (-1)^a \frac{x^\pm_{k,n+a}}{a!}
\otimes \frac{t_{j,b}}{b!}\right)\\
&= \sum_{a=0}^\infty x_{k,n+a}^{\pm}\otimes \partial_s^{(a)}(t_j(s)).
\end{aligned}
\end{equation}
Since $\sop=e^{-(\hbar/2)\partial_s}$, the first relation of \eqref{tau:comm} follows immediately. \qedhere
\end{pf}

\subsection{The operator $\Rg(s)$}\label{ssec:Rg}
We now define
\begin{equation}\label{eq:Rg}
\Rg(s) \ceq   \sum_{i,j\in\hbfI} \bfB(\sop)^\ast_{ji}\cdot \Rtau_{ij}(s).
\end{equation}
The following is a direct corollary of Proposition \ref{pr:tau}
and the symmetry of $\bfB(\sop)^\ast$.

\begin{cor}\label{cor:Rg} The element $\Rg(s)$ has the following properties:

\begin{enumerate}[font=\upshape]\itemsep0.25cm
\item\label{Rg:1} $\exp(\Rg(s))$ is a rational function of $s$, taking value $1$ at $s=\infty$.
\item \label{Rg:2}$\Rg^{\op}(s) = \Rg(-s)$.
\item\label{Rg:3} $\tau_a\otimes \tau_b (\Rg(s)) = \Rg(s+a-b)$ for each $a,b\in \C$.
\item\label{Rg:4} Let $k\in\hbfI,n\in\N$ and let $y\in\Ysupkm{0}$. Then, we have the following
commutation relations:
\begin{align*}
[\Rg(s),x_{k,n}^\pm\otimes y] &= \pm
(\sop-\sop^{-1})\det(\bfB(\sop)) \cdot 
\sum_{a=0}^\infty x^{\pm}_{k,n+a}\otimes y\partial_s^{(a)}(t_k(s)),\\
[\Rg(s),y\otimes x_{k,n}^\pm] &=
\mp (\sop-\sop^{-1})\det(\bfB(\sop)) \cdot 
\sum_{b=0}^{\infty} y(-\partial_s)^{(b)}(t_k(s))\otimes x_{k,n+b}^{\pm}.
\end{align*}
\end{enumerate}
\end{cor}

\subsection{Expansion and regularization of $\Rg$}
\label{ssec:Rg-reg}

The following lemma computes the expansion of $\Rg(s)$ in $s^{-1}$,
and is crucial in carrying out a ``regularization" argument later.

\begin{lem}\label{lem:Rg-terms}
$\Rg(s)$ admits the expansion 
\begin{align*}
\hbar^{-2}\Rg(s) &=
\sczero \caff_0\otimes\caff_0 s^{-2}
+ 2\sczero \lp -\caff_1\otimes \caff_0 + \caff_0\otimes \caff_1\rp s^{-3} \\
&\hspace*{0.2in}
+ 6 \lp \sczero \lp \frac{\caff_2\otimes \caff_0}{2} -  \caff_1\otimes \caff_1 +
\frac{\caff_0\otimes\caff_2}{2}\rp + \frac{\hbar^2}{4} \ctzero \rp s^{-4} + O(s^{-5}),
\end{align*}
where $\sczero$ and $\ctzero$ are given as follows:
\begin{itemize}\itemsep0.25cm
\item $\sczero\in\Z_{>0}$ is a constant depending on $\gkm$ defined by
$\sczero a_ia_j = \bfB(1)^\ast_{i,j}$ for every $i,j\in\hbfI$
(see Table \ref{table} below). 

\item $\ctzero=\sum_{i,j} \bfB^{\ast,(2)}_{ij}
t_{i,0}\otimes t_{j,0}\in\hkm\otimes\hkm$, where $\bfB^{\ast,(2)}$ is the coefficient
of $t^2$ in the Taylor series expansion of $\bfB(e^t)^\ast$ near $t=0$:
\[
\bfB(e^t)^* = \bfB^\ast + t^2 \bfB^{\ast,(2)} + O(t^4).
\]
\end{itemize}

\end{lem}

\begin{pf}
Using the definition of $\Rg(s)$, and the formula \eqref{eq:tau-borel}
for $\Rtau_{ij}(s)$, we obtain
\begin{equation*}
\Rg(s) = \sum_{i,j} \bfB(\sop)^{\ast}_{ij} \cdot \Rtau_{ij}(s) = \left.\sum_{i,j} \bfB(e^{\frac{\hbar}{2}z})^* B_i(-z)\otimes B_j(z) \right|_{z=-\partial_s}
\cdot s^{-2}.
\end{equation*}
For notational simplicity, let us write
\[
\Rtau_{ij}^{(n)} \ceq   \sum_{a=0}^n (-1)^a \frac{t_{i,a}}{a!}\otimes
\frac{t_{j,n-a}}{(n-a)!}\ ,
\]
so that $B_i(-z)\otimes B_j(z) = \hbar^2\sum_{n\geqslant 0} \Rtau_{ij}^{(n)}z^n$.
With this notation at hand, we can write the expansion
of $\Rg(s)$ as
\begin{equation}\label{Rg-exp-proof}
\Rg(s) = \hbar^2 \sum_{N\geqslant 0} (N+1)!s^{-N-2} 
\lp
\sum_{n=0}^{\lfloor\frac{N}{2}\rfloor} \frac{\hbar^{2n}}{2^{2n}}
\sum_{i,j} \bfB^{*,(2n)}_{ij} \Rtau_{ij}^{(N-2n)}
\rp,
\end{equation}
where $\bfB^{*,(2n)}$ is the coefficient of $t^{2n}$ in the Taylor expansion of $\bfB(e^t)^\ast$ at $t=0$:
\[
\bfB(e^t)^\ast = \sum_{n=0}^{\infty} \bfB^{*,(2n)} t^{2n}.
\]
Thus, the coefficient of $s^{-2}$ in $\hbar^{-2}\Rg(s)$
is given by $\sum_{ij} \bfB^{*,(0)} t_{i,0}\otimes t_{j,0}$.
Note that $\bfB^{*,(0)} = \bfB^\ast$ is the adjoint of the corank $1$,
symmetric matrix $\bfB=D\hbfA$, and hence is of rank $1$. Thus, its rows (and columns, as
it is symmetric) are scalar multiples of $\caff_0$. In other words,
there is a constant $\sczero$ such that
\[
\bfB^*_{ij} = \sczero a_ia_j\, \text{ for every } i,j\in\hbfI\ .
\]
This shows that the coefficient of $s^{-2} $ in $\hbar^{-2}\Rg(s)$ is $\sum_{ij} \bfB^*_{ij}t_{i,0}\otimes t_{j,0} = \sczero\caff_0\otimes\caff_0$, as claimed in the statement of the lemma.

Similarly, we obtain from \eqref{Rg-exp-proof} that the coefficients of $s^{-3}$ and $s^{-4}$ in  $\hbar^{-2}\Rg(s)$ are given by 
\begin{equation*}
2 \sum_{ij} \bfB^*_{ij} \Rtau_{ij}^{(1)}
\quad \text{ and }\quad 
%
6 \sum_{ij} \bfB^*_{ij} \Rtau_{ij}^{(2)}+\frac{6\hbar^2}{4}\ctzero,
\end{equation*}
respectively. Since $\bfB^*_{ij} = \sczero a_ia_j$ for each $i,j\in \bfI$, we have
\begin{gather*}
2 \sum_{ij} \bfB^*_{ij} \Rtau_{ij}^{(1)}= 2 \sczero (\caff_0\otimes\caff_1 - \caff_1\otimes\caff_0),\\
6 \sum_{ij} \bfB^*_{ij} \Rtau_{ij}^{(2)}=6  \sczero \left(\frac{\caff_0\otimes
\caff_2}{2}- \caff_1\otimes
\caff_1+ \frac{\caff_2\otimes
\caff_0}{2}\right).
\end{gather*}
Combining these facts, we can conclude that the $s^{-3}$ and $s^{-4}$ coefficients of $\hbar^{-2}\Rg(s)$ are as stated in the lemma. 

We are left to compute the explicit values of the constant $\sczero$, which we list in Table \ref{table}. 

\begin{table}[h]
	\begin{tabular}{|c|c|c|c|c|c|c|c|c|c|}
		\hline
		Type of $\gkm$       & $\mathsf{A}_n^{(1)}$ & $\mathsf{B}_n^{(1)}$ & $\mathsf{C}_n^{(1)}$ &
		$\mathsf{D}_n^{(1)}$ & $\mathsf{E}_6^{(1)}$ & $\mathsf{E}_7^{(1)}$ & $\mathsf{E}_8^{(1)}$ & 
		$\mathsf{F}_4^{(1)}$ & $\mathsf{G}_2^{(1)}$ \\ 
		\hline
		$\sczero$ & $n+1$ & $2^n$ & $4$ & $4$ & $3$ & $2$ & $1$ & $4$ & $3$ \\
		\hline
	\end{tabular}
	
	\vspace*{0.1in}
	
	\begin{tabular}{|c|c|c|c|c|c|}
		\hline
		Type of $\gkm$ & $\mathsf{A}_{2n}^{(2)}$ & $\mathsf{A}_{2n-1}^{(2)}$ & 
		$\mathsf{D}_{n+1}^{(2)}$ & $\mathsf{E}_6^{(2)}$ & $\mathsf{D}_4^{(3)}$ \\
		\hline
		$\sczero$ & $2^n$ & $4$ & $2^n$ & $4$ & $3$\\
		\hline
	\end{tabular}
	\hspace*{1.275in}
	\vspace{0.25cm}
	\caption{Table of values of $\sczero$}\label{table} 
\end{table}

The values of $\sczero$ are computed as follows. Observe that for any fixed $i\in \hbfI$, we have 
\[
\sczero = \frac{\det(\bfB(1)^{(i)})}{a_i^2}, 
\]
where $\bfB(1)^{(i)}$ is the submatrix obtained by removing  the $i$--th
row and $i$--th column from $\bfB$. By choosing $i$
so that $a_i=1$ and $\bfB(1)^{(i)}$ is a (connected) finite--type Dynkin diagram, whose
determinants are known (see, for instance, \cite[\S 4.8, Table Fin]{kac}), we obtain the explicit values of $\sczero$.
This type-by-type computation completes the proof of the lemma. \qedhere
%

\end{pf}

Next, we define $\Rg_{\reg}(s)$ by setting
\begin{equation}\label{eq:Rg-reg}
\Rg_{\reg}(s) \ceq   \Rg(s) - 
\hbar^2\lp
\sczero \caff_0\otimes\caff_0 s^{-2}
+ 2\sczero \lp -\caff_1\otimes \caff_0 + \caff_0\otimes \caff_1\rp s^{-3}
\rp.
\end{equation}
Since the terms removed from $\Rg(s)$ are central (see Corollary
\ref{cor:c123} above), we have $\ad(\Rg(s)) = \ad(\Rg_{\reg}(s))$.

\begin{cor}\label{cor:Rg-reg}
Properties \eqref{Rg:2}--\eqref{Rg:4} of Corollary \ref{cor:Rg} hold for
$\Rg_{\reg}(s)$. In addition,
\[
\Rg_{\reg}(s) =
6 \hbar^2 \lp \sczero \lp \frac{\caff_2\otimes \caff_0}{2} -  \caff_1\otimes \caff_1 +
\frac{\caff_0\otimes\caff_2}{2}\rp + \frac{\hbar^2}{4} \ctzero \rp s^{-4} + O(s^{-5}),
\]
where $\sczero$ and $\ctzero$ are as in Lemma \ref{lem:Rg-terms}.
\end{cor}

\subsection{The difference equation and its
formal solution}\label{ssec:R0-diff}
Recall the difference equation \eqref{eq:L0-diff}:
\[
(\sop-\sop^{-1})\det(\bfB(\sop)) \cdot \Rl(s) = \Rg(s)\ .
\]
This is an irregular, additive difference equation. Meaning, the difference
operator has order of vanishing $3$, while the right--hand side is
$O(s^{-2})$. In more elementary terms, this equation has no solution
in $\Pseries{\Ysupkm{0}^{\otimes 2}}{s^{-1}}$, since the difference
operator applied to any such formal power series results in a series
starting with $s^{-4}$. The following simple lemma makes this more precise.

\begin{lem}\label{lem:diff-ops}
Let $\mathcal{A}$ be an arbitrary vector space over $\C$,
and let $F(s) = \sum_{n\geqslant 0} F_n s^{-n-1}\in\Pseries{\mathcal{A}}{s^{-1}}$.
Let $\hbar\in\nC$ and $\sop\cdot F(s) = F(s-\hbar/2)$, as above. Then,
for any polynomial $D(\sop)\in\C[\sop^{\pm 1}]$, we have
\[
D(\sop)\cdot F(s) = \sum_{N\geqslant 0} s^{-N-1} 
\lp
\sum_{\ell=0}^N \cbin{N}{\ell} F_{N-\ell} \frac{\hbar^\ell}{2^\ell}
\left.\lp (\sop\partial_\sop)^{\ell}\cdot D(\sop)\rp\right|_{\sop=1} 
\rp\ .
\]
\end{lem}

\begin{pf}
The proof of this lemma is a direct verification. Namely,
write $D(\sop) = \sum_{r\in\Z} d_r\sop^r$. Then, we have:
\begin{align*}
D(\sop)\cdot F(s) &= \sum_{n\geqslant 0} F_n s^{-n-1} 
\lp \sum_r d_r\lp 1-s^{-1}r\frac{\hbar}{2}\rp^{-n-1}\rp \\
&= \sum_{n\geqslant 0} F_n s^{-n-1} \lp
\sum_r d_r \sum_{\ell\geqslant 0} \cbin{n+\ell}{\ell}
r^\ell\frac{\hbar^\ell}{2^\ell} s^{-\ell}
\rp\\
&= \sum_{N\geqslant 0} s^{-N-1} \lp
\sum_{\ell=0}^N F_{N-\ell} \cbin{N}{\ell} \frac{\hbar^\ell}{2^\ell}
\lp\sum_r d_rr^\ell \rp
\rp\ .
\end{align*}
The lemma now follows from the observation that
\begin{equation*}
\sum_r d_rr^\ell = \left.(\sop\partial_\sop)^\ell\cdot D(\sop)
\right|_{\sop=1}. \qedhere
\end{equation*}
\end{pf}

\begin{rem}\label{rem:irregular}
The lemma above implies that if $D(\sop)$ has order of vanishing $k\in\N$
at $\sop=1$, then $D(\sop) \cdot F(s) \in s^{-k-1}\Pseries{\mathcal{A}}{s^{-1}}$,
for every $F(s)\in s^{-1}\Pseries{\mathcal{A}}{s^{-1}}$. In fact,
the operation is invertible, thus yielding an isomorphism of
vector spaces
\[
D(\sop) : s^{-1}\Pseries{\mathcal{A}}{s^{-1}}
\longisom s^{-k-1}\Pseries{\mathcal{A}}{s^{-1}}\ .
\]
\end{rem}

\begin{cor}\label{cor:Rl-formal}
There exists a unique $\Rl(s)\in s^{-1}\Pseries{\Ysupkm{0}^{\otimes 2}}{s^{-1}}$
such that
\begin{equation*}
(\sop-\sop^{-1})\det(\bfB(\sop)) \cdot \Rl(s) = \Rg_{\reg}(s)\ .
\end{equation*}
\end{cor}

\subsection{Properties of $\Rl(s)$}\label{ssec:Rl-props}

We now turn to establishing the key properties satisfied by the element $\Rl(s)$ from the previous corollary. Recall that $\Prim(\Ykm)$ is the linear span of $\{t_{i,r}\}_{i\in\hbfI,r\in\N}$.

\begin{prop}\label{pr:Rl}

$\Rl(s)$ has the following properties:

\begin{enumerate}[font=\upshape]\itemsep0.25cm
\item\label{Rl:prim} $\ds \Rl(s)\in \Pseries{\Prim(\Ykm)^{\otimes 2}}{s^{-1}}$.

\item\label{Rl:LT} The leading term of $\Rl(s)$ is given by
\[
\Rl(s) = s^{-1} \lp \frac{\sczero}{\qdzero\hbar} \lp \frac{\caff_2\otimes\caff_0}{2}
-\caff_1\otimes\caff_1 + \frac{\caff_0\otimes \caff_2}{2}\rp + 
\frac{\hbar}{4\qdzero} \ctzero \rp + O(s^{-2}),
\]
where $\sczero,\ctzero$ are as given in Lemma \ref{lem:Rg-terms} and 
$\qdzero$ is given in Lemma \ref{lem:order}.

\item\label{Rl:unitary} $
\Rl^{\op}(s) = -\Rl(-s)$.

\item\label{Rl:shift} $(\tau_a\otimes \tau_b)(\Rl(s)) = \Rl(s+a-b)$ for each $a,b\in \C$.

\item\label{Rl:comm} For each $k\in\hbfI$, $n\in\N$ and $y\in\Ysupkm{0}$, we have
the following commutation relations:
\begin{align*}
[\Rl(s),x_{k,n}^\pm\otimes y] &= \pm 
\sum_{a=0}^{\infty} x_{k,n+a}^{\pm}
\otimes \partial_s^{(a)}(t_k(s)) y, \\
[\Rl(s),y\otimes x_{k,n}^\pm] &=
\mp \sum_{b=0}^{\infty} y(-\partial_s)^{(b)}(t_k(-s))\otimes x^{\pm}_{k,n+b}.
\end{align*}

\end{enumerate}
\end{prop}

\begin{pf}
Note that \eqref{Rl:prim} is obvious from the construction and \eqref{Rl:shift}
follows directly from Part \eqref{tau:shift} of Proposition \ref{pr:tau}.

Let us prove \eqref{Rl:LT}. Let $\Rl_0$ denote the coefficient
of $s^{-1}$ in $\Rl(s)$. In our case, the difference operator is
\[
\Rd(\sop) = (\sop-\sop^{-1})\det(\bfB(\sop)) = (\sop-\sop^{-1})^3 \frac{\det(\bfB(\sop))}{(\sop-\sop^{-1})^2}.
\]
Let $\ds C(\sop)= \frac{\det(\bfB(\sop))}{(\sop-\sop^{-1})^2}$. The following is an
easy computation:
\[
\left. (\sop\partial_\sop)^{\ell}\cdot \Rd(\sop)\right|_{\sop=1} = 
\sum_{j=0}^{\ell} \cbin{\ell}{j} 3(3^{j-1}-1)(1+(-1)^{j-1}) 
\left.\left((\sop\partial_\sop)^{\ell-j}\cdot C(\sop)\right)\right|_{\sop=1}\ .
\]
The relevant coefficient is the $\ell=j=3$ term, where we get $48\qdzero$.
Using Lemma \ref{lem:diff-ops}, we compare the coefficients
of $s^{-4}$ in $\Rd(\sop)\cdot\Rl(s)$ and $\Rg_{\reg}(s)$ from Lemma \ref{lem:Rg-terms} 
to obtain
\[
48 \frac{\hbar^3}{8} \qdzero \Rl_0 = 
6 \hbar^2\lp \sczero \lp \frac{\caff_2\otimes \caff_0}{2} -  \caff_1\otimes \caff_1 +
\frac{\caff_0\otimes\caff_2}{2}\rp + \frac{\hbar^2}{4} \ctzero \rp,
\]
which is precisely \eqref{Rl:LT}.

For \eqref{Rl:unitary}, we flip 
the tensor factors in the difference equation \eqref{eq:Rl-diff} to get
\[
\Rd(\sop)\cdot \Rl^{\op}(s) = \Rg_{\reg}(-s).
\]
On the other hand, $\sop\cdot f(-s) = (\sop^{-1}\cdot f(s))|_{s\mapsto -s}$, and
$\Rd(\sop^{-1}) = -\Rd(\sop)$. This gives:
\[
\Rd(\sop)\cdot \Rl(-s) = -\Rg_{\reg}(-s).
\]
Thus, by uniqueness, we have  $\Rl^{\op}(s) = -\Rl(-s)$.

The proof of \eqref{Rl:comm} is also based on a uniqueness argument. Namely,
apply $\Rd(\sop)$ to both sides of the commutation relation. The resulting
equation holds by Corollaries \ref{cor:Rg} and \ref{cor:Rg-reg}.
Thus both sides solve the same difference equation and hence must be equal.
\end{pf}

\subsection{Growth properties of the Borel transform of $\Rg_{\reg}(s)$}
\label{ssec:borel-growth}

Recall from Corollary \ref{cor:Rg} that $\exp(\Rg(s))$ is a rational function of $s$, 
taking value $1$ at $s=\infty$. In particular, $\Rg(s)$ and $\Rg_{\reg}(s)$
are holomorphic in a neighbourhood of $\infty$.
Thus, the following general fact applies to our situation,
which is essential to take the Laplace transform and use Watson's lemma
(see Theorem \ref{athm:watson} below).

\Omit{\noindent {\bf Claim.}
The matrix entries of $\Rg(s)$, and hence of $\Rg_{\reg}(s)$, 
are of the form
\[
\log\lp \prod_j \frac{s-a_j}{s-b_j}\rp + r_0(s),
\]
where $\{a_j,b_j\}_j\subset\C$ is a finite set of complex numbers and
$r_0(s)$ is a rational function vanishing at $s=\infty$.

\begin{pf}[Proof of Claim]Write $X(s) = \exp(\Rg(s))$ in its mulitplicative
Jordan decomposition, whose entries are again rational functions
of $s$ (see \cite[Lemma 4.12]{sachin-valerio-2}): $X(s) = X_d(s)(1+X_n(s))$.
Note that entries of $X_d(s)$ are rational functions taking
value $1$ at $s=\infty$, while those of $X_n(s)$ vanish
at $s=\infty$. Hence,
\[
\Rg(s) = \log(X(s)) = \log(X_d(s)) + \sum_{r\geqslant 0} (-1)^r \frac{X_n(s)^{r+1}}{r+1}.
\]
Here, the second sum on the last line is finite, since $X_n(s)$ is nilpotent.
Our claim follows.
\end{pf}
}

\begin{lem}\label{lem:borel-growth}
If $a(s)=\sum_{n=0}^{\infty} a_ns^{-n-1}\in s^{-1}\Pseries{\C}{s^{-1}}$ is convergent, then
its Borel transform $\Borel{a}(t) = \sum_{n=0}^{\infty} a_n \frac{t^n}{n!}$
is an entire function of sub--exponential growth.
That is, there exist constants $M,R\in\R_{>0}$ such that
\[
\left|\Borel{a}(t)\right| \leqslant Me^{R|t|},
\text{ for all } t\in\C\ .
\]
\end{lem}

\begin{pf}
By Cauchy--Hadamard's theorem, $a(s)$ having non--zero radius of convergence
implies that its coefficients grow geometrically. That is,
there exist constants $M,R$ such that $|a_n|\leqslant M R^n$,
for all $n\geqslant 0$. Thus, we obtain by the triangle inequality
\[
\left|\Borel{a}(t)\right| \leqslant \sum_{n=0}^\infty M\frac{R^n|t|^n}{n!}
= Me^{R|t|}
\]
as claimed.
\end{pf}

\Omit{
\begin{cor}\label{cor:borel-growth}
Let $g(s)$ be a matrix entry of $\Rg_{\reg}(s)$.
Write $g(s) = \sum_{n=0}^{\infty} g_n s^{-n-1}$ and let
$ \Borel{g}(t) = \sum_{n=0}^{\infty} g_n \frac{t^n}{n!}$
be its Borel transform. Then, $\Borel{g}(t)$ is an entire
function of $t$, and there exist constants $C_1,C_2,R$
such that:
\[
\left|\Borel{g}(t)\right| \leqslant C_1 e^{C_2|t|},
\text{ for every $t$ with } |t|>R\ .
\]
\end{cor}
}


\appendix

\section{Laplace transform and regular difference equations}\label{app:LapDE}

We collect some of the well--known techniques to solve linear,
additive, regular difference equations. The material of this
section is fairly standard, and can be found in any advanced
text on complex analysis, for instance \cite{ablowitz-fokas,
costin, whittaker-watson}.

\subsection{Set up and statement of the main theorem}
\label{appA:setup}

Let $D(\sop)\in\C[\sop^{\pm 1}]$ and let $k\in\N$ be its order of vanishing
at $\sop=1$. Let $g(s)$ be a meromorphic, $\C$--valued function of
$s\in\C$, which is holomorphic in a neighbourhood of $s=\infty$, and has
order of vanishing $k+1$ there. Thus,
$g(s) = \sum_{n=k}^{\infty} g_n s^{-n-1}$. Let $\Borel{g}$ denote
the Borel transform of $g(s)$:
\[
\Borel{g}(t) = \sum_{n=k}^{\infty} g_n \frac{t^n}{n!}\ .
\]

Let us fix $\chi\in\nC$ a non--zero step and let $\sop$ act as shift
$\sop\cdot F(s) = F(s-\chi)$. Let $\theta=\arg(\chi)\in (-\pi,\pi]$.\\

Below, we use the following notation for rays and half--planes. Let
$\psi\in\R$ and $c\in\R_{>0}$. Let $\ell_\psi \ceq   \R_{\geqslant 0}e^{\iota\psi}$
be the ray at phase $\psi$. Let
\[
\mathbb{H}_{\psi,c} \ceq   \{z\in\C : \Re(ze^{-\iota\psi})>c\}\ ,
\]
be the half--plane orthogonal to $\ell_\psi$, located to the right
of the line perpendicular to $\ell_\psi$, passing through $ce^{\iota\psi}$.

\begin{figure}[h]
\includegraphics[height=1in]{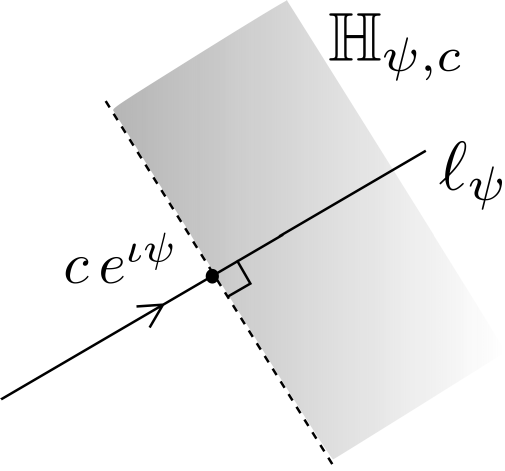}
\caption{The ray $\ell_\psi$ and corresponding half--plane $\mathbb{H}_{\psi,c}$}
\end{figure}

According to Lemma \ref{lem:diff-ops},
we have a unique formal power series $f(s)\in s^{-1}\Pseries{\C}{s^{-1}}$,
such that $\ds D(\sop)\cdot f(s) = \sum_{n=k}^\infty g_ns^{-n-1}$.

\begin{thm}\label{appA:thm}
Assume that the non--zero roots of $D(\sop)=0$ lie on the unit circle.
Then, there exist two meromorphic functions $f^\eta(s)$, $\eta\in\{\uparrow,\downarrow\}$,
uniquely determined by the following conditions.
\begin{enumerate}[font=\upshape]\itemsep0.25cm
\item\label{fn:diff} $D(\sop)\cdot f^\eta(s) = g(s)$.

\item\label{fn:hol} $f^\eta(s)$ is holomorphic in $\mathbb{H}^\eta$, where
$\mathbb{H}^{\uparrow} = \{s : \Re(s/\chi)\gg 0\} = -\mathbb{H}^{\downarrow}$.

\item\label{fn:asym} $f^\eta(s)\sim f(s)$ as $\pm\Re(s/\chi)\to \infty$.
\end{enumerate}

Moreover, $f^\eta(s)$ is holomorphic in a larger domain $\PP^\eta_\delta$,
for any $\delta>0$, where,
\[
\PP^{\uparrow}_\delta \ceq   \bigcup_{\psi\in \lp \theta-\frac{\pi}{2}+\delta,
\theta+\frac{\pi}{2}-\delta\rp} \mathbb{H}_{\psi,C}\ , \qquad
\PP^{\downarrow}_{\delta} = -\PP^{\uparrow}_{\delta}\ .
\]
Here $C$ is the constant for which $|\Borel{g}(t)|\leqslant Me^{C|t|}$,
as in Lemma \ref{lem:borel-growth} above.

The asymptotic expansion $f^\eta(s)\sim f(s)$ is valid in a larger sector 
$\Sigma^\eta_\delta$, for any $\delta>0$, 
\[
\Sigma^{\uparrow}_\delta \ceq   \{re^{\iota\psi} : r\in\R_{>0}, \psi \in \lp \theta-\pi+\delta,
\theta+\pi-\delta\rp \} , \qquad
\Sigma^{\downarrow}_{\delta} = -\Sigma^{\uparrow}_{\delta}\ .
\]

(see Figures \ref{afig:domain} and \ref{afig:sector} below).
\end{thm}

\begin{figure}[h]
\includegraphics[height=2in]{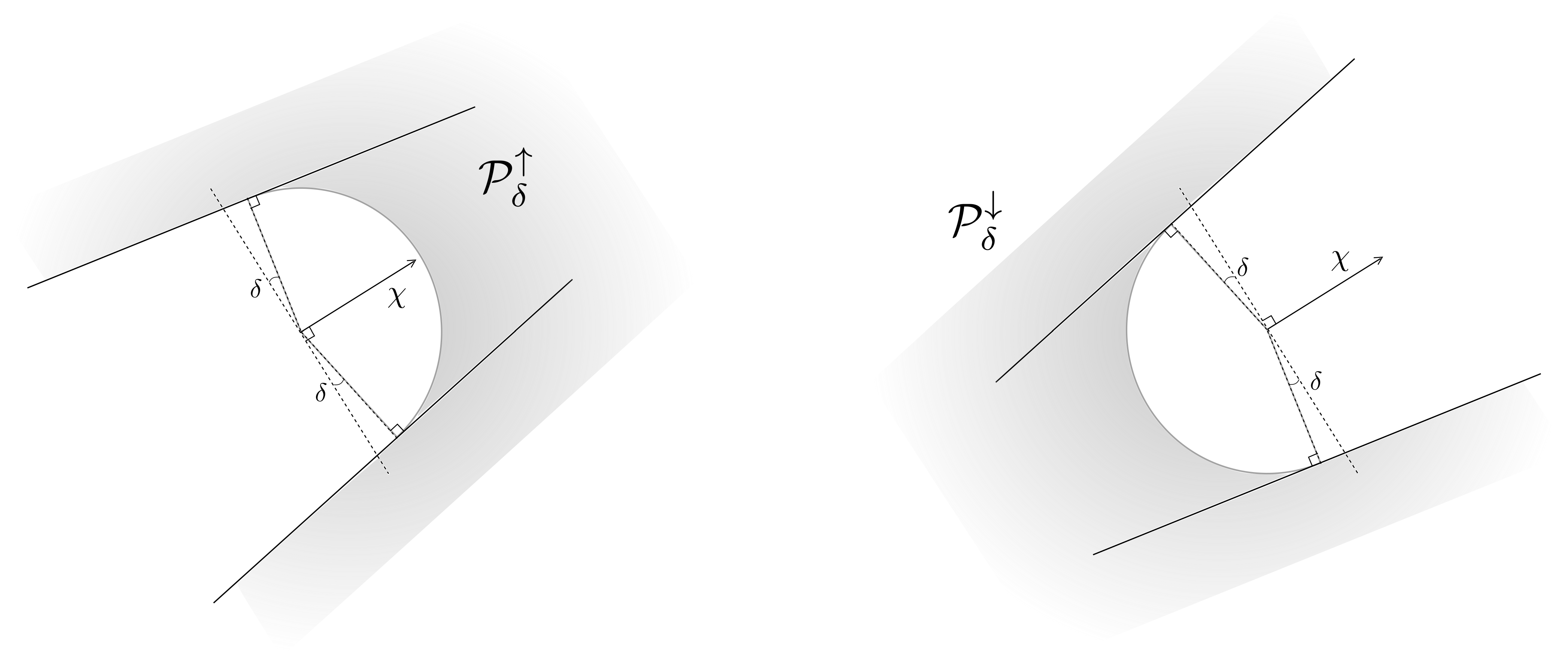}
\caption{Domains $\mathcal{P}^{\eta}_\delta$ of holomorphy}
\label{afig:domain}
\end{figure}

\begin{figure}[h]
\includegraphics[height=2in]{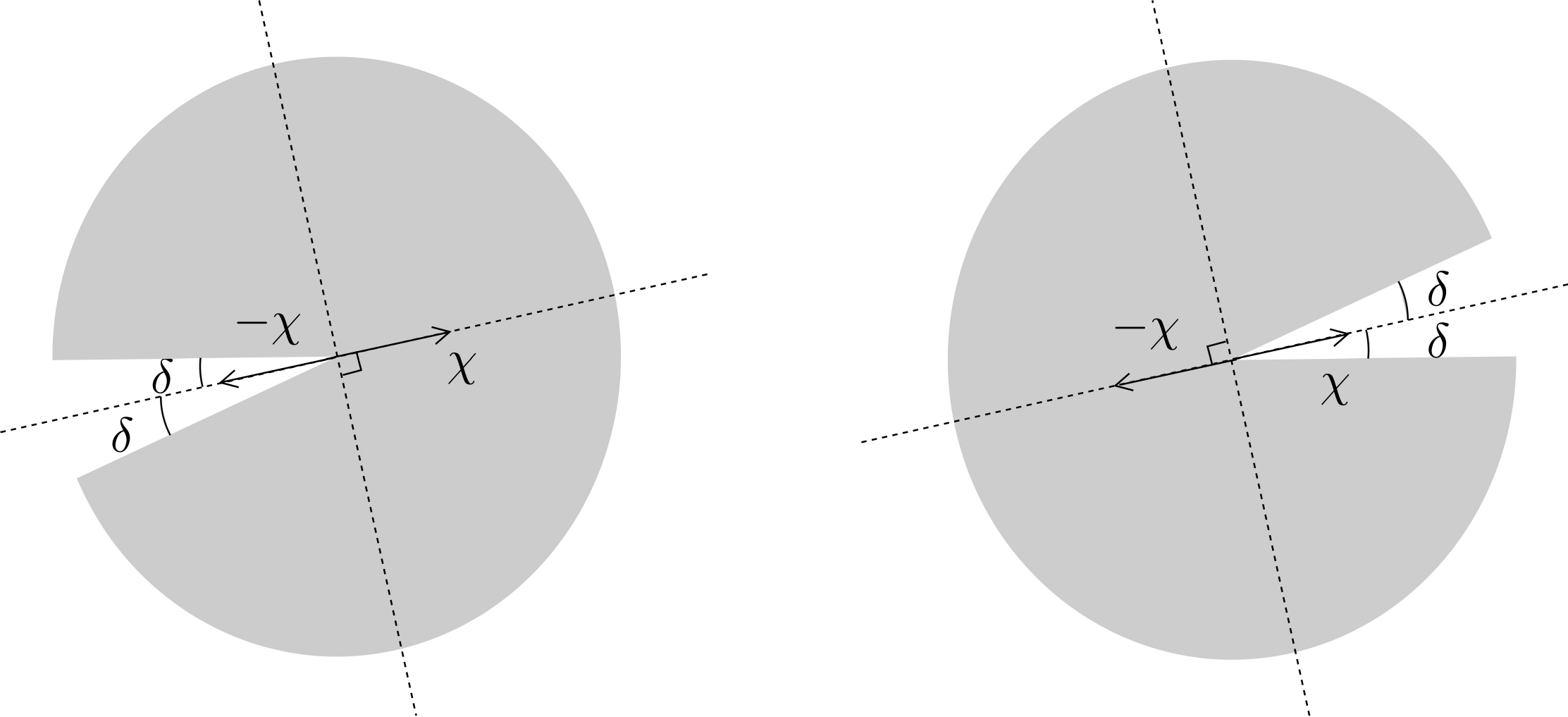}
\caption{Sectors $\Sigma^{\eta}_{\delta}$ of asymptotic expansion}
\label{afig:sector}
\end{figure}

\begin{pf}
For notational convenience, we set $\chi=1$. The reader can easily
verify that the statement of the theorem can be obtained from
its $\chi=1$, $\theta=0$ counterpart, by a counterclockwise rotation
by $\theta$.\\

The proof is given in the rest of this section. We show
uniqueness of $f^\eta$ in Section \ref{appA:uniqueness}.
The existence is based on a general technique of Laplace transforms
and their asymptotic expansions obtained via Watson's lemma. We
review these results in Section \ref{appA:watson}, and use them
to show the existence of $f^{\eta}$ in Section \ref{appA:existence}.
We show that the domain of holomorphy and sector of validity
of asymptotic expansion can be enlarged, as stated in the
theorem, in Section \ref{appA:doms}.
\end{pf}

\subsection{Uniqueness}\label{appA:uniqueness}
The uniqueness of $f^\eta(s)$ follows from the following
general lemma. For this, we drop the hypothesis that
the roots of $D(\sop)=0$ lie on the unit circle, as it appears
naturally from the conclusion of the lemma. We will assume,
without loss of generality, that $0$ is not a root of $D(\sop)$.

Let $\ol{\rho}(D)$ (resp. $\ul{\rho}(D)$) be the modulus of the longest 
(resp. shortest) root of $D(\sop)=0$. Note that $0<\ul{\rho}(D)\leqslant \ol{\rho}(D)$.

\begin{lem}\label{alem:uniqueness}
Let $\Sigma\subset\C$ be an unbounded open set satisfying:
\begin{gather*}
\text{For every } z\in\C, \text{ there exists } N\gg 0,
\text{ such that } \\ 
z+n\in\Sigma,\ \text{ (resp. $z-n\in\Sigma$)}\ \ , 
\forall\ n\geqslant N\ .
\end{gather*}
Assume that $\phi(s)$ is a meromorphic function of $s$, holomorphic
on $\Sigma$ such that:
\begin{enumerate}\itemsep0.25cm
\item $\phi(s)$ is asymptotically zero in $\Sigma$. That is, 
\[
\text{For every } n\in\N,\ 
\lim_{\begin{subarray}{c} s\to\infty \\ s\in\Sigma \end{subarray}} s^n\phi(s) = 0\ .
\]

\item $D(\sop)\cdot \phi(s) = 0$.
\end{enumerate}

If $\ol{\rho}(D)\leqslant 1$ (resp. $\ul{\rho}(D)\geqslant 1$), then $\phi \equiv 0$.
\end{lem}

\begin{pf}
For the purposes of the proof, we write $D(\sop) = \sop^N - \sum_{r=0}^{N-1}d_r \sop^r$,
and $d_0\neq 0$ (up to an overall shift, this is the general case).
We work with the right fundamental domain.
Thus, for every $s\in \C$,
\[
\phi(s) = \sum_{n=0}^{N-1} d_n \phi(s+N-n)
\]
In vector notation $\vec{\Phi}(s) = (\phi(s)\ \phi(s+1)\ \cdots \ \phi(s+N-1))^T$,
we get:
$\ds
\vec{\Phi}(s) = \mathcal{D}\cdot \vec{\Phi}(s+1)$, where,

\[
\mathcal{D} = \left[
\begin{array}{ccccc}
d_{N-1} & d_{N-2} & \cdots & \cdots & d_0 \\
1 & 0 & \cdots & \cdots & 0 \\
0 & 1 & 0 & \cdots & 0 \\
\vdots & \ddots & \ddots & \ddots & \vdots \\
0 & \cdots & 0 & 1 & 0 
\end{array}
\right]
\]

is the companion matrix of $D(\sop)$. The following identity is known
as Gelfand's formula. Its proof can be found in any standard functional
analysis book, for instance \cite[\S 149]{riesz-nagy}. 
Here $|\cdot|$ is an arbitrary, fixed norm on the space of $N\times N$
matrices:

\[
\lim_{n\to \infty} |\mathcal{D}^n|^{\frac{1}{n}} = \ol{\rho}(D)\ .
\]

We will assume that $\ol{\rho}(D)\leqslant 1$. For the left domain, we have
to use the fact that $\mathcal{D}$ is invertible ($d_0\neq 0$), thus,
$|\mathcal{D}^{-n}|$ grows as $\ul{\rho}(D)^{-n}$.\\

The remainder of the argument is standard. Let $s_0\in\C$ be fixed,
and assume that we are given $\varepsilon>0$. We will show that
$|\vec{\Phi}(s_0)|<\varepsilon$, hence showing that it has to be zero.\\

Choose an $\ell\geqslant 1$ and $\alpha>0$. Using the asymptotically zero condition,
we get
$R>0$ such that $|\vec{\Phi}(s)|<\alpha |s|^{-\ell}$ for every
$s\in\Sigma$ with $|s|>R$. Using the hypothesis on $\Sigma$,
we choose an $N_0>0$ so that $s_0+n\in\Sigma$ and $\Re(s_0+n)>R$, for
every $n\geqslant N_0$. The following inequality follows for every
$n\geqslant N_0$:

\[
\begin{aligned}
\left|\vec{\Phi}(s_0)\right| &= \left| \mathcal{D}^n \cdot \vec{\Phi}(s_0+n)\right|
\leqslant |\mathcal{D}|^n \alpha |s_0+n|^{-\ell}\ .
\end{aligned}
\]
As $|\mathcal{D}|^n$ grows proportionally to $\ol{\rho}(D)^n$, it is
enough to observe that
\[
\lim_{n\to\infty} \rho^n |s_0+n|^{-\ell} = 
\begin{array}{cl}
0 & \text{ if } \rho\leqslant 1\ ,\\
\infty & \text{ otherwise}.
\end{array}
\]
\end{pf}

\begin{rem}
It is worth mentioning that the uniqueness fails in general. For instance,
consider the equation $f(t+1) - 2f(t) = 0$. It has infinitely many
solutions $\lambda 2^{t}$ ($\lambda\in\C$) which are all 
asymptotically zero as $\Re(t)\to-\infty$.
\end{rem}

\subsection{Laplace's theorem and Watson's lemma}\label{appA:watson}

The following theorem summarizes the fundamental results of
Laplace transforms and their asymptotic expansions.
For a proof, see \cite[\S 6.2.2]{ablowitz-fokas}
or \cite[\S 3.3, 3.4]{costin}.

\begin{thm}\label{athm:watson}
Assume that $\psi\in\R$ and $\ell_{-\psi} = \R_{>0}e^{-\iota\psi}$
is the ray at phase $-\psi$. Assume given a continuous
function $F(t)$, $t\in\ell_{-\psi}$ such that
\begin{itemize}\itemsep0.25cm
\item $|F(t)|$ grows slower than exponential as $t\to\infty$.
That is, there are constants $R,C_1,C_2\in\R_{>0}$ such that
\[
|F(t)| < C_1 e^{C_2|t|}, \text{ for } |t|>R,\ t\in\ell_{-\psi}\ .
\]
\item $F(t)$ has at worst logarithmic singularity as $t\to 0$.
That is, there are constants $r,c_1,c_2\in\R_{>0}$ such that
\[
|F(t)|<c_1 t^{c_2-1}, \text{ for } |t|<r,\ t\in\ell_{-\psi}\ .
\]
\end{itemize}

Then, the following formula defines a function of $z$, holomorphic
on the half plane $\mathbb{H}_{\psi,C_2}$:
\[
\Laplace{F}_\psi(z) \ceq   \int_{\ell_{-\psi}} F(t) e^{-tz} dt\ .
\]

Moreover, if $F(t)\sim \sum F_n \frac{t^n}{n!}$ as $t\to 0$ along
$\ell_{-\psi}$, then we have the following asymptotic
expansion as $\Re(ze^{-\iota\psi})\to\infty$:
\[
\Laplace{F}_\psi(z) \sim \sum_{n=0}^{\infty} F_n z^{-n-1}\ .
\]
This expansion remains valid as $|z|\to\infty$ in
\[
z\in\Sigma^\delta_{\psi} \ceq   \{re^{\iota t} : r\in\R_{>0},
t\in (\psi-\pi/2+\delta,\psi+\pi/2-\delta)\}\ ,
\]
for any $\delta>0$.
\end{thm}

\begin{figure}[h]
\includegraphics[height=1.5in]{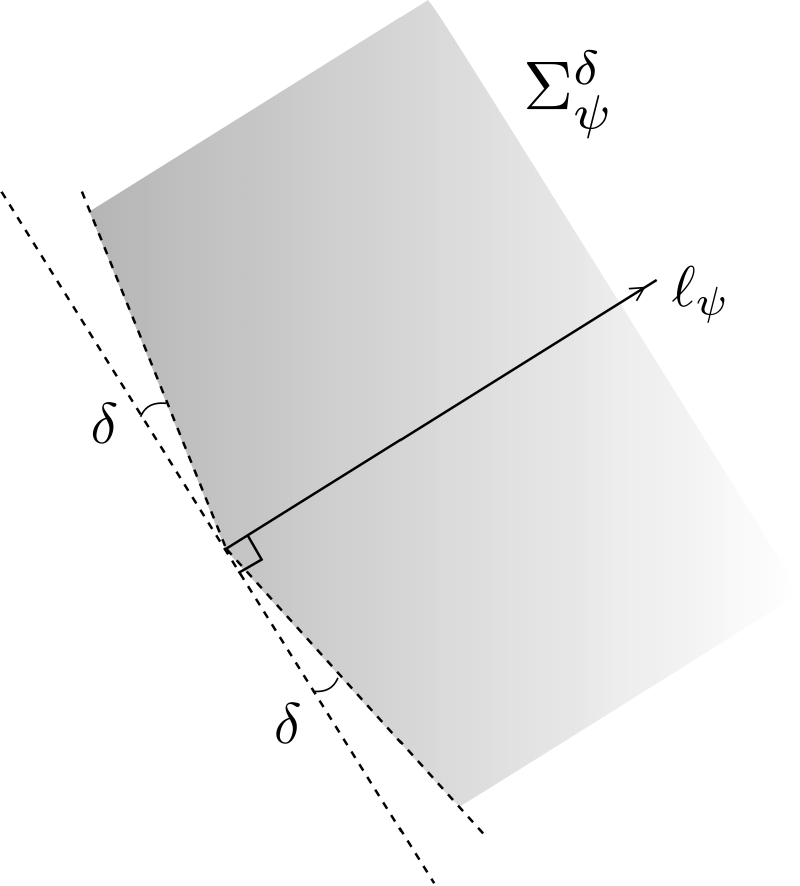}
\caption{Sector $\Sigma_\psi^{\delta}$}
\label{afig:asym}
\end{figure}

\begin{rem}
Note that the last part of the theorem is a triviality, since
\[
|s|\to\infty,\ s\in\Sigma_\psi^{\delta} \text{ if, and only if }
\Re(se^{-\iota\psi})\to\infty\ .
\]
\end{rem}

\subsection{Existence}\label{appA:existence}

Let $\psi\in\R$, and define:
\[
f_\psi(s) \ceq   \int_{\ell_{-\psi}} \frac{\Borel{g}(t)}{D(e^t)} e^{-ts}\, dt\ .
\]

We have to verify that the integrand satisfies the
conditions of Theorem \ref{athm:watson} above, as long as
$\ell_{-\psi}$ does not contain any roots of $D(e^t)=0$.
As roots of $D(\sop)=0$ are on the unit circle, this excludes
rays with $\psi \in \pm \frac{\pi}{2}+2\pi\Z$.\\

As $t\to\infty$ along $\ell_{-\psi}$, either $\Re(t)\to\infty$
or $-\Re(t)\to\infty$. It is not difficult to see that
$|D(e^t)|$ approaches a finite, positive number, or $\infty$,
as $\pm \Re(t)\to\infty$. In either case, we can choose $c>0$
and $R>0$, such that $|D(t)|>c$ for every $t\in\ell_{-\psi}$,
with $|t|>R$. Combined with our hypothesis on $\Borel{g}(t)$,
the sub--exponential growth condition is met for our kernel.
We note that, by our hypothesis on the order of vanishing,
$\frac{\Borel{g}(t)}{D(e^t)}$ is holomorphic near $t=0$.\\

Thus for $\psi\not=\pm\frac{\pi}{2}$ modulo $2\pi$, we get a holomorphic function
on $\mathbb{H}_{\psi,C_2}$. Moreover,
\[
D(\sop)\cdot f_\psi(s) = \int_{\ell_{-\psi}} \frac{\Borel{g}(t)}{D(e^t)} D(e^t) e^{-ts}\, dt
= \int_{\ell_{-\psi}} \Borel{g}(t) e^{-ts}\, dt = g(s)\ .
\]

This functional equation allows us to extend $f_\psi(s)$
as a meromorphic function of $s\in\C$. 

Moreover, as $|s|\to\infty$,
$s\in\Sigma^\delta_\psi$ (see Figure \ref{afig:asym}), 
we get the following asymptotic expansion of $f_\psi$:
let $\ds \frac{\Borel{g}(t)}{D(e^t)} = \sum_{n=0}^{\infty} \beta_n \frac{t^n}{n!}$
be its Taylor series expansion. Then,
\[
f_\psi(s) \sim \sum_{n=0}^{\infty} \beta_n s^{-n-1}\ .
\]
By uniqueness of the formal solution to $D(\sop)\cdot f(s)=g(s)$, 
the series on the right--hand side is $f(s)$ for all $\psi$.\\

Hence, we obtain $f^\uparrow$ for $\psi=0$ and $f^{\downarrow}$
for $\psi=\pi$. This prove that existence of the claimed meromorphic
solutions.

\subsection{Extension of domains}\label{appA:doms}
We will now show that $f_{\psi_1}(s) = f_{\psi_2}(s)$ for any $\psi_1,\psi_2\in (-\pi/2,\pi/2)$.
Similarly, for $\psi_1,\psi_2\in (\pi/2,3\pi/2)$.
Therefore, $f^{\uparrow} = f_\psi$ for any $\psi\in (-\pi/2,\pi/2)$,
and $f^{\downarrow} = f_\psi$ for any $\psi\in (\pi/2,3\pi/2)$.
This shows that $f^\uparrow$ is holomorphic in $\PP^\uparrow_\delta$,
and $f^\uparrow\sim f$ is valid in $\Sigma^\uparrow_\delta$ for any
$\delta>0$.\\

Now, let $\psi_1,\psi_2\in (-\pi/2,\pi/2)$, $\psi_1<\psi_2$.
As $f_{\psi_1}(s)=f_{\psi_2}(s)$ is an identity between two meromorphic 
functions solving the same difference
equation, it is enough to verify it for $s\in\mathbb{H}_{\psi_1,C_2}\cap \mathbb{H}_{\psi_2,C_2}$.
Note that this set satisfies the hypotheses of Lemma \ref{alem:uniqueness}.

Fix a constant $C_3>C_2$, and assume that 
$s\in \mathbb{H}_{\psi_1,C_2}\cap \mathbb{H}_{\psi_2,C_2}$ is such
that $\Re(se^{-\iota\theta})\geqslant C_3$ for every $\theta\in [\psi_1,\psi_2]$.

Then, we get:
\[
f_{\psi_1}(s) - f_{\psi_2}(s) = \int_{\ell_{-\psi_1}} - \int_{\ell_{-\psi_2}} 
\lp \frac{\Borel{g}(t)}{D(e^t)} e^{-ts}\rp\, dt\ .
\]

Let $R>0$ and let $\mathcal{C}_R$ be the closed contour consisting of three smooth
components: (a) straight (directed) line segment $L_{1,R}$ from $0$ to $Re^{-\iota\psi_1}$,
(b) circular arc $\gamma_R$ from $Re^{-\iota\psi_1}$ to $Re^{-\iota\psi_2}$, and (c)
directed segment from $Re^{-\iota\psi_2}$ to $0$, denoted by $-L_{2,R}$. 

\begin{figure}[h]
\includegraphics[height=1in]{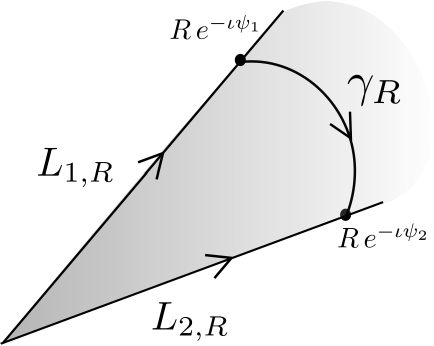}
\caption{The contour $\mathcal{C}_R$}
\end{figure}

As the integrand is
holomorphic on the right half--plane, by Cauchy's theorem:
\[
\int_{\mathcal{C}_R} \frac{\Borel{g}(t)}{D(e^t)} e^{-ts}\, dt = 0
\]

Therefore, we get:

\[
\begin{aligned}
f_{\psi_1}(s) - f_{\psi_2}(s) &= \lim_{R\to\infty} \int_{L_{1,R}} - \int_{L_{2,R}}
\lp \frac{\Borel{g}(t)}{D(e^t)} e^{-ts}\, dt\rp \\
&= \lim_{R\to\infty}
\int_{\mathcal{C}_R} - \int_{\gamma_R} 
\lp \frac{\Borel{g}(t)}{D(e^t)} e^{-ts}\, dt\rp \\
&= -\lim_{R\to\infty} \int_{\gamma_R} \frac{\Borel{g}(t)}{D(e^t)} e^{-ts}\, dt
\end{aligned}
\]

Thus, it remains to show that the last limit is zero. Recall that
we have the following bound on our kernel: there exist constants
$C_1,C_2$ and $R_1$ such that

\[
\left|  \frac{\Borel{g}(t)}{D(e^t)} \right| < C_1 e^{C_2|t|},
\text{ for every $t$ with } |t|>R_1\ .
\]

And, $s\in\C$ is chosen so that $|e^{-st}| = e^{-\Re(st)} \leqslant e^{-C_3R}$
($C_3>C_2$ was a fixed constant).
Hence, for $R\gg 0$, we obtain:

\[
\left| \int_{\gamma_R} \frac{\Borel{g}(t)}{D(e^t)} e^{-ts}\, dt\right|
\leqslant C_1 e^{(C_2-C_3)R} (\psi_2-\psi_1)R.
\]

The last quantity clearly approaches $0$, as $R\to \infty$, and
we are done.

\section{The augmented symmetrized Cartan matrix}\label{app:Augmented}

In this appendix, we give a full proof of Proposition \ref{P:full-rank}, which states  that the augmented matrix $(\bfB\,|\, \underline{\mu})$ has full rank,
where $\underline{\mu}$ is the column vector with $j$-th component $\mu_j=\sum_{i\in \hbfI} a_i (d_ia_{ij})^3$ (see  Corollary \ref{cor:c123}) and $\bfB=(d_i a_{ij})_{i,j\in \hbfI}$ is the symmetrized Cartan matrix of $\gkm$.

Recall that the integers $(a_j)_{j\in\bfI}$ are listed explicitly in \cite[Ch. 4, Tables Aff 1-2-3]{kac}. For convenience, we list below the values of the symmetrizing integers $(d_j)_{j\in\bfI}$, following \emph{loc. cit.} for the ordering of $\bfI$.

\begin{table}[h]
	\begin{tabular}{|c|c|c|c|}
		\hline
		Type of $\gkm$     & $(d_j)_{j\in\bfI}$ & Type of $\gkm$     & $(d_j)_{j\in\bfI}$\\
		\hline\hline
		 $\mathsf{A}_\ell^{(1)}, \mathsf{D}_\ell^{(1)}, \mathsf{E}^{(1)}_{6,7,8}$ & (1,\dots, 1) & $\mathsf{A}_{2\ell}^{(2)}$ & (1,2,\dots,2,4)\\
		\hline
		$\mathsf{B}_{\ell}^{(1)}$ & (2,\dots,2,1) & $\mathsf{A}_{2\ell-1}^{(2)}$  & (1,\dots,1,2)\\
		\hline
		$\mathsf{C}_\ell^{(1)}$ & (2,1,\dots,1,2) & $\mathsf{D}_{\ell+1}^{(2)}$& (1,2,\dots,2,1)\\
		\hline
		 $\mathsf{F}_4^{(1)}$ & (2,2,2,1,1)& $\mathsf{E}_6^{(2)}$& (1,1,1,2,2)\\
		\hline
		$\mathsf{G}_2^{(1)}$ & (3,3,1) & $\mathsf{D}_4^{(3)}$ & (1,1,3)\\
		\hline
	\end{tabular}
	\vspace{0.25cm}
	\caption{Symmetrizing integers}\label{table-d}
\end{table}

\subsection{A characterization of $\mu_j$}\label{appB-mu_j}

 We begin with the following simple lemma. 
\begin{lem}\label{L:full-rank}
For each $i,j\in \hbfI$, set $m_{ji}\ceq  a_{ji}(1-a_{ji}^2)$. Then 
\begin{equation*}
\mu_j/d_j^3=6a_j-\sum_{\substack{i\in \hbfI\\a_{ji}<-1}}m_{ji}a_i \quad \forall\; j\in \hbfI.
\end{equation*}
In particular, $\mu_j=6d_j^3 a_j$ unless there is $i_j\in \hbfI$ such that $a_{j,i_j}<-1$. If such an index $i_j$ exists then, provided $\gkm$ is not of type $\mathsf{C}_2^{(1)}$, it is unique and one has 
\begin{equation*}
\mu_j=
\begin{cases}
6d_j^3(a_j-a_{i_j}) \quad &\text{ if }\; a_{j,i_j}=-2,\\
6d_j^3(a_j-4a_{i_j}) \quad &\text{ if }\; a_{j,i_j}=-3.
\end{cases}
\end{equation*}
\end{lem}
\begin{pf}
Fix $j\in \hbfI$. Then, since $\mu_j=\sum_{i\in \hbfI} a_i (d_ia_{ij})^3$, we have 
\begin{equation*}
\mu_j/d_j^3=\sum_{i\in \hbfI} a_i a_{ji}^3=\sum_{i\in \hbfI} a_i a_{ji}+\sum_{i\in \hbfI}a_i (a_{ji}^3-a_{ji})=\sum_{i\in \hbfI}a_i (a_{ji}^3-a_{ji}),
\end{equation*}
where we have used that $\hbfA \underline{\delta}=0$, with $\underline{\delta}\ceq  (a_j)_{j\in \hbfI}\in \Z^{\hbfI}$. Since $a_{jj}^3-a_{jj}=6$ and $m_{ji}=a_{ji}-a_{ji}^3$ is zero when $a_{ji}=-1$, this proves the first part of the lemma.

 As for the second assertion, if $i\in \hbfI$ is such that $a_{ji}<-1$, then the vertices $i$ and $j$ of the Dynkin diagram of $\gkm$ are connected by multiple edges, with an arrow going from $j$ to $i$. By \cite[\S4.8, Tables Aff 1-2-3]{kac}, any such index $i=i_j$ is unique except in type $\mathsf{C}_2^{(1)}$, where there is a single short simple root connected to two long simple roots.  The claimed formulas now follow from the fact that, for any such type we have 
 \begin{equation*}
 \sum_{\substack{i\in \hbfI\\a_{ji}<-1}}m_{ji}a_i
=
 \begin{cases}
 6a_i &\text{ if } a_{ji}=-2\\
 24a_i&\text{ if } a_{ji}=-3.
 \end{cases} \qedhere
 \end{equation*}
\end{pf}

\subsection{The rank of $(\bfB\,|\, \underline{\mu})$}\label{appB-full-rank}

We now come to the following result, which is a restatement of Proposition \ref{P:full-rank}.
\begin{prop}\label{P:app-full-rank}
The augmented matrix $(\mathbf{B}\,|\, \underline{\mu})$ has rank $|\hbfI|$. 
\end{prop}
\begin{pf}
Let $\bfB_j$ denote the $j$-th column of $\bfB$. We first note that there does not exist a vector $\underline{\gamma}\in \Q^{\hbfI}$ such that 
$
\bfB\underline{\gamma}\in (\Q_{>0})^{\hbfI}.
$
Indeed, if there did, then $\underline{\delta}=(a_j)_{j\in \hbfI}\in (\Z_{>0})^{\hbfI}$ would satisfy $\underline{\delta}^t (\bfB\underline{\gamma})>0$, but this contradicts that 
$
\underline{\delta}^t(\bfB\underline{\gamma})=(\underline{\delta}^t \bfB)\underline{\gamma}=0.
$ 

Hence, to prove the proposition it is sufficient to show that there is a vector $\underline{\gamma}\in \Z^{\hbfI}$ such that 
\begin{equation}\label{ul-gamma}
\underline{\gamma}\in (\Z\underline{\mu} + \mathrm{range}(\bfB)) \cap (\Z_{>0})^{\hbfI}. 
\end{equation}
If the Dynkin diagram associated to $\hbfA$ is simply laced, then Lemma \ref{L:full-rank} yields $\mu_j=6a_j>0$ for all $j\in \hbfI$, so we may take $\underline{\gamma}=\underline{\mu}$. 

Suppose instead that the Dynkin diagram of $\hbfA$ belongs to the following list: 
\begin{equation*}
\mathsf{C}_\ell^{(1)}\, (\ell>2), \; \mathsf{F}_4^{(1)},\; \mathsf{A}_{2\ell-1}^{(2)},\; \mathsf{E}_6^{(2)}.
\end{equation*}
Then it follows from Lemma \ref{L:full-rank} and \cite[Tables Aff 1-2]{kac} that $\mu_j=6d_j^3 a_j$ for all $j\in \hbfI$ such that $a_{ji}\geqslant -1$ for all $i\in \hbfI$, and $\mu_j=6d_j^3(a_j-a_{i})=6d_j^3$ if $j\in \hbfI$ is such $a_{ji}=-2$ for some $i\in \hbfI$. Thus, for all these types we may again take $\underline{\gamma}=\underline{\mu}$.  

To complete the proof, it remains to show that there is $\underline{\gamma}$ as in \eqref{ul-gamma} associated to each of the following affine types:
\begin{equation*}
\mathsf{C}_2^{(1)}, \;  \mathsf{G}_2^{(1)}, \; \mathsf{D}_4^{(3)},\; \mathsf{B}_\ell^{(1)},\; \mathsf{A}_{2\ell}^{(2)},\; \mathsf{D}_{\ell+1}^{(2)}.
\end{equation*}
We do this case by case, using the classification provided by \cite[Thm.~4.8]{kac} and appealing to the notation from Tables Fin and Aff of \cite[\S4.8]{kac} as needed. 

\noindent \textit{The $\mathsf{B}_\ell^{(1)}$ case}. In this case, the Dynkin diagram of $\gkm$ is 
\begin{equation*}
\dynkin [extended,labels={1,1,2,2,2,2,2},scale=1.5]B[1]{}
\end{equation*}
and we take $\hbfI=\{0,l,\ldots, \ell\}$ with $\{1,\ldots,\ell\}$ labelling the type $\mathsf{B}_\ell$ subdiagram with vertices given by the blackened nodes (labelled from left to right in increasing order), and $i=0$ labels the extending node. In the above diagram, the node $j\in \hbfI$ is marked by the number $a_j$.  By Lemma \ref{L:full-rank} we have $\mu_j=6d_j^3a_j=48 a_j$ for all $0\leqslant j<\ell$ and $\mu_\ell=0$. Equivalently, 
\begin{equation*}
\underline{\mu}^t=48\cdot(1,1,2,\ldots, 2,0).
\end{equation*}
On the other hand, we have $\bfB_{\ell}^t=(0,\ldots,-2,2)$, so we may take 
\begin{equation*}
\underline{\gamma}\ceq  \underline{\mu}+\bfB_{\ell}\in (\Z_{>0})^{\hbfI}. 
\end{equation*}

\noindent \textit{The $\mathsf{D}_{\ell+1}^{(2)}$ case}. 
In this case the Dynkin diagram is
\begin{equation*}
	\dynkin [extended, labels={1,1,1,1,1,1,1,1},scale=1.5]D[2]{}
\end{equation*}
As before, we take $\hbfI=\{0,1,\ldots,\ell\}$, and label the nodes from left to right in increasing order, marking the node $j\in \hbfI$ by the number $a_j$. 
Since $a_j=1$ for all $j\in \hbfI$, Lemma \ref{L:full-rank} yields
\begin{equation*}
\underline{\mu}^t=(0,6d_1^3,\ldots, 6d_{\ell-1}^3,0)=48\cdot (0,1,\ldots,1,0).
\end{equation*}
As $\mathbf{B}_0^t=(2,-2,0,\ldots,0)$ and $\mathbf{B}_\ell^t=(0,\ldots,0,-2,2)$, we deduce that 
\begin{equation*}
\underline{\gamma}\ceq  \underline{\mu}+\mathbf{B}_0+\mathbf{B}_\ell\in (\Z_{>0})^{\hbfI}. 
\end{equation*}
\noindent \textit{The $\mathsf{A}_{2\ell}^{(2)}$ case}. In this case, the Dynkin diagram is 
\begin{equation*}
\dynkin [extended, labels={2,2,2,2,2,2,1}, scale=1.5]A[2]{even}
\end{equation*}
and we choose the same labeling convention as in the $\mathsf{D}_{\ell+1}^{(2)}$ case. We have $d_0=1$, $d_j=2$ for $0\leqslant j<\ell$, and $d_\ell=4$.  Lemma \ref{L:full-rank} thus yields 
\begin{equation*}
\underline{\mu}^t=(0,12d_1^3,\ldots,12d_{\ell-2}^3,6d_{\ell-1}^3(a_{\ell-1}-a_\ell), 6 d_\ell^3)=48\cdot (0, 2, \ldots,2, 1, 8).
\end{equation*}
Since $\bfB_0^t=(2,-2,0,\ldots,0)$, the vector $\underline{\gamma}\ceq  \underline{\mu}+\bfB_0$ lies in $(\Z_{>0})^{\hbfI}$. 

\noindent\textit{The cases $\mathsf{C}_2^{(1)}$, $\mathsf{G}_2^{(1)}$ and  $\mathsf{D}_4^{(3)}$.}
In these cases, all the relevant data, together with a choice of $\underline{\gamma}$ satisfying \eqref{ul-gamma}, is given in the following table:

\begin{center}
\begin{tabular}{|c|c|c|c|c|}
\hline 
Type & Diagram &  $\bfB$ &  $\underline{\mu}$ & $\underline{\gamma}$\\[3pt]
\hline \hline 
& & & & \\[-.5em]
$\mathsf{C}_2^{(1)}$ & $\dynkin [extended, labels={1,2,1}]C[1]{2} $ & 
\scalebox{0.8}{$\begin{pmatrix} 4 & -2 & 0\\-2 & 2 & -2\\ 0 & -2 & 4\end{pmatrix}$}
& \scalebox{0.8}{$\begin{bmatrix} 48 \\0\\48\end{bmatrix}$} &$\underline{\mu}+\mathbf{B}_2$ \\
& & & & \\[-.5em]
 \hline
& & & & \\[-.5em]
$\mathsf{G}_2^{(1)}$ & $\dynkin [extended, labels={1,2,3}]G[1]{2}$ & \scalebox{0.8}{$\begin{pmatrix} 6 & -3 & 0\\-3 & 6 & -3 \\ 0 & -3 & 2\end{pmatrix}$}  & \scalebox{0.8}{$\begin{bmatrix} 162\\324\\-30\end{bmatrix}$} & $\underline{\mu}+16\mathbf{B}_3$\\
& & & & \\[-.5em]
 \hline 
& & & & \\[-.5em]
$\mathsf{D}_4^{(3)}$ & $\dynkin [extended,labels={1,2,1}]D[3]{4}$  & \scalebox{0.8}{$\begin{pmatrix} 2 & -1 & 0\\-1 & 2 & -3 \\ 0 & -3 & 6\end{pmatrix}$} & 
\scalebox{0.8}{$\begin{bmatrix} 6\\-12\\162\end{bmatrix}$} & $\underline{\mu}-5\mathbf{B}_3$\\[-.5em]
& & & & \\ \hline
\end{tabular}
\end{center}
Therefore, we may conclude that  $(\mathbf{B}\,|\, \underline{\mu})$ has rank $|\hbfI|$ for all affine types, where we continue to exclude types $\mathsf{A}_1^{(1)}$ and $\mathsf{A}_2^{(2)}$.
\qedhere
\end{pf}

\section{Determinant of affine quantum Cartan matrices}\label{app:QCM}

In this short appendix, we prove Lemma~\ref{lem:order}. 
We provide explicit formulae for both the determinant of the
symmetrized affine quantum Cartan matrix $\bfB(\sop)$ and 
\[
\qdzero= \left.
\frac{\det(\bfB(\sop))}{(\sop-\sop^{-1})^2}
\right|_{\sop=1}\,.
\]
We exclude type $\sfA_1^{(1)}$.
For the exceptional types, the formulae are obtained by a direct computation
while, for the other types,  it is enough to proceed by a corank two Lagrange expansion.
The resulting formulae are given in Table~\ref{table-ADE} for untwisted simply laced types;
in Table~\ref{table-BCFG} for untwisted non simply laced types; 
in Table~\ref{table-twisted} for twisted types.

By direct inspection, this shows that $\det(\bfB(\sop))$ has order of vanishing $2$ at $\sop=1$
and all its zeros lie in $U(1)$.

\begin{rem}
For simply laced types , the explicit formulae for $\det(\bfB(\sop))$
first appeared in \cite[Sec.~4]{suter-07}. Note however that in 
untwisted type BCFG our formulae differ from those given in \emph{loc. cit.},
due to a different definition of the quantum Cartan matrix in these types.
\end{rem}

\begin{table}[h]
	\begin{tabular}{|c|c|c|}
		\hline
		Type of $\gkm$     & $\det(\bfB(\sop))$ & $\qdzero$ \\
		\hline\hline
		$\mathsf{A}_{\ell}^{(1)}$ & $\sop^{-\ell-1}(\sop^{\ell+1}-1)^2$ & $\frac{(\ell+1)^2}{4}$\\ 
		\hline
		$\mathsf{D}_\ell^{(1)}$ & $\sop^{-\ell-1}(\sop^2+1)(\sop^{2(\ell-2)}-1)(\sop^4-1)$ & $4(\ell-2)$\\
		\hline
		$\mathsf{E}_6^{(1)}$ & $\sop^{-7}(\sop^2+1)(\sop^{6}-1)^2$ & 18\\
		\hline
		$\mathsf{E}_7^{(1)}$ & $\sop^{-8}(\sop^2+1)(\sop^{8}-1)(\sop^6-1)$ & 24\\
		\hline
		$\mathsf{E}_8^{(1)}$ & $\sop^{-9}(\sop^2+1)(\sop^{10}-1)(\sop^6-1)$ & 30 \\
		\hline
	\end{tabular}
	\vspace{0.25cm}
	\caption{Untwisted type ADE}\label{table-ADE}
\end{table}	

\begin{table}[h]
	\begin{tabular}{|c|c|c|}
		\hline
		Type of $\gkm$     & $\det(\bfB(\sop))$ & $\qdzero$\\
		\hline\hline
		$\mathsf{B}_\ell^{(1)}$ & $\sop^{-3\ell-1}(\sop^2+1)^{\ell}(\sop^{2(2\ell-3)}-1)(\sop^8-1)$ & $2^{\ell+2}(2\ell-3)$\\
		\hline
		$\mathsf{C}_\ell^{(1)}$ & $\sop^{-3\ell-11}(\sop^2+1)^{\ell+1}(\sop^4+1)(\sop^{4(\ell+2)}-1)(\sop^8-1)$ & $2^{\ell+5}(\ell+2)$\\
		\hline
		$\mathsf{F}_4^{(1)}$ & $\sop^{-11}(\sop^2+1)^3(\sop^{10}-1)(\sop^6-1)$ & $120$\\
		\hline
		$\mathsf{G}_2^{(1)}$ & $[3]_\sop\sop^{-9}(\sop^2+1)(\sop^{10}-1)(\sop^6-1)$ & $90$\\
		\hline
	\end{tabular}
	\vspace{0.25cm}
	\caption{Untwisted type BCFG}\label{table-BCFG}
\end{table}

\begin{table}[h]
	\begin{tabular}{|c|c|c|}
		\hline
		Type of $\gkm$     & $\det(\bfB(\sop))$  & $\qdzero$\\
		\hline\hline
		$\mathsf{A}_{2\ell}^{(2)}$ & $\sop^{-6\ell-5}(\sop^2+1)^{2\ell}(\sop^{2(4\ell+1)}-1)(\sop^8-1)$  & $4^{\ell+1}(4\ell+1)$\\
		\hline
		$\mathsf{A}_{2\ell-1}^{(2)}$  & $\sop^{-2(\ell+1)}(\sop^2+1)(\sop^{2(\ell-2)}-1)(\sop^4-1)$& $4(\ell-2)$ \\
		\hline
		$\mathsf{D}_{\ell+1}^{(2)}$& $\sop^{-3\ell-2}(\sop^2+1)^\ell(\sop^{4\ell}-1)(\sop^4-1)$ & $2^{\ell+2}\ell$\\
		\hline
		$\mathsf{E}_6^{(2)}$& $\sop^{-9}(\sop^2+1)^2(\sop^{8}-1)(\sop^6-1)$ & $48$\\
		\hline
		$\mathsf{D}_4^{(3)}$ & $\sop^{-7}(\sop^2+1)(\sop^{6}-1)^2$ & $18$\\
		\hline
	\end{tabular}
	\vspace{0.25cm}
	\caption{Twisted type}\label{table-twisted}
\end{table}


\noindent {\bf Acknowledgments.} This project started during a
\emph{Research in Pairs} visit at the Centro Internazionale per la
Ricerca Matematica in Trento, May 21 - June 9, 2023.
We are grateful to the entire CIRM team for providing us an
exceptional research environment. We would like to thank Mar\'{i}a Ang\'{e}lica Cueto
for her help with the more computational aspects of the paper, Ovidiu Costin
and Valerio Toledano Laredo for their insightful comments.

AA was supported in part by the INdAM Project 2024 (CUP E53C23001740001) and the PRIN grant 2022HMTBLL. SG was supported through the Simons foundation
collaboration grant 526947. CW gratefully acknowledges the support of the Natural Sciences
and Engineering Research Council of Canada (NSERC), provided via the Discovery Grants Program
(Grant RGPIN-2022-03298 and DGECR-2022-00440).


 \newcommand{\noop}[1]{}
\providecommand{\bysame}{\leavevmode\hbox to3em{\hrulefill}\thinspace}
\providecommand{\MR}{\relax\ifhmode\unskip\space\fi MR }
\providecommand{\MRhref}[2]{%
	\href{http://www.ams.org/mathscinet-getitem?mr=#1}{#2}
}
\providecommand{\href}[2]{#2}

\renewcommand{\MR}[1]{}

\end{document}